\documentclass[a4paper,11pt,draft]{amsart}

\usepackage{amsmath}
\usepackage{amsthm}
\usepackage{amsfonts}
\usepackage{amssymb}
\usepackage[all]{xy}

\newcommand{\myAA}{\mathcal{A}}
\newcommand{\myFF}{\mathcal{F}}
\newcommand{\FF}{\mathcal{F}}
\newcommand{\myCC}{\mathcal{C}}
\newcommand{\ZZ}{\mathbb{Z}}
\newcommand{\PP}{\mathbb{P}}
\newcommand{\QQ}{\mathbb{Q}}
\newcommand{\RR}{\mathbb{R}}
\newcommand{\BB}{\mathbf{B}}
\newcommand{\yy}{\mathrm{y}}
\newcommand{\bfy}{\mathbf{y}}

\newtheorem{theorem}{Theorem}[section]
\newtheorem{proposition}{Proposition}[section]
\newtheorem{corollary}{Corollary}[section]
\newtheorem{definition}{Definition}[section]
\newtheorem{lemma}{Lemma}[section]
\newtheorem{remark}{Remark}[section]
\newtheorem{example}{Example}[section]

\begin{document}
\title[Cluster algebras of type $A_{2}^{(1)}$]{Cluster algebras of type $A_2^{(1)}$}

\author{Giovanni Cerulli Irelli}
\address{Universit\`a degli studi di Padova, Dipartimento di Matematica Pura ed Applicata. Via Trieste 63, 35121, Padova (ITALY)}
\email{gcerulli@math.unipd.it}

\begin{abstract}
In this paper we study cluster algebras $\myAA$ of type $A_2^{(1)}$. We solve the recurrence relations among the cluster variables (which form a T--system of type $A_2^{(1)}$).  We solve the recurrence relations among the coefficients of $\myAA$ (which form a Y--system of type $A_2^{(1)}$). In $\myAA$ there is a natural notion of positivity.  We find linear bases $\BB$ of $\myAA$ such that positive linear combinations of elements of $\BB$ coincide with the cone of positive elements. We call these bases \emph{atomic bases} of $\myAA$.  These are the analogue of the ``canonical bases'' found by Sherman and Zelevinsky in type $A_{1}^{(1)}$. Every atomic basis consists of cluster monomials together with extra elements. We provide explicit expressions for the elements of such bases in every cluster. We prove that the elements of $\BB$ are parameterized by $\ZZ^3$ via their $\mathbf{g}$--vectors in every cluster. We prove that the denominator vector map in every acyclic seed of $\myAA$ restricts to a bijection between $\BB$ and $\ZZ^3$. In particular this gives an explicit algorithm to determine the ``virtual'' canonical decomposition of every element of the root lattice of type $A_2^{(1)}$. We find explicit recurrence relations to express every element of $\myAA$ as linear combinations of elements of $\BB$.
\end{abstract}

\maketitle
\tableofcontents

\section{Introduction}\label{sec:Intro}

Cluster algebras were introduced in \cite{FZI} by Fomin and Zelevinsky in order to define a combinatorial framework for studying canonical bases in quantum groups. In this perspective, it is an important problem, still widely open, to construct explicitly linear bases of cluster algebras.  Among  possible bases of a cluster algebra $\myAA$  (see \cite{CK2}, \cite{CZ}, \cite{dupont-2009}, \cite{DupontGeneric}, \cite{GLS10}), there may exists a specific one of particular interest, that we call  \emph{atomic basis}. It is not known in general if such a basis exists. For cluster algebras of type $A_2$ and $A_1^{(1)}$, Sherman and Zelevinsky were able to construct atomic bases  in \cite{Sherman} but these are the only known examples at this time. 

In this paper we  prove the existence of  atomic bases   $\BB$ of every cluster algebra of type $A_2^{(1)}$. We provide a completely explicit description of the elements of $\BB$ which consists of cluster monomials together with extra elements. We find explicit straightening relations which provide a recursive way to compute the linear expansion in $\BB$ of every element of $\myAA$.  Our proof of the existence of an atomic basis $\BB$ of $\myAA$ follows \cite{Sherman} but with an important technical difference:  we prove that the powerful $\mathbf{g}$--vector parametrization of cluster monomials found in \cite{FZIV} extends  to the elements of $\BB$. In particular we work directly in the cluster algebra with coefficients in an arbitrary semifield. This allows us to provide explicit formulas for the elements of $\BB$. Moreover we think that this approach can help in the study of the existence of atomic bases of cluster algebras of higher rank in view of results of \cite{FuKeller} and \cite{DWZII} and of the geometric realization of $\BB$ given in \cite{CEsp} and \cite{dupont-2009}.  

The paper is organized as follows: in section~\ref{sec:MainResults} we collect  the main results  and the remaining sections are devoted to proofs. 
 
\section{Main results}\label{sec:MainResults}
 
A \emph{semifield} $\PP=(\PP,\cdot,\oplus)$ is an abelian multiplicative group $(\PP,\cdot)$ endowed with an auxiliary addition $\oplus:\PP\times\PP\rightarrow\PP$ which is associative, commutative and $a(b \oplus c)=ab\oplus ac$ for every $a,b,c\in\PP$.  

The main example of a semifield is a \emph{tropical semifield}: by definition a tropical semifield $Trop(y_{j}:j\in J)$ is an abelian multiplicative group freely generated by the elements $\{y_{j}:j\in J\}$ (for some set of indexes $J$) endowed with the auxiliary addition $\oplus$ given by:
$$
\prod_{j} y_{j}^{a_{j}}\oplus\prod_{j} y_{j}^{b_{j}}:=\prod_{j} y_{j}^{\text{min}(a_{j},b_{j})}.
$$ 
It can be shown (see \cite[Section~5]{FZI}) that every semifield $\PP$ is torsion-free as a multiplicative group and hence its group ring $\ZZ\PP$ is a domain.  Given a semifield $\PP$, let $\QQ\PP$ be the field of fractions of $\ZZ\PP$ and $\FF_{\PP}=\QQ\PP(x_{1},x_{2},x_{3})$ be the field of rational functions in three independent variables with coefficients in $\QQ\PP$. We consider the $\ZZ\PP$--subalgebra $\myAA_\PP$ of $\myFF\PP$ defined as follows: we choose arbitrarily three elements $y_1$, $y_2$ and $y_3$ of $\PP$. We consider the family $\{y_{1;m}:m\in\ZZ\}\subset\PP$ of coefficients defined by the initial conditions:
$$
\begin{array}{ccc}
y_{1;0}=\frac{1}{y_{3}},&y_{1;1}=y_{1},&y_{1;2}=\frac{y_{1}y_{2}}{y_{1}\oplus1}
\end{array}
$$
together with the recursive relations:
\begin{equation}\label{Eq:Y1mY1m+3}
y_{1;m}y_{1;m+3}=\frac{y_{1;m+2}y_{1;m+1}}{(y_{1;m+1}\oplus1)(y_{1;m+2}\oplus 1)}
\end{equation}
We consider the following rational functions $w$, $z$  of $\myFF_\PP$: 
\begin{equation}\label{Eq:DefinitionWZ}
\begin{array}{cc}
 w:=\frac{y_2x_1+x_3}{(y_{2}\oplus1)x_2},&z:=\frac{y_1y_3x_1x_2+y_1+x_2x_3}{(y_{1}y_{3}\oplus y_{1}\oplus1)x_1x_3}
\end{array}
\end{equation}
We consider the elements $x_{m}$, $m\in\ZZ$, of $\myFF_\PP$ defined recursively by:
\begin{equation}\label{ExchRelIntro}
x_{m}x_{m+3}=\frac{x_{m+1}x_{m+2}+y_{1;m}}{y_{1;m}\oplus1}.
\end{equation}
By definition $\myAA_{\PP}$ is the $\ZZ\PP$--subalgebra of $\FF_{\PP}$ generated by the rational functions $w$, $z$ and  $x_{m}$, $m\in\ZZ$. 

In section~\ref{Sec:AlgStructure} we prove that $\myAA_\PP$ is the cluster algebra with initial seed (see sections~\ref{Sec:Background} for some background on cluster algebras)
\begin{equation}\label{Eq:InitialSeed}
\Sigma:=\{H=\left(\def\objectstyle{\scriptstyle}\def\labelstyle{\scriptstyle}\vcenter{\xymatrix@R=0pt@C=0pt{0&1&1\\-1&0&1\\-1&-1&0}}\right), \displaystyle{\mathbf{x}=\{x_1,x_2,x_3\}, \mathbf{y}=\{y_1,y_2,y_3\}}\}.
\end{equation}
The rational functions $w$, $z$ and $x_m$, $m\in\ZZ$, are the \emph{cluster variables} of $\myAA_\PP$. To the matrix $H$ is naturally associated the quiver
\begin{equation}\label{Eq:Q}
\xymatrix @R=5pt@C=10pt{
         &&2\ar[dl]&\\
Q:=&1&  &3\ar[ul]\ar[ll]}
\end{equation}
whose underlying graph is an extended Dynkin graph of type $A_2^{(1)}$. The algebra $\myAA_\PP$ is hence called a cluster algebra of type $A_2^{(1)}$ (see \cite{FZII}). For $\PP=\{1\}$ such cluster algebra appears in \cite[example~7.8]{FZI}.
 
The family of coefficients $\{y_{1;m}\}$ form a $Y$--system of type $A_2^{(1)}$ (see \cite{Ysyst}, \cite{NakanishiKeller}). We solve the recursion~\eqref{Eq:Y1mY1m+3} in \eqref{Eq:Y1mGeneral}. The family of cluster variables form a $T$--system of type $A_2^{(1)}$ (see \cite{NakanishiKeller}). We solve the recurrence \eqref{ExchRelIntro} in theorem~\ref{Thm:ExplicitExpressions}

Every cluster variable $s_{1}$ can be completed to a  set $\myCC=\{s_{1},s_{2},s_{3}\}$ of cluster variables which form a free generating set of the field $\FF_{\PP}$, so that $\FF_{\PP}\simeq\QQ\PP(s_{1},s_{2},s_{3})$. Such a set $\myCC$ is called a \emph{cluster} of $\myAA_{\PP}$.  The clusters
of $\myAA_{\PP}$ are the sets
$\{x_m,x_{m+1},x_{m+2}\}$, $\{x_{2m+1},w,x_{2m+3}\}$ and the set $\{x_{2m},z,x_{2m+2}\}$ for every $m\in\ZZ$. Figure~\ref{fig:ExchangeGraph}  shows (a piece of) the  ``exchange graph" of $\myAA_{\PP}$. By definition it has clusters as
\begin{figure}[htbp]
\begin{center}
\xymatrix{
*{\bullet}\ar@{-}[d]+0\ar@{-}[rrrrrrrrrrrr]-0^{w}\ar@{}|{\cdots}[rrd]&*{}          &*{\bullet}\ar@{-}[d]\ar@{}|{x_{-1}}[rrd]&*{}          &*{\bullet}\ar@{-}[d]\ar@{}|{x_{1}}[rrd]& *{}         &*{\bullet}\ar@{-}[d]\ar@{}|{x_{3}}[rrd]&*{}         &*{\bullet}\ar@{-}[d]\ar@{}|{x_{5}}[rrd]&  *{}       &*{\bullet}\ar@{-}[d]\ar@{}|{\cdots}[rrd]&*{}         &*{\bullet}\ar@{-}[d]&*{}           \\
       *{\bullet}\ar@{-}[rrrrrrrrrrrr]   &*{\bullet}\ar@{-}[d]\ar@{}|{\cdots}[rrd]&*{\bullet}          &*{\bullet}\ar@{-}[d]\ar@{}|{x_{0}}[rrd]&*{\bullet}          &*{\bullet}\ar@{-}[d]\ar@{}|{x_{2}}[rrd]&*{\bullet}          &*{\bullet}\ar@{-}[d]\ar@{}|{x_{4}}[rrd]&*{\bullet}          &*{\bullet}\ar@{-}[d]\ar@{}|{\cdots}[rrd]&*{\bullet}          &*{\bullet}\ar@{-}[d]&*{\bullet}          \\
  *{}        &*{\bullet}\ar@{-}[rrrrrrrrrr]_{z}&*{}          &*{\bullet}&*{}          &*{\bullet}&*{}          &*{\bullet}&*{}          &*{\bullet}&          &*{\bullet}&*{}
}
\caption{The exchange graph of $\myAA_{\PP}$ }\label{fig:ExchangeGraph}
\end{center}
\end{figure}
vertexes and an edge between two clusters if they share exactly two cluster variables.
In this figure cluster variables are associated with regions: there are infinitely many
bounded regions labeled by the $x_m$'s, and there are two unbounded regions labeled
respectively by $w$ and $z$. Each cluster $\{s_{1},s_{2},s_{3}\}$ corresponds to the common vertex of the three
regions $s_{1}$, $s_{2}$ and $s_{3}$. 

We have already observed that $\FF_{\PP}\simeq\QQ\PP(s_{1},s_{2},s_{3})$ for every cluster $\{s_{1},s_{2},s_{3}\}$ of $\myAA_{\PP}$ and hence every element of $\myAA_{\PP}$ can be expressed as a rational function in $\{s_{1},s_{2},s_{3}\}$. By the \emph{Laurent phenomenon} proved in \cite{FZI}, such a rational function is actually a Laurent polynomial. 
Following \cite{Sherman} we say that an element  of $\myAA_{\PP}$ is \emph{positive} if its Laurent expansion in every cluster has coefficients in $\ZZ_{\geq0}\PP$. Positive elements form a semiring, i.e. sums and products of positive elements are positive. We say that a $\ZZ\PP$--basis $\BB$ of $\myAA_{\PP}$ is \emph{atomic} if the semiring of positive elements coincides with the $\ZZ_{\geq0}\PP$--linear combinations of elements of it. If an atomic basis   exists, it is composed by positive indecomposable elements, i.e. positive elements that cannot be written as a sum of positive elements. Moreover such a basis is unique up to rescaling by elements of $\PP$. 
\begin{definition}\label{def:Un}
Let $W:=(y_{2}\oplus1)w$ and $Z:=(y_{1}y_{3}\oplus y_{1}\oplus1)z$.
We define elements $\{u_n|\,n\geq 0\}$ of $\myAA_{\PP}$ by the
initial conditions:
\begin{equation}\label{Eq:DefinitionU0u_1U2}
\begin{array}{ccc}
u_0=1,&u_1=ZW-y_1y_3-y_2,&u_2=u_1^2-2y_1y_2y_3
\end{array}
\end{equation}                                                         
together with the recurrence relation
\begin{equation}\label{Eq:DefinitionUn+1}
\begin{array}{cc}
u_{n+1}=u_1u_n-y_1y_2y_3u_{n-1},&n\geq2.
\end{array}
\end{equation} 
\end{definition}

\begin{theorem}\label{Thm:CanonicalBasis}
The set 
\begin{equation}\label{Def:B}
\BB:=\{\text{cluster monomials}\}\cup\{u_nw^k,u_nz^k|\,n\geqslant 1,\,k\geqslant 0\}
\end{equation}
is an atomic basis   of $\myAA_{\PP}$. It is unique up to rescaling by elements of $\PP$.
\end{theorem}
Theorem~\ref{Thm:CanonicalBasis} shows that the atomic bases of $\myAA_\PP$ consist of cluster monomials together with some extra elements. This is precisely the form of other linear basis constructed for cluster algebras of affine type: the generic basis introduced by G.~Dupont (see \cite{DupontGeneric}, \cite{GLS10}, \cite{DingXiaoXu}); the semicanonical basis introduced by Caldero and Zelevinsky (\cite{CZ}) and the basis obtained by Dupont by using transverse quiver Grassmannians (\cite{dupont-2009}). 

The proof of theorem~\ref{Thm:CanonicalBasis} is given in section~\ref{Sec:CanonicalBasis}. We now collect some properties of the basis $\BB$.

\subsection{Parametrization of $\BB$ by denominator vectors}

By the already mentioned Laurent phenomenon \cite[Theorem 3.1]{FZI} every element $b$ of $\myAA_{\PP}$ is a Laurent polynomial in $\{x_1,x_2,x_3\}$ of the form $\frac{N_{b} (x_1 , x_2 , x_3 )}{x_1^{d_1}x_2^{d_2}x_3^{d_3}}$ for some primitive, i.e. not divisible by any $x_i$, polynomial $N_b\in\ZZ\PP[x_1,x_2,x_3]$ in $x_1$,
$x_2$ and $x_3$, and some integers $d_1$, $d_2$, $d_3$. We consider the root lattice $Q$ generated by an affine root system of type $A_2^{(1)}$. The choice of a simple system $\{\alpha_{1},\alpha_{2},\alpha_{3}\}$, with coordinates $\{e_{1},e_{2},e_{3}\}$, identifies $Q$ with $\ZZ^3$. We usually write elements of $Q$ as column vectors with integer coefficients and we denote them by a bold type letter.  The map $b\mapsto\mathbf{d}(b)=(d_1,d_2,d_3)^t$ is hence a map between $\myAA_{\PP}$ and $Q$; it is called the
\emph{denominator vector map in the cluster $\{x_1,x_2,x_3\}$}. The following result provides a parameterization of  $\BB$  by $Q$ via the denominator vector map. Recall  that given $\delta:= (1, 1, 1)^{t}$, the minimal                             
positive imaginary root, and $\Pi^{\circ}= \{\alpha_{1},\alpha_{3}\}$, a basis of simple roots for a root
system $\Delta^{\circ}$ of type $A_2$, the positive real roots of $Q$ are of the form $\alpha+n\delta$ with $n\geq0$ if $\alpha$ is a positive root of $\Delta^{\circ}$ and $n\geq 1$ if $\alpha$ is a negative root of $\Delta^{\circ}$ (see e.g. \cite[Proposition~6.3]{KacBook}).

\begin{theorem}\label{Prop:BijectionDenominators}
The denominator vector map $\mathbf{d}:\myAA_\PP\rightarrow Q$ $:b\mapsto\mathbf{d}(b)$ in the cluster $\{x_1,x_2,x_3\}$
restricts to a bijection between  $\BB$  and $Q$. Under this bijection positive real roots of the
root system of type $A_2^{(1)}$ correspond to the set of cluster variables different from $x_1$, $x_2$ and $x_3$ together with the
set $\{u_{n} w, u_{n}$z$|\, n\geq1\}$. Moreover for every cluster $\myCC =\{c_{1} , c_{2} , c_{3}\}$, the set
$\{\mathbf{d}(c_{1} ), \mathbf{d}(c_{2} ), \mathbf{d}(c_{3} )\}$ is a $\ZZ$--basis of $Q$.
\end{theorem}
The denominator vector map in the cluster $\{x_1,w,x_3\}$ does not restrict to a bijection between $\BB$ and the whole lattice $\ZZ^3$ (see remark~\ref{Rem:NoBijectionDen}).
\begin{figure}[htbp]
\xymatrix@R=90pt@C=70pt@L=5pt{            &           
&*{\bullet}\ar@{-}[rrdd]|(.5)*{\bullet}="X-1"\ar@{-}[dd]_<{\mathbf{d}(x_{0})}
|(.15)*{\bullet}="X-2"_(.15){\mathbf{d}(x_{-2})}
|(.27)*{\bullet}="X-4"_(.27){\mathbf{d}(x_{-4})}
|(.37)*{\bullet}="X-6"_(.37){\mathbf{d}(x_{-6})}
|(.45)*{\bullet}="X-8"
|(.50)*{\bullet}="X-10"
|(.53)*{\bullet}
|(.55)*{\bullet}
_(.55){\mathbf{d}(x_{-2m})}
|(1)*{\bullet}="X5"^(1.02){\mathbf{d}(x_{5})}
|(.96)*{\bullet}="X7"^(.96){\mathbf{d}(x_{7})}
|(.90)*{\bullet}="X9"^(.90){\mathbf{d}(x_{9})}
|(.83)*{\bullet}="X11"
|(.80)*{\bullet}="X13"
|(.78)*{\bullet}="X15"
|(.76)*{\bullet}
^(.80){\mathbf{d}(x_{2m+1})}
&             &\\
&*{\bullet}\ar@{..};"X-2"\ar@{..};"X-4"\ar@{..};"X-6"\ar@{..};"X-8"
&\ar@{}[u]+<0mm,5mm>*{(0,0,1)}                              &*{\bullet}\ar@{-}[dlll]^(-.03){\mathbf{d}(x_{-1})}
|(.06)*{\bullet}="X-3"^(.06){\mathbf{d}(x_{-3})}
|(.12)*{\bullet}="X-5"^(.14){\mathbf{d}(x_{-5})}
|(.17)*{\bullet}="X-7"
|(.2)*{\bullet}="X-9"
|(.22)*{\bullet}="X-11"
|(.24)*{\bullet}
|(.26)*{\bullet}
^(.23){\mathbf{d}(x_{-(2m+1)})}
_(1){\mathbf{d}(x_{4})}
|(.80)*{\bullet}="X6"_(.80){\mathbf{d}(x_{6})}
|(.65)*{\bullet}="X8"_(.65){\mathbf{d}(x_{8})}
|(.55)*{\bullet}="X10"_(.55){\mathbf{d}(x_{10})}
|(.48)*{\bullet}="X12"
|(.46)*{\bullet}="X14"
|(.44)*{\bullet}
|(.42)*{\bullet}
|(.40)*{\bullet}
_(.44){\mathbf{d}(x_{2m})}
&\\
*{\bullet}\ar@{-}[rruu]|(.5)*{\bullet}="Z"\ar@{-}[rrrr]\ar@{}[rrrr]+<0mm,-5mm>*{(0,1,0)}\ar@{}[rrrr]+<3mm,3mm>*{\scriptstyle{\mathbf{d}(w)}}&              &             *{\bullet}                   &           &*{\bullet}\ar@{..};"X-3"\ar@{..};"X-5"\ar@{..};"X-7"\ar@{..};"X-9"
\ar@{..};"X7"\ar@{..};"X9"\ar@{..};"X11"\ar@{..};"X13"
\ar@2{-}[lllu]\ar@{}[llll]+<0mm,-5mm>*{(1,0,0)}\ar@{}[lllu]+<-4mm,3mm>*{\scriptstyle{\mathbf{d}(z)}}
\ar@{..}"Z";"X6"\ar@{..}"Z";"X8"\ar@{..}"Z";"X10"\ar@{..}"Z";"X12"
\ar@{..}"X5";"X6"\ar@{..}"X6";"X7"\ar@{..}"X7";"X8"\ar@{..}"X8";"X9"\ar@{..}"X9";"X10"\ar@{..}"X10";"X11"\ar@{..}"X11";"X12"
\ar@{..}"X-1";"X-2"\ar@{..}"X-2";"X-3"\ar@{..}"X-3";"X-4"\ar@{..}"X-4";"X-5"\ar@{..}"X-5";"X-6"\ar@{..}"X-6";"X-7"\ar@{..}"X-7";"X-8"
}
\caption{ ``Cluster triangulation'' of the intersection between the positive octant $Q_{+}$ and the plane $\mathcal{P} = \{e_{1}+e_{2}+e_{3} = 1\}$. A dotted line joins two points corresponding to cluster variables belonging to the same cluster. The double line between $w$ and $z$ denotes the intersection with the ``regular'' cone which contains  the denominator vector of all the elements $u_{n}w^{k}$ and $u_{n}z^{k}$.}
\label{Fig:ClusterTriangulation}
\end{figure}
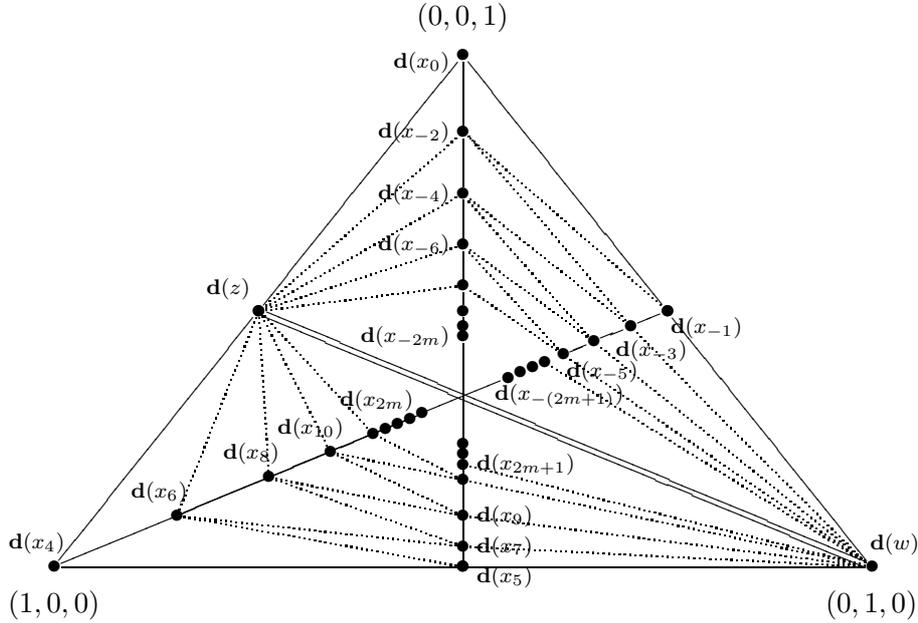

In lemma~\ref{Lemma:DenominatorVectors} we show the denominator vectors of the elements of $\BB$ in $\{x_1,x_2,x_3\}$. 
Figure~\ref{Fig:ClusterTriangulation} shows  the qualitative positions of denominator vectors of cluster variables (different from the initial ones) in $Q$: it shows the intersection between the plane $\mathcal{P}=\{e_{1} + e_{2} + e_{3} = 1\}$ and the
positive octant $Q_{+}$ of the real vector space $Q_{\RR}$ generated by $Q$; a point labeled by $\mathbf{d}(s)$  denotes the intersection between $\mathcal{P}$ and the line generated by the denominator vector of the cluster variable $s$. A dotted line joins two cluster variables that belong to the same cluster. In view of theorem~\ref{Prop:BijectionDenominators}, these lines form mutually disjoint triangles. We notice that some of the points of $\mathcal{P}\cap Q_{+}$ do not lie in one of these triangles. All such points belong to the line between $\mathbf{d}(w)$ and $\mathbf{d}(z)$ that is denoted by a double line in the figure. This is the intersection between $\mathcal{P}$ and the ``regular'' cone generated by $\mathbf{d}(w)$ and $\mathbf{d}(z)$. The denominator vector of $u_{n}w^{k}$ and $u_{n}z^{k}$, for all $n,k\geq0$, lie in the regular cone. Figure~\ref{Fig:ClusterTriangulation} appears also in \cite{DW} where the authors analyze the canonical decomposition of the elements of $Q$. Indeed denominator vectors of cluster variables are positive Schur roots (see \cite{DWZ}, \cite{DWZII}, \cite{CR}, \cite{BMRRT}).

\subsection{Explicit expressions of the elements of $\BB$}

Our next result provides explicit formulas for the elements of $\BB$ in every cluster of $\myAA_{\PP}$.  By the symmetry of the exchange relations it is sufficient to consider only the two clusters $\{x_{1},x_{2},x_{3}\}$ and $\{x_{1},w,x_{3}\}$ and only cluster variables $x_{m}$ with $m\geq2$ (see remark~\ref{Rem:Symmetries}). We use the notation $\mathbf{x}^\mathbf{e}:=x_1^{e_1}x_2^{e_2}x_3^{e_3}$.

\begin{theorem}\label{Thm:ExplicitExpressions}
Let $\PP$ be any semifield. For a cluster variable with denominator vector $\mathbf{d}=(d_1,d_2,d_3)^t$ in the cluster  $\{x_{1},x_{2},x_{3}\}$ we use the notation $$\varepsilon(\mathbf{e}):=(e_2+e_3,d_1-e_1+e_3,d_1+d_2-e_1-e_2)$$.  For every $m,n \geq 1$ the following formulas hold:
$$
x_{2m+1}=\frac{\sum_{\mathbf{e}}{e_{1}-e_{3}\choose e_{2}-e_{3}}{m-1-e_{3}\choose m-1-e_{1}}{e_{1}-1\choose e_{3}}\mathbf{y}^{\mathbf{e}}\mathbf{x}^{\varepsilon(\mathbf{e})}
+x_{2}^{m-1}x_{3}^{2m-2}}{\left(\scriptstyle{\bigoplus_{\mathbf{e}}{e_{1}-e_{3}\choose e_{2}-e_{3}}{m-1-e_{3}\choose m-1-e_{1}}{e_{1}-1\choose e_{3}}}\mathbf{y}^{\mathbf{e}}\oplus1\right)\;\;x_{1}^{m-1}x_{2}^{m-1}x_{3}^{m-2}}
$$
$$
x_{2m+2}\!=\!\frac{\sum_{\mathbf{e}}{e_{1}-1\choose e_{3}}{m-e_2\choose e_1-e_2}{m-1-e_3\choose e_2-e_3}\mathbf{y}^{\mathbf{e}}\mathbf{x}^{\varepsilon(\mathbf{e})}\!+\!x_{2}^{m}x_{3}^{2m-1} }{\left(\bigoplus_{\mathbf{e}}{e_{1}-1\choose e_{3}}{m-e_2\choose e_1-e_2}{m-1-e_3\choose e_2-e_3}\mathbf{y}^{\mathbf{e}}\oplus1\right)\;\;x_{1}^{m}x_{2}^{m-1}x_{3}^{m-1}}
$$
$$
u_{n}\!=\!\frac{y_{1}^{n}y_{2}^{n}y_{3}^{n}x_{1}^{2n}x_{2}^{n}\!+\!x_{2}^{n}x_{3}^{2n}\!\!+\!\!\sum_{\mathbf{e}}\!{e_{1}-e_{3}\choose e_{1}-e_{2}}\!\left[\!{n-e_{3}\choose n-e_{1}}{e_{1}-1\choose e_{3}}\!\!+\!\!{n-e_{3}-1\choose n-e_{1}}{e_{1}-1\choose e_{3}-1}\!\right]\mathbf{y}^{\mathbf{e}}\mathbf{x}^{\varepsilon(\mathbf{e})}}{x_{1}^{n}x_{2}^{n}x_{3}^{n}}
$$

Let $\yy_{1}:=y_{1}(y_{2}\oplus1)$, $\yy_{2}:=\frac{1}{y_{2}}$ and $\yy_{3}:=\frac{y_{2}y_{3}}{y_{2}\oplus1}$. For every $m,n \geq 1$:
$$
x_{2m+1}=\frac{\sum_{e_{1},e_{3}}{m-1-e_{3}\choose m-1-e_{1}}{e_{1}-1\choose e_{3}}\yy_{1}^{e_{1}}\yy_{3}^{e_{3}}x_{1}^{2e_{3}}w^{e_{1}-e_{3}}x_{3}^{2m-2e_{1}-2}+x_{3}^{2m-2}}{\left(\bigoplus_{e_{1},e_{3}}\scriptstyle{{m-1-e_{3}\choose m-1-e_{1}}{e_{1}-1\choose e_{3}}}\yy_{1}^{e_{1}}\yy_{3}^{e_{3}}\oplus 1\right)x_{1}^{m-1}x_{3}^{m-2}}
$$
$$
x_{2m+2}=\frac{\sum_{\mathbf{e}}\scriptstyle{{m-1-e_{3}+e_{2}\choose m-1-e_{1}+e_{2}}{e_{1}-1\choose e_{3}}{1\choose e_{2}}}\mathbf{\yy}^{\mathbf{e}}x_{1}^{2e_{3}+1-e_{2}}w^{e_{1}-e_{3}}x_{3}^{2m-2e_{1}-1+e_{2}}+x_{3}^{2m-1}(\yy_{2}x_{3}+x_{1})}{\left(\bigoplus_{\mathbf{e}}\scriptstyle{{m-1-e_{3}+e_{2}\choose m-1-e_{1}+e_{2}}{e_{1}-1\choose e_{3}}{1\choose e_{2}}}\mathbf{\yy}^{\mathbf{e}}\oplus \yy_{2}\oplus1\right)\;\;x_{1}^{m}wx_{3}^{m-1}}
$$
$$
u_{n}\!=\!\frac{\displaystyle{\yy_{1}^{n}\yy_{3}^{n}x_{1}^{2n}\!\!+\!\!x_{3}^{2n}\!\!+\!\!\sum_{e_{1},e_{3}}\!\!\left[\!\scriptstyle{{n-e_{3}\choose n-e_{1}}{e_{1}-1\choose e_{3}}\!+\!{n-e_{3}-1\choose n-e_{1}}{e_{1}-1\choose e_{3}-1}}\!\right]\!\yy_{1}^{e_{1}}\yy_{3}^{e_{3}}x_{1}^{2e_{3}}w^{e_{1}-e_{3}}x_{3}^{2n-2e_{1}}}}{x_{1}^{n}x_{3}^{n}}
$$
\begin{equation}\label{Eq:ZCyclic}
z=\frac{\yy_{1}\yy_{2}\yy_{3}x_{1}^{2}+\yy_{1}\yy_{2}^{2}\yy_{3}x_{1}x_{3}+\yy_{1}\yy_{2}w+\yy_{2}x_{3}^{2}+x_{1}x_{3}}{\left(\yy_{1}\yy_{2}\yy_{3}\oplus\yy_{1}\yy_{2}^{2}\yy_{3}\oplus\yy_{1}\yy_{2}\oplus\yy_{2}\oplus1\right)\;x_{1}wx_{3}}
\end{equation}
\end{theorem}
Our first proof of theorem~\ref{Thm:ExplicitExpressions} uses the  theory of cluster categories by computing explicitly cluster characters \cite{Thesis}. We do not use this approach here. Instead the strategy of our proof uses a parametrization of the elements of $\BB$ shown in the next section. 

\subsection{Parametrization of $\BB$ by $\mathbf{g}$--vectors}

In \cite[section~7]{FZIV} it is shown that cluster monomials are parametrized by some integer vectors called $\mathbf{g}$--vectors. In this section we extend such parametrization to all the elements of $\BB$. Given a polynomial $F\in\ZZ_{\geq_{0}}[z_{1},\cdots,z_{n}]$ with positive coefficients in $n$ variables and a semifield $\PP=(\PP,\cdot,\oplus)$, its \emph{evaluation} $F|_{\PP}(y_{1},\cdots,y_{n})$ at $(y_{1},\cdots,y_{n})\in\PP\times\cdots\times\PP$ is the element of $\PP$ obtained by replacing the addition of $\ZZ[z_{1},\cdots,z_{n}]$ with the auxiliary addition $\oplus$ of $\PP$ in the expression $F(y_{1},\cdots,y_{n})$. For example the evaluation of $F(z_{1},z_{2}):=z_{1}+1$ at $(y_{1},y_{2})\in Trop(y_{1},y_{2})\times Trop(y_{1},y_{2})$ is $1$. Let $\PP$ be a semifield and let $\myAA_\PP$ be a cluster algebra of type $A_2^{(1)}$ with coefficients in $\PP$ and let $\Sigma=\{H,\myCC, \{y_1,y_2,y_3\}\}$ be a seed of $\myAA_\PP$, so that $H$ is an exchange matrix, $\myCC$ is a cluster and $y_1$, $y_2$ and $y_3$ are arbitrarily chosen elements of $\PP$ (see section~\ref{Sec:Background}). We use the notation $\mathbf{s}^\mathbf{e}:=s_1^{e_1}s_2^{e_2}s_3^{e_3}$.

\begin{proposition}\label{Prop:Homogeneity}
 For every cluster monomial $b$ of $\myAA_\PP$ there exist a (unique) polynomial $F_{b}^{\myCC}\in \ZZ_{\geq0}[z_1,z_2,z_3]$ in three variables, with non--negative coefficients and constant term $1$ and a (unique) integer vector $\mathbf{g}_{b}^{\myCC}\in\ZZ^{3}$ such that the expansion of $b$ in  the cluster $\myCC=\{s_{1},s_{2},s_{3}\}$ is:
\begin{equation}
\label{Eq:b=FbgbInEveryCluster}
b =\frac{F_{b}^{\myCC} (y_1\mathbf{s}^{\mathbf{h}_{1}}, y_2\mathbf{s}^{\mathbf{h}_{2}},y_3 \mathbf{s}^{\mathbf{h}_{3}})}{F_{b}^{\myCC}|_{\PP} (y_1,y_2,y_3)}\mathbf{s}^{\mathbf{g}_{b}^{\myCC}}
\end{equation}
where $\mathbf{h}_i$ is the $i$--th column vector of the exchange matrix $H$. For every $n\geq1$  there exist a (unique) polynomial $F_{u_n}^{\myCC}\in \ZZ_{\geq0}[z_1,z_2,z_3]$ in three variables, with non--negative coefficients and constant term $1$ and a (unique) integer vector $\mathbf{g}_{u_n}^{\myCC}\in\ZZ^{3}$ such that the expansion of $u_n$ in  the cluster $\myCC$ is:
\begin{equation}\label{Eq:b=Fbgb}
u_n= F_{u_n}^{\myCC} (y_1\mathbf{s}^{\mathbf{h}_{1}}, y_2\mathbf{s}^{\mathbf{h}_{2}},y_3 \mathbf{s}^{\mathbf{h}_{3}})\mathbf{s}^{\mathbf{g}_{b}^{\myCC}}.
\end{equation}
\end{proposition}
In view of proposition~\ref{Prop:Homogeneity}, every element $b\in\BB$ determines a polynomial $F_b^\myCC$  and an integer vector $\mathbf{g}_b^\myCC$, for every cluster $\myCC$ of $\myAA_\PP$,. The polynomial $F_{b}^{\myCC}$ and the vector $\mathbf{g}_{b}^{\myCC}$ are called respectively the $F$--polynomial and the $\mathbf{g}$--vector of $b$ in the cluster $\myCC$. The previous proposition is the key result for our proof of theorem~\ref{Thm:ExplicitExpressions}. In section~\ref{subsec: Homogeneity} we find the explicit expression of $F$--polynomials and $\mathbf{g}$--vectors of every element of $\BB$ in every cluster.

We can define more intrinsically $F$--polynomials and $\mathbf{g}$--vectors as follows: consider the tropical semifield $\PP=Trop(y_{1},y_{2},y_{3})$ generated by the coefficients of $\Sigma$ and expand a cluster monomial $b\in\myAA_{\PP}$ in the cluster $\myCC$. Since $F_b^\myCC$ has constant term $1$, $F_b^\myCC|_\PP(y_1,y_2,y_3)=1$. It hence follows that if we replace $s_{1}$, $s_{2}$ and $s_{3}$ by $1$ in \eqref{Eq:b=FbgbInEveryCluster} we get $F_b^\myCC(y_1,y_2,y_3)$. 

The $\mathbf{g}$--vector of $b$ can be defined as follows: let $\PP=Trop(y_{1},y_{2},y_{3})$; following \cite{FZIV} we consider the \emph{principal $\ZZ^3$--grading of $\myAA_{\PP}$} given by
\begin{equation}\label{Eq:PrincipalGrading}
\begin{array}{ccc}
deg(s_{i})= \mathbf{e}_{i},& deg(y_{i}) = -\mathbf{h}_{i},&i=1, 2, 3
\end{array}
\end{equation}                               
($\mathbf{e}_{i}$ is the $i$--th basis vector of $\ZZ^3$). The element $\hat{y}_{i}:=y_{i} \mathbf{s}^{\mathbf{h}_{i}^{\myCC}}$, $i=1,2,3$, has  degree zero  with
respect to such grading. Therefore, since $F_b^\myCC|_\PP(y_1,y_2,y_3)=1$, every element $b$ of $\BB$ is homogeneous with respect to such grading and the $\mathbf{g}$--vector $\mathbf{g}_{b}^{\myCC}$ is its degree. 

The relation between denominator vectors and $\mathbf{g}$--vector in the cluster $\{x_{1},x_{2},x_{3}\}$ is given by the following proposition. 
 
\begin{proposition}\label{Prop:gbDb}
Let $b$ be an element of $\BB$ not divisible by $x_{1}$, $x_{2}$ or $x_{3}$. The $\mathbf{g}$--vector $\mathbf{g}_{b}$ of $b$ and its denominator vector $\mathbf{d}(b)$ both in the cluster $\{x_{1},x_{2},x_{3}\}$  are related by
\begin{equation}
\label{Eq:gbEQDb}
\mathbf{g}_{b}=\left(\def\objectstyle{\scriptscriptstyle}\def\labelstyle{\scriptscriptstyle}\vcenter{\xymatrix@R=0pt@C=0pt{-1 & 0 & 0\\1 & -1 & 0\\1 & 1 & -1}}\right)\mathbf{d}(b)
\end{equation}
If $b$ is a cluster monomial in the initial cluster $\{x_{1},x_{2},x_{3}\}$ then $\mathbf{g}_{b}=-\mathbf{d}(b)$.
\end{proposition}

Formula \eqref{Eq:gbEQDb} between the $\mathbf{g}$--vector and the denominator vector of a cluster variable in an ``acyclic'' seed can be deduced from results of \cite{FuKeller} and \cite{DWZII}. We notice that the map \eqref{Eq:gbEQDb} is bijective, and hence by combining proposition~\ref{Prop:BijectionDenominators} and theorem~\ref{Prop:BijectionDenominators}, we get that the map $b\mapsto \mathbf{g}_b$ restricts to a bijection between $\BB$ and $\ZZ^3$. What happens if we change cluster? The denominator vector map in the cluster $\{x_{1} , w, x_{3}\}$ is not surjective; on the other hand, as expected for cluster monomials (\cite[Conjecture 7.10]{FZIV}, proved in \cite[theorem~6.3]{FuKeller}), there is a bijective map (see \eqref{eq:RecursionGCyclic}) between the $\mathbf{g}$--vectors of the elements of $\BB$ in $\{x_1,x_2,x_3\}$ and in $\{x_1,w,x_3\}$. We have hence the following parametrization of $\BB$ by $\ZZ^3$ via the $\mathbf{g}$--vector map.

\begin{proposition}\label{Prop:GVectorPar}
Given a cluster $\mathcal{C}$ of $\myAA$, the map $b\mapsto \mathbf{g}_{b}^{\mathcal{C}}$ which sends an element $b$ of $\BB$ to its $\mathbf{g}$--vector $\mathbf{g}_b^{\mathcal{C}}$ in the
cluster $\mathcal{C}$, is a bijection between $\BB$ and $\ZZ^3$ 
\end{proposition}

The last comment about \eqref{Eq:gbEQDb} is the following: we note that the matrix in \eqref{Eq:gbEQDb} equals $-E_Q$ where $E_Q$ denotes the Euler matrix of the quiver $Q$ given in \eqref{Eq:Q}.   In \cite{SchofieldSemiInvariants} it is shown that the weight of the Schofield's semi--invariant $c_{V}$, where $V$
is an indecomposable $Q$--representation of dimension vector $\mathbf{d}=\mathbf{d}(b)$, is given by $-E\mathbf{d}$ and in view of \eqref{Eq:gbEQDb},  is the $\mathbf{g}$--vector of $b$. In particular one could deduce proposition~\ref{Prop:GVectorPar} from results of \cite{ITW}.
\newpage
\section{Cluster algebras of type $A_{2}^{(1)}$}\label{sec:AlgebraicStructure}

\subsection{Background on cluster algebras}\label{Sec:Background}
Let $\PP$ be a semifield (see section~\ref{sec:MainResults}).  Let $\FF_{\PP}=\QQ\PP(x_{1},\cdots,x_{n})$ be the field of rational functions in $n$ independent variables $x_{1}$, $\cdots$, $x_{n}$. A \emph{seed} in $\FF_{\PP}$ is a triple $\Sigma=\{H,\myCC,\mathbf{y}\}$  where $H$ is an $n\times n$ integer matrix which is skew--symmetrizable, i.e. there exists a diagonal matrix $D=diag(d_{1},\cdots,d_{n})$ with $d_{i}>0$  for all $i$ such that $DB$ is skew--symmetric; $\myCC=\{s_{1},\cdots,s_{n}\}$ is an n--tuple of elements of $\FF_{\PP}$ which is a free generating set for $\FF_{\PP}$ so that $\FF_{\PP}\simeq\QQ\PP(s_{1},\cdots,s_{n})$; and finally $\mathbf{y}=\{y_{1},\cdots,y_{n}\}$ is an n-tuple of elements of $\PP$. The matrix $H$ is called the \emph{exchange matrix} of $\Sigma$, the set $\myCC$ is called the \emph{cluster} of $\Sigma$ and its elements are called \emph{cluster variables} of $\Sigma$ and the set $\mathbf{y}$ is called the coefficients tuple of the seed $\Sigma$. 

We fix an integer $k\in [1,n]$. Given a seed $\Sigma$ of $\FF_{\PP}$ it is defined another seed $\Sigma_{k}:=\{H_{k},\myCC_{k},\mathbf{y}_{k}\}$ by the following \emph{mutation rules} (see \cite{FZIV}):   
\begin{enumerate}
\item\label{MatrixMutation} the exchange matrix $H_{k}=(h'_{ij})$ is obtained from $H=(h_{ij})$ by
\begin{equation}\label{Eq:MatrixMutation}
h'_{ij}=\left\{\begin{array}{cc} -h_{ij}&\text{ if }i=k\text{ or }j=k\\
h_{ij}+sg(h_{ik})[h_{ik}h_{kj}]_{+}&\text{ otherwise}\end{array}\right.
\end{equation}
where $[c]_{+}:=\text{max}(c,0)$ for every integer $c$;
\item the new coefficients tuple $\mathbf{y}_{k}=\{y'_{1},\cdots,y'_{n}\}$ is given by:
\begin{equation}\label{eq:CoefficientsMutation}
y'_{j}:=\left\{\begin{array}{cc}\frac{1}{y_{k}}&\text{ if }j=k\\ y_{j}y_{k}^{[h_{kj}]_{+}}(y_{k}\oplus1)^{-h_{kj}}&\text{ otherwise.}\end{array}\right.
\end{equation}
\item the new cluster $\myCC_{k}$ is given by $\myCC_{k}=\myCC\setminus \{s_{k}\}\cup\{s'_{k}\}$ where
\begin{equation}\label{eq:ClusterMutation}
s'_{k}:=\frac{y_{k}\prod_{i}s_{i}^{[h_{ik}]_{+}}+\prod_{i}s_{i}^{[-h_{ik}]_{+}}}{(y_{k}\oplus 1) s_{k}};
\end{equation} 
\end{enumerate}
It is not hard to verify that $\myCC_{k}$ is again a seed of $\FF_{\PP}$. We say that the seed $\Sigma_{k}$ is obtained from the seed $\Sigma$ by a \emph{mutation in direction $k$}. Every seed can be mutated in all the directions.  

Given a seed $\Sigma$ we consider the set $\chi(\Sigma)$ of all the cluster variables of all the seeds obtained by a sequence of mutations.  
The \emph{rank $n$ cluster algebra with initial seed $\Sigma$ with coefficients in $\PP$} is by definition the $\ZZ\PP$--subalgebra of $\FF_{\PP}$ generated by $\chi(\Sigma)$; we denote it by $\myAA_{\PP}(\Sigma)$.

The cluster algebra $\myAA_{\PP}(\Sigma)$ is said to have \emph{principal coefficients at $\Sigma$} and it is denoted by $\myAA_{\bullet}(\Sigma)$ if $\PP=Trop(y_{1},\cdots y_{n})$. 

If the semifield $\PP=Trop(c_1,\cdots,c_r)$ is a tropical semifield, the elements of $\PP$ are monomials in the $c_{j}$'s. Therefore the coefficient $y_{j}$ of a seed $\Sigma=\{H,\myCC, \{y_{1},\cdots,y_{n}\}\}$ is a monomial of the form $y_{j}=\prod_{i=1}^{r}c_{i}^{h_{n+r,j}}$ for some integers $h_{n+1,j}$, $\cdots$, $h_{n+r,j}$. It is convenient to ``complete" the exchange $n\times n$ matrix $H$ to a rectangular $(n+r)\times n$ matrix $\tilde{H}$ whose $(i,j)$--th entry is $h_{ij}$. The seed $\Sigma$ can hence be seen as a couple $\{\tilde{H},\myCC\}$ and the mutation of the coefficients tuple \eqref{eq:CoefficientsMutation} translates into the mutation \eqref{Eq:MatrixMutation} of the rectangular matrix $\tilde{H}$. We sometimes use this formalism.   

We say that two seeds $\Sigma=\{H, \myCC, \mathbf{y}\}$ and $\Sigma'=\{H', \myCC', \mathbf{y}'\}$ of a cluster algebra of rank $n$ are equivalent if there exists a permutation $\sigma$ of the index set $[1,n]:=\{1,\cdots,n\}$ such that  $h_{\sigma(i)\sigma(j)}=h'_{ij}$, $s_{\sigma(i)}=s'_i$ and $y_{\sigma(i)}=y'_i$ for every $i,j\in[1,n]$. A class of equivalent seeds is called an  \emph{unlabeled} seeds. We often consider unlabeled seeds in this paper.  We use the cyclic representation of a permutation $\sigma$ so that $\sigma=(i_{1},i_{2},\cdots,i_{t})$ denotes the permutation $\sigma$ such that $\sigma(i_{k})=i_{k+1}$ if $k=1,\cdots,t-1$ and $\sigma(i_{t})=i_{1}$; all the other indices are fixed by $\sigma$.

\begin{example}\label{ex:UnlabeledSeeds}
The seed $$\Sigma:=\{H=\left(\def\objectstyle{\scriptscriptstyle}\def\labelstyle{\scriptscriptstyle}\vcenter{\xymatrix@R=0pt@C=0pt{0&1&1\\-1&0&1\\-1&-1&0}}\right),\mathbf{x}=\{x_1,x_2,x_3\},\mathbf{y}=\{y_1,y_2,y_3\}\}$$ is equivalent to the seed $$\Sigma':=\{H'=\left(\def\objectstyle{\scriptscriptstyle}\def\labelstyle{\scriptscriptstyle}\vcenter{\xymatrix@R=0pt@C=0pt{0&-1&-1\\1&0&1\\1&-1&0}}\right),\mathbf{x}'=\{x_3,x_1,x_2\},\mathbf{y}'=\{y_3,y_1,y_2\}\}$$ by the permutation  $(132)$.
\end{example}
We remark that even if every mutation of a seed in a fixed direction is involutive, this is not true for unlabeled seeds.

\subsection{Algebraic structure of $\myAA_{\PP}$}\label{Sec:AlgStructure}
Let $\PP=(\PP,\cdot,\oplus)$ be a semifield and $\myAA_{\PP}$ be the cluster algebra with initial seed 
$$
\Sigma:=\{H=\left(\def\objectstyle{\scriptstyle}\def\labelstyle{\scriptstyle}\vcenter{\xymatrix@R=0pt@C=0pt{0&1&1\\-1&0&1\\-1&-1&0}}\right), \displaystyle{\mathbf{x}=\{x_1,x_2,x_3\}, \mathbf{y}=\{y_1,y_2,y_3\}}\}.
$$ 
The following lemma gives the algebraic structure of $\myAA_{\PP}$.
\begin{lemma}\label{LemmaAlgStructure}
The unlabeled seeds of the cluster algebra $\myAA_{\PP}$ with initial seed $\Sigma=\Sigma_{1}$ are the following:
\begin{eqnarray}\label{Def:SeedSigmam}
\Sigma_{m}&:=&\{H_{m}, \{x_{m}, x_{m+1}, x_{m+2}\},\{y_{1;m},y_{2;m},y_{3;m}\}\},\\
\label{Def:SeedSigmaW}
\Sigma^{cyc}_{2m-1}\!\!\!&:=&\{H_{m}^{cyc}, \{x_{2m-1},w, x_{2m+1}\},\{y_{1;2m-1}^{cyc},y_{2;2m-1}^{cyc},y_{3;2m-1}^{cyc}\}\},\\
\label{Def:SeedSigmaZ}
\Sigma^{cyc}_{2m}\!\!\!&:=&\{H_{m}^{cyc}, \{x_{2m}, z, x_{2m+2}\},\{y_{1;2m}^{cyc},y_{2;2m}^{cyc},y_{3;2m}^{cyc}\}\}
\end{eqnarray}
for every $m\in\ZZ$; they are mutually related by the following diagram of
mutations:
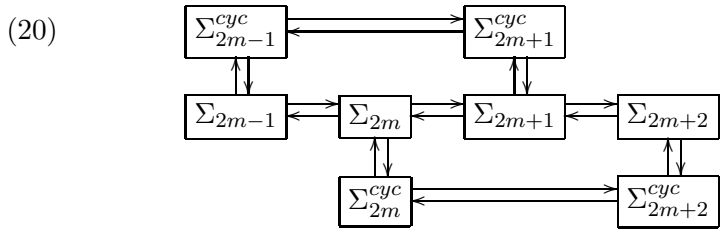
\begin{equation}\label{DiagramMutations}
\xymatrix@!R=8pt@C=20pt{                       
*+[F]{\Sigma^{cyc}_{2m-1}}\ar@<1ex>[rr]\ar@<1ex>[d]&   &  *+[F]{\Sigma^{cyc}_{2m+1}}\ar@<1ex>[d]\ar[ll]    &                 &\\
*+[F]{\Sigma_{2m-1}}\ar@<1ex>[r]\ar[u]    &*+[F]{\Sigma_{2m}}\ar@<1ex>[r]\ar@<1ex>[d]\ar[l]         &   *+[F]{\Sigma_{2m+1}}\ar@<1ex>[r]\ar[l]\ar[u]       & *+[F]{\Sigma_{2m+2}}\ar[l]\ar@<1ex>[d]\\
                              & *+[F]{\Sigma^{cyc}_{2m}}\ar@<1ex>[rr]\ar[u]   &                                   & *+[F]{\Sigma^{cyc}_{2m+2}}\ar[u]\ar[ll]
}
\end{equation}      
where  arrows from left to right (resp. from right to left) are mutations in direction~$1$
(resp.~$3$) and vertical arrows (in both directions) are mutations in direction~$2$. The seed $\Sigma_{m}$ is not equivalent to the seed $\Sigma_{n}$ if $m\neq n$, in particular the exchange graph of $\myAA_{\PP}$ is given by figure~\ref{fig:ExchangeGraph} and every cluster $\myCC$ determines a unique seed $\Sigma^{\myCC}$.
The exchange matrices $H_{m}$ and $H_{m}^{cyclic}$ are the following: 
\begin{equation}\label{eq:ExchangeMatrices}
\begin{array}{cc}
H_{m}=\left(\def\objectstyle{\scriptstyle}\def\labelstyle{\scriptstyle}\vcenter{\xymatrix@R=0pt@C=0pt{0&1&1\\-1&0&1\\-1&-1&0}}\right),&
H_{m}^{Cyc}=\left(\def\objectstyle{\scriptstyle}\def\labelstyle{\scriptstyle}\vcenter{\xymatrix@R=0pt@C=0pt{0&-1&2\\1&0&-1\\-2&1&0}}\right)
\end{array}
\end{equation}
for every $m\in\ZZ$. The exchange relations for the coefficient tuples are the following:
\begin{equation}\label{eq:InitialCoefficients}
\begin{array}{ccc}
y_{1;1}=y_{1};&y_{2;1}=y_{2};&y_{3;1}=y_{3}
\end{array}
\end{equation}
\begin{equation}\label{CoeffMutationDirection1}
\begin{array}{ccc}
y_{1;m+1}=\frac{y_{2;m}y_{1;m}}{y_{1;m}\oplus1};&y_{2;m+1}=\frac{y_{3;m}y_{1;m}}{y_{1;m}\oplus1};&y_{3;m+1}=\frac{1}{y_{1;m}}
\end{array}
\end{equation}
\begin{equation}\label{CoeffMutationDirection3}
\begin{array}{ccc}
y_{1;m-1}=\frac{1}{y_{3;m}};&y_{2;m-1}=y_{1;m}(y_{3;m}\oplus1);&y_{3;m-1}=y_{2;m}(y_{3;m}\oplus1)
\end{array}
\end{equation}
\begin{equation}\label{CoeffMutationDirection2}
\begin{array}{ccc}
y_{1;m}^{cyc}=\frac{y_{1;m}}{y_{2;m}\oplus1};&y_{2;m}^{cyc}=\frac{1}{y_{2;m}};&y_{3;m}^{cyc}=\frac{y_{3;m}y_{2;m}}{y_{2;m}\oplus1}
\end{array}
\end{equation}
\begin{equation}\label{CoeffMutationCyclicDirection1}
\begin{array}{ccc}
y_{1;m+2}^{cyc}=y_{3;m}^{cyc}(y_{1;m}^{cyc}\oplus1)^{2};&y_{2;m+2}^{cyc}=\frac{y_{2;m}^{cyc}y_{1;m}^{cyc}}{y_{1;m}^{cyc}\oplus1};&y_{3;m+2}^{cyc}=\frac{1}{y_{1;m}^{cyc}}\end{array}
\end{equation}
\begin{equation}\label{CoeffMutationCyclicDirection2}
\begin{array}{ccc}
y_{1;m}=\frac{y_{1;m}^{cyc}y_{2;m}^{cyc}}{y_{2;m}^{cyc}\oplus1};&y_{2;m}=\frac{1}{y_{2;m}^{cyc}};&y_{3;m}=y_{3;m}^{cyc}(y_{2;m}^{cyc}\oplus1)
\end{array}
\end{equation}
\begin{equation}\label{CoeffMutationCyclicDirection3}
\begin{array}{ccc}
y_{1;m-2}^{cyc}=\frac{1}{y_{3;m}^{cyc}};&y_{2;m-2}^{cyc}=\frac{y_{2;m}^{cyc}y_{3;m}^{cyc}}{y_{3;m}^{cyc}\oplus1};&y_{3;m-2}^{cyc}=y_{1;m}^{cyc}(y_{3;m}^{cyc}\oplus1)^{2}
\end{array}
\end{equation}

The exchange relations for the cluster variables are the following:
\begin{eqnarray}\label{ExchRel:XmXm+3}
x_{m}x_{m+3}\!\!\!&=&\!\!\!\frac{x_{m+1} x_{m+2} + y_{1;m}}{y_{1;m}  \oplus 1}=\frac{y_{3;m+1} x_{m+1} x_{m+2} + 1}{y_{3;m+1}  \oplus 1}\\ 
\label{ExchRel:WX2m}wx_{2m}\!\!\!&=&\!\!\!\frac{y_{2;2m-1} x_{2m-1} + x_{2m+1}}{y_{2;2m-1}  \oplus 1}=\frac{x_{2m-1} + y_{2;2m-1}^{cyc} x_{2m+1}}{y_{2;2m-1}^{cyc}  \oplus 1}\\
\label{ExchRel:ZX2m+1} zx_{2m+1}\!\!\!&=&\!\!\!\frac{y_{2;2m} x_{2m} + x_{2m+2}}{y_{2;2m}  \oplus 1}=\frac{x_{2m} + y_{2;2m}^{cyc} x_{2m+2}}{y_{2;2m}^{cyc}  \oplus 1}\\
\label{ExchRel:X2m-2X2m+2} x_{2m-2} x_{2m+2}\!\!\!&=&\!\!\!\frac{x_{2m}^{2} + y_{1;2m-2}^{cyc}z}{y_{1;2m-2}^{cyc}  \oplus 1}=\frac{y_{3;2m}^{cyc} x_{2m}^{2} + z}{y_{3;2m}^{cyc}  \oplus 1}\\
\label{ExchRel:X2m-1X2m+3} x_{2m-1} x_{2m+3}\!\!\!&=&\!\!\!\frac{ x_{2m+1}^{2}+ y_{1;2m-1}^{cyc}w}{y_{1;2m-1}^{cyc}  \oplus 1}=\frac{y_{3;2m+1}^{cyc} x_{2m+1}^{2}+w}{y_{3;2m+1}^{cyc}  \oplus 1}
\end{eqnarray}
\end{lemma}
\begin{proof}[Proof of Lemma~\ref{LemmaAlgStructure}]  
We need to prove the diagram~\eqref{DiagramMutations} for every $m\in\ZZ$. Clearly $\Sigma_{1}$ coincides with the initial seed $\Sigma$. For $m\geq2$ let $\Sigma_{m}$ (resp. $\Sigma_{-m+1}$) be the seed of $\myAA_{\PP}$ obtained from $\Sigma_{1}$ by applying $m-1$ times the following operation: first we mutate in direction $1$ (resp. $3$) and then we reorder the index set of the obtained seed by the permutation $(132)$ (resp. $(123)$) of the index set (as in example~\ref{ex:UnlabeledSeeds}). Assume that $\Sigma_{m}$ (resp.$\Sigma_{-m+1}$) has the form \eqref{Def:SeedSigmam}, i.e. its exchange matrix is $H_{m}$ (resp. $H_{-m+1}$), its cluster is $\{x_{m},x_{m+1},x_{m+2}\}$ and its coefficient tuple is $\{y_{1;m},y_{2;m},y_{3;m}\}$. Then it is 
straightforward to check by induction on $m$ that $H_{m}$ is given by \eqref{eq:ExchangeMatrices} and the exchange relations passing from $\Sigma_{m}$ to $\Sigma_{m+1}$ are given by \eqref{CoeffMutationDirection1} (resp. \eqref{CoeffMutationDirection3}) for the coefficient tuple and by \eqref{ExchRel:XmXm+3} for the cluster variables. 

Now it is straightforward to verify that, for every $m\in\ZZ$, $\Sigma_{m}$ is obtained from $\Sigma_{m+1}$ by mutating in direction $3$ and then by reordering with the permutation $(123)$. 

The central line of the diagram is hence proved. 

Let $\Sigma_{m}^{cyc}$ be the seed obtained from $\Sigma_{m}$ by the mutation in direction $2$. Suppose that $\Sigma_{m}^{cyc}$ has the form $\{H_{m}^{cyc},\{x_{m},c_{m},x_{m+2}\},\{y_{1;m}^{cyc},y_{2;m}^{cyc},y_{3;m}^{cyc}\}\}$. Then it is straightforward to verify the following: that by \eqref{Eq:MatrixMutation} the exchange matrix $H_{m}^{cyc}$ is given by \eqref{eq:ExchangeMatrices}; that by \eqref{eq:CoefficientsMutation} the coefficient tuple satisfy \eqref{CoeffMutationDirection2} and by the exchange relation \eqref{eq:ClusterMutation} the cluster variable $c_{m}$ is given by:
\begin{equation}\label{eq:c=xM+xM+2}
c_{m}=\frac{y_{2;m} x_{m} + x_{m+2}}{(y_{2;m}  \oplus 1)x_{m+1}}.
\end{equation}
By using \eqref{CoeffMutationDirection1}, \eqref{CoeffMutationDirection3} and \eqref{ExchRel:XmXm+3}, it is straightforward to verify that $c_{m}=c_{m+2}$ for every $m\in\ZZ$. We define $w:=c_{1}$ and $z:=c_{2}$ and we hence find that $\Sigma_{2m-1}^{cyc}$ has the form \eqref{Def:SeedSigmaW} and $\Sigma_{2m}^{cyc}$ has the form \eqref{Def:SeedSigmaZ}. Now since two cluster variables can belong to at most two clusters, we conclude that the mutation in direction $1$ (resp. $3$) of $\Sigma_{m}^{cyc}$ is $\Sigma_{m+2}^{cyc}$ (resp. $\Sigma_{m-2}^{cyc}$). 

We have hence proved the diagram~\eqref{DiagramMutations} and the fact that every cluster determines a unique seed.

It remains to prove that the seeds $\{\Sigma_{m}\}$ (resp.$\{\Sigma_{m}^{cyc}\}$) are not equivalent for every $m\in\ZZ$. It is sufficient to prove that $x_{m}$ is different from $x_{1}$, $x_{2}$ and $x_{3}$ for every $m\geq 4$. In view of the exchange relation \eqref{ExchRel:XmXm+3}, the denominator vector of $x_{m}$, $m\geq4$, in the cluster $\{x_{1},x_{2},x_{3}\}$ satisfies the initial conditions $\mathbf{d}(x_{4})=\mathbf{e}_{1}$, $\mathbf{d}(x_{5})=\mathbf{e}_{1}+\mathbf{e}_{2}$ and
\begin{equation}\label{eq:RecurrenceDenominator}
\mathbf{d}(x_{m+3})+\mathbf{d}(x_{m})=\mathbf{d}(x_{m+1})+\mathbf{d}(x_{m+2}).
\end{equation}   
The solution of this recursion is given in lemma~\ref{Lemma:DenominatorVectors} below and it is clearly not periodic.
\end{proof}

\begin{remark}\label{Rem:Symmetries}
The expansion of a cluster variable $x_{m+n}$ (resp. $x_{2m+n}$) in the cluster $\{x_{m},c,x_{m+2}\}$ for $c=w$ or $c=x_{m+1}$, (resp. $\{x_{2m},z,x_{2m+2}\}$) is obtained by the expansion of $x_{1+n}$ (resp. $x_{2m+1+n}$) in the cluster $\{x_{1},c,x_{3}\}$ (resp. $\{x_{2m+1},w,x_{2m+3}\}$) by replacing $x_{1}$ with $x_{m}$, $c$ with $x_{2}$ when $c\neq w$, $x_{3}$ with $x_{m+2}$ and $y_{i}$ with $y_{i;m}$ (resp. $x_{2m+1}$ with $x_{2m}$, $w$ with $z$, $x_{2m+3}$ with $x_{2m+2}$ and $y_{i;2m+1}^{cyc}$ with $y_{i;2m}^{cyc}$) for $i=1,2,3$ and $n$, $m\in\ZZ$. Moreover the expansion of $x_{-m+2}$ is obtained from the expansion of $x_{m+2}$ by replacing $x_{1}$ with $x_{3}$, $x_{3}$ with $x_{1}$ and $y_{1}$ with $y_{3}^{-1}$, $y_{2}$ with $y_{2}^{-1}$ and $y_{3}$ with $y_{1}^{-1}$. 
\end{remark} 
In the sequel we abbreviate and we write:
\begin{equation}\label{Def:CyclicSeed}
\Sigma^{cyc}:=\{\left(\def\objectstyle{\scriptscriptstyle}\def\labelstyle{\scriptscriptstyle}\vcenter{\xymatrix@R=0pt@C=0pt{0&-1&2\\1&0&-1\\-2&1&0}}\right),\{x_{1},w,x_{3}\},\{\yy_{1},\yy_{2},\yy_{3}\}\}
\end{equation}
for the seed $\Sigma^w_1$ obtained from the seed \eqref{Eq:InitialSeed} by a mutation in direction two. We sometimes call $\Sigma^{cyc}$ the cyclic seed of $\myAA_\PP$.

\subsubsection{Denominator vectors}
In this subsection we compute denominator vectors of the elements of $\BB$ in every cluster of $\myAA_{\PP}$. By the symmetry in the exchange relations it is sufficient to consider only the two clusters $\{x_{1},x_{2},x_{3}\}$ and $\{x_{1},w,x_{3}\}$.

\begin{lemma}\label{Lemma:DenominatorVectors}
Denominator vectors in the cluster $\{x_{1},x_{2},x_{3}\}$ of cluster variables different from $x_{1}$, $x_{2}$ and $x_{3}$ are the following:  for $m \geq 1$
\begin{eqnarray}\label{eq:D(XmPositive)}
\mathbf{d}(x_{2m+1})=\left[\def\objectstyle{\scriptstyle}\def\labelstyle{\scriptstyle}\vcenter{\xymatrix@R=0pt@C=0pt{m-1\\m-1\\m-2}}\right]
&&\mathbf{d}(x_{2m+2})=\left[\def\objectstyle{\scriptstyle}\def\labelstyle{\scriptstyle}\vcenter{\xymatrix@R=0pt@C=0pt{m\\m-1\\m-1}}\right]\\ \label{eq:D(XmNegative)}
\mathbf{d}(x_{-2m+1})=\left[\def\objectstyle{\scriptstyle}\def\labelstyle{\scriptstyle}\vcenter{\xymatrix@R=0pt@C=0pt{m-1\\m\\m}}\right]
&&\mathbf{d}(x_{-2m+2})=\left[\def\objectstyle{\scriptstyle}\def\labelstyle{\scriptstyle}\vcenter{\xymatrix@R=0pt@C=0pt{m-1\\m-1\\m}}\right]\\
\label{eq:DenWZ}
\mathbf{d}(w)=\left[\def\objectstyle{\scriptstyle}\def\labelstyle{\scriptstyle}\vcenter{\xymatrix@R=0pt@C=0pt{0\\1\\0}}\right]&&\mathbf{d}(z)=\left[\def\objectstyle{\scriptstyle}\def\labelstyle{\scriptstyle}\vcenter{\xymatrix@R=0pt@C=0pt{1\\0\\1}}\right]
\end{eqnarray}
For every $n\geq 1$ the denominator vector (in $\{x_{1},x_{2},x_{3}\}$) of $u_{n}$ is given by
\begin{equation}\label{eq:D(un)}
\mathbf{d}(u_{n})=\left[\def\objectstyle{\scriptstyle}\def\labelstyle{\scriptstyle}\vcenter{\xymatrix@R=0pt@C=0pt{n\\n\\n}}\right].
\end{equation}
In particular the denominator vector of all the cluster variables and of all the elements $\{u_{n}w,\,u_{n}z|\,n\geq1\}$ are all the positive real roots of the root system of type $A_{2}^{(1)}$.
\end{lemma}
\begin{proof}
In view of remark~\ref{Rem:Symmetries} we have $\mathbf{d}(x_{-m+2})=(13)\mathbf{d}(x_{m+2})$ for every $m\geq2$. Then \eqref{eq:D(XmNegative)} follows from \eqref{eq:D(XmPositive)}. We hence  prove \eqref{eq:D(XmPositive)}, \eqref{eq:DenWZ} and \eqref{eq:D(un)}.
By induction on $m\geq4$, one verifies that \eqref{eq:D(XmPositive)} is the solution of the recurrence relation  \eqref{eq:RecurrenceDenominator} together with the initial condition $\mathbf{d}(x_{4})=\mathbf{e}_{1}$, $\mathbf{d}(x_{5})=\mathbf{e}_{1}+\mathbf{e}_{2}$  (being $\mathbf{e}_{i}$ the $i$--th standard basis vector of $\ZZ^{3}$).  

The equalities in \eqref{eq:DenWZ} follow from \eqref{Eq:DefinitionWZ}.

It remains to prove \eqref{eq:D(un)}.  We notice that the denominator vector map $s\mapsto\mathbf{d}(s)$ in a cluster $\mathcal{C}$ is additive, i.e. $\mathbf{d}(s_1s_2)=\mathbf{d}(s_1)+\mathbf{d}(s_2)$. From its definition \eqref{Eq:DefinitionU0u_1U2}, it hence follows that the denominator vector of $u_{1}$ is the sum of $\mathbf{d}(w)$ and of $\mathbf{d}(z)$, i.e. $\mathbf{d}(u_{1})=\delta=(1,1,1)^{t}$. Then from \eqref{Eq:DefinitionUn+1}, $\mathbf{d}(u_{n})=n\mathbf{d}(u_{1})=n\delta$. 

The proof follows now by knowing the structure of a root system of type $A_{2}^{(1)}$, recalled in section~\ref{sec:MainResults}. 
\end{proof}
 
In view of a general well known result due to V.~Kac \cite{KacInfinite2}, for every positive real root $\mathbf{d}$ of a root system of type $A_2^{(1)}$ there exists a unique indecomposable representation $M(\mathbf{d})$ of the quiver $Q$ given in \eqref{Eq:Q} whose dimension vector is $\mathbf{d}$. In particular, in view of lemma~\ref{Lemma:DenominatorVectors}, for every element $b$ which is either a cluster variable or  $u_nw$ or $u_nz$ ( $n\geq0$) there exists a unique $Q$--representation  $M(b)$ whose dimension vector is $\mathbf{d}(b)$. If $b$ is a cluster variable, this is a special case of \cite[theorem~3]{CK2}.
\begin{corollary}\label{Cor:D(x)+D(u)}
For every $m\geq4$, $\mathbf{d}(x_m)+\mathbf{d}(u_1)=\mathbf{d}(x_{m+2})$. For every $m\leq0$, $\mathbf{d}(x_m)+\mathbf{d}(u_1)=\mathbf{d}(x_{m-2})$.
\end{corollary}
The following lemma shows  the denominator vectors of the elements of $\BB$ in the cyclic seed.

\begin{lemma}\label{lemma:DenominatorsCyclic}
Let $\mathbf{d}^{w}(b)$ be the denominator vector of $b\in\BB$ in the cluster $\{x_{1},w,x_{3}\}$. The following formulas hold:  for $m \geq 1$
\begin{eqnarray}\label{eq:D(XmPositiveCyclic)}
\mathbf{d}^{w}(x_{2m+1})=\left[\def\objectstyle{\scriptstyle}\def\labelstyle{\scriptstyle}\vcenter{\xymatrix@R=0pt@C=0pt{m-1\\0\\m-2}}\right]
&&\mathbf{d}^{w}(x_{2m+2})=\left[\def\objectstyle{\scriptstyle}\def\labelstyle{\scriptstyle}\vcenter{\xymatrix@R=0pt@C=0pt{m\\1\\m-1}}\right]\\ \label{eq:D(XmNegative)Cyclic}
\mathbf{d}^{w}(x_{-2m+1})=\left[\def\objectstyle{\scriptstyle}\def\labelstyle{\scriptstyle}\vcenter{\xymatrix@R=0pt@C=0pt{m-1\\0\\m}}\right]
&&\mathbf{d}^{w}(x_{-2m+2})=\left[\def\objectstyle{\scriptstyle}\def\labelstyle{\scriptstyle}\vcenter{\xymatrix@R=0pt@C=0pt{m-1\\1\\m}}\right]\\
\label{eq:DenWZCyclic}
\mathbf{d}^{w}(x_{2})=\left[\def\objectstyle{\scriptstyle}\def\labelstyle{\scriptstyle}\vcenter{\xymatrix@R=0pt@C=0pt{0\\1\\0}}\right]&&\mathbf{d}^{w}(z)=\left[\def\objectstyle{\scriptstyle}\def\labelstyle{\scriptstyle}\vcenter{\xymatrix@R=0pt@C=0pt{1\\1\\1}}\right]
\end{eqnarray}
For every $n\geq 1$ the denominator vector (in $\{x_{1},w,x_{3}\}$) of $u_{n}$ is given by
\begin{equation}\label{eq:D(un)Cyclic}
\mathbf{d}^w(u_{n})=\left[\def\objectstyle{\scriptstyle}\def\labelstyle{\scriptstyle}\vcenter{\xymatrix@R=0pt@C=0pt{n\\0\\n}}\right].
\end{equation}
\end{lemma}
\begin{proof}
In view of remark~\ref{Rem:Symmetries}, $\mathbf{d}^w(x_{-m+2})=(13)\mathbf{d}^w(x_{m+2})$ for every $m\geq2$. In particular \eqref{eq:D(XmNegative)Cyclic} follows from \eqref{eq:D(XmPositiveCyclic)}. We hence prove \eqref{eq:D(XmPositiveCyclic)}, \eqref{eq:DenWZCyclic} and \eqref{eq:D(un)Cyclic}. From the exchange relations it follows that  $\mathbf{d}^w(x_4)=\mathbf{e}_1+\mathbf{e}_2$ and $\mathbf{d}^w(x_5)=\mathbf{e}_1$. From the exchange relation \eqref{ExchRel:XmXm+3} it follows that the sequence of denominator vectors $\mathbf{d}^w(x_m)$ with $m\geq4$ satisfies the relation:
$$
\mathbf{d}^w(x_{m+3})+\mathbf{d}^w(x_{m})=\mathbf{d}^w(x_{m+1})+\mathbf{d}^w(x_{m+2}).
$$ 
By induction on $m\geq4$ we verify that \eqref{eq:D(XmPositiveCyclic)} is the unique solution of this recurrence. 

The equalities in \eqref{eq:DenWZCyclic} follows by direct check. 

The equality \eqref{eq:D(un)Cyclic} follows from the definition~\ref{Eq:DefinitionU0u_1U2}.
\end{proof}

Many properties of denominator vectors of cluster variables can be found in \cite{BMR} and \cite{BMRden}. Here we notice two  of them. 

\begin{remark}\label{Rem:NoBijectionDen}
The denominator vector map in the cluster $\{x_1,w,x_3\}$ does not restrict to a bijection between $\BB$ and $\ZZ^3$. Indeed there is a fractal region of $\ZZ^3$ which cannot be covered \cite{DW}. For example the element $(1,2,3)^t$ cannot be expressed as a non--negative linear combination of denominator vectors of adjacent (in the exchange graph ) cluster variables.  
\end{remark}

\begin{definition}\cite[Definition~6.12]{FZIV}\label{Def:SignCoherence}
A collection of vectors in $\ZZ^{n}$ (or in $\mathbb{R}^{n}$) are
\emph{sign--coherent} (to each other) if, for every $i \in \{1, \cdots , n\}$, the $i$--th coordinates of all of these vectors are either all non--negative or all non--positive.
\end{definition}
\begin{corollary}\label{Cor:SignCoherence} 
Denominator vectors (in every cluster) of cluster variables  belonging to the same cluster are sign--coherent.
\end{corollary}
\begin{proof}
By lemmas~\ref{Lemma:DenominatorVectors} and \ref{lemma:DenominatorsCyclic} one verifies directly that for every cluster $\{s_1,s_2,s_3\}$ the sets $\{\mathbf{d}(s_1),\mathbf{d}(s_2),\mathbf{d}(s_3)\}$ and  $\{\mathbf{d}^w(s_1),\mathbf{d}^w(s_2),\mathbf{d}^w(s_3)\}$ of corresponding denominator vectors are sign--coherent.
\end{proof}

\subsection{Explicit expression of the coefficient tuples}\label{sec:Coefficients}

In this section we solve the recurrence relations \eqref{eq:InitialCoefficients}--\eqref{CoeffMutationCyclicDirection3} for the coefficient tuples of the seeds of $\myAA_{\PP}$, for every semifield $\PP$. Such relations form a Y--system and are studied in several papers (see \cite{Ysyst}, \cite{NakanishiKeller} and references therein). The solution of this recurrence is given in terms of denominator vectors and of $F$--polynomials in the cluster $\{x_{1},x_{2},x_{3}\}$ of the cluster variables of $\myAA_{\PP}$. We assume proposition~\ref{Prop:Homogeneity}, even if its proof will be given later in section~\ref{subsec: Homogeneity}. Then $F_{s}$ is a polynomial with positive coefficients. In particular we can consider its evaluation $F_{s}|_{\PP}$ at $\PP$.  

\begin{proposition}\label{Prop:CoefficientsRecurrence}
\begin{itemize}
\item
The family $\{y_{1;m}:m\in\ZZ\}$  is given by
\begin{equation}\label{Eq:Y1mGeneral}
y_{1;m}=\frac{\mathbf{y}^{\mathbf{d}(x_{m+3})}}{F_{m+1}|_{\PP}(\mathbf{y})F_{m+2}|_{\PP}(\mathbf{y})}.
\end{equation} 
\item
The family $\{y_{2;m}:\, m\in\ZZ\}$ is given by: 
\begin{equation}\label{eq:Y20Y2-1}
\begin{array}{ccc}
y_{2;1}=y_{2},&y_{2;0}=(y_{3}\oplus1)y_{1},&y_{2;-1}=(y_{2}y_{3}\oplus y_{2}\oplus1)\frac{1}{y_{3}}.
\end{array}
\end{equation}
For every $m\geq1$
\begin{eqnarray}\label{eq:Y22m}
       y_{2;2m}&=&\frac{F_{2m}|_{\PP}(\mathbf{y})}{F_{2m+2}|_{\PP}(\mathbf{y})}y_{1}y_{3}\\\label{eq:Y2-2m}
       y_{2;-2m}&=&\frac{F_{-2m}|_{\PP}(\mathbf{y})}{F_{-2m+2}|_{\PP}(\mathbf{y})}\frac{1}{y_{2}}\\\label{eq:Y22m+1}
       y_{2;2m+1}&=&\frac{F_{2m+1}|_{\PP}(\mathbf{y})}{F_{2m+3}|_{\PP}(\mathbf{y})}y_{2}\\\label{eq:Y2-2m-1}
       y_{2;-2m-1}&=&\frac{F_{-2m-1}|_{\PP}(\mathbf{y})}{F_{-2m+1}|_{\PP}(\mathbf{y})}\frac{1}{y_{1}y_{3}}       
\end{eqnarray}
\item
The family $\{y_{3;m}|\,m \in\ZZ\}$ is given by $y_{3;m+1}=1/y_{1;m}$.
\item
The families $\{y_{i;m}^{cyc}|\,m \in\ZZ\}$ for $i=1, 2, 3$ are given by
\begin{equation}\label{eq:Yic}
\begin{array}{ccc}
y_{1;m}^{cyc}= \frac{y_{1;m}}{y_{2;m}  \oplus 1};& y_{2;m}^{cyc} = 1/y_{2;m};& y_{3;m+2}^{cyc} = 1/y_{1;m}^{cyc}
\end{array}
\end{equation}
\end{itemize}
\end{proposition}
\begin{proof}[Proof of proposition~\ref{Prop:CoefficientsRecurrence}] Formulas \eqref{eq:Y20Y2-1}--\eqref{eq:Yic} follows directly from \eqref{Eq:Y1mGeneral}. We hence prove \eqref{Eq:Y1mGeneral}. 

By \eqref{CoeffMutationDirection1}  it follows that the family $\{y_{1;m} : m \in\ZZ\}$ is the sequence of elements of $\PP$ uniquely determined by the
initial data: $y_{1;m} = 1/y_{m+3}$ if $m=0,-1,-2$, $y_{1;1}=y_1$,  $y_{1;-3}=y_3$ together with the recurrence relations
\begin{equation}\label{eq:Y1m}
y_{1;m}y_{1;m+3} =\frac{y_{1;m+2}y_{1;m+1}}{(y_{1;m+2}\oplus1)(y_{1;m+1}\oplus1)}
\end{equation}                                      

Let us first suppose $\PP=Trop(y_{1},y_{2},y_{3})$. By induction on $m\geq1$ and $m+3\leq1$ and by the relation \eqref{eq:RecurrenceDenominator} one proves that   the solution of \eqref{eq:Y1m} is the following: for every $m\in\ZZ$
\begin{equation}\label{eq:RecurrenceSolution}        
y_{1;m} = \mathbf{y}^{\mathbf{d}(x_{m+3})}
\end{equation}
where $\mathbf{d}(x_{m+3})$ is the denominator vector the cluster variable $x_{m+3}$ in the cluster $\{x_{1},x_{2},x_{3}\}$ given in lemma~\ref{Lemma:DenominatorVectors}. 

Let us assume now that $\PP$ is any semifield. Recall that the $F$--polynomial $F_{m}$ of the cluster variable $x_{m}$ in the cluster $\{x_{1},x_{2},x_{3}\}$ is obtained from the Laurent expansion of $x_{m}$ in $\{x_{1},x_{2},x_{3}\}$ by specializing $x_{1}=x_{2}=x_{3}=1$. In view of \eqref{ExchRelPrincipal}, the family $\{F_{m}:\,m\in\ZZ\}$ is hence recursively defined by the initial data: $F_{1}=F_{2}=F_{3}=1$ 
\begin{equation}\label{F0F-1F-2}
\begin{array}{c}
F_{0}(\mathbf{y})=y_{3}+1\qquad
F_{-1}(\mathbf{y})=y_{2}y_{3}+y_{2}+1\\
F_{-2}(\mathbf{y})=y_{1}y_{2}y_{3}^{2}+2y_{1}y_{2}y_{3}+y_{1}y_{3}+y_{1}y_{2}+y_{1}+1
\end{array}
\end{equation}
together with the recurrence relations for $m\geq1$ and $m\leq-3$:
\begin{equation}\label{Eq:FPolyRecursion}
F_{m}(\mathbf{y})F_{m+3}(\mathbf{y})=F_{m+1}(\mathbf{y})F_{m+2}(\mathbf{y})+\mathbf{y}^{\mathbf{d}(x_{m+3})}
\end{equation}
We use the notation $F(\mathbf{y}):=F(y_{1},y_{2},y_{3})$.  Then an easy induction gives the desired \eqref{Eq:Y1mGeneral}.
\end{proof}

As a corollary of the previous proposition we get  relations between principal cluster variables: we recall that given a cluster variable $s$ the corresponding principal cluster variable in the cluster $\{x_{1},x_{2},x_{3}\}$ is the element $S:=F_{s}|_{\PP}(y_{1},y_{2},y_{3})s$. 

\begin{corollary} The principal cluster variables in the cluster $\{x_{1},x_{2},x_{3}\}$ of $\myAA_{\PP}$ satisfy the following relations:
\begin{equation}\label{ExchRelPrincipalX0}
\begin{array}{ccc}
X_{0}X_{3}\!=\!y_{3}X_{1}X_{2}\!+\!1;\!&\!
X_{-1}X_{2}\!=\!y_{2}X_{0}X_{1}\!+\!1;\!&\!
X_{-2}X_{1}\!=\!y_{1}X_{-1}X_{0}\!+\!1.
\end{array}
\end{equation}
\begin{equation}\label{ExchRelPrincipal}
X_{m}X_{m+3}=X_{m+1}X_{m+2}+\mathbf{y}^{\mathbf{d}(x_{m+3})}\; \text{ for }m\geq1\text{ and }m\leq-3
\end{equation}
\begin{equation}\label{ExRelWX2mExplicit}
WX_{2m}=\left\{\begin{array}{cc} y_{2}X_{2m-1}+X_{2m+1}&\text{ if }m\geq1\\
X_{-1}+y_{3}X_{1}&\text{ if }m=0\\
X_{2m-1}+y_{1}y_{3}X_{2m+1}&\text{ if }m\leq-1
\end{array}
\right.
\end{equation}
\begin{equation}\label{ExRelZX2m+1Explicit}
ZX_{2m+1}=\left\{\begin{array}{cc} y_{1}y_{3}X_{2m}+X_{2m+2}&\text{ if }m\geq1\\
y_{1}X_{0}+X_{2}&\text{ if }m=0\\
X_{2m}+y_{2}X_{2m+2}&\text{ if }m\leq-1
\end{array}
\right.
\end{equation}
\begin{eqnarray}
\label{ExchRel:X2m-2X2m+2Principal} 
X_{2m-2} X_{2m+2}&=&X_{2m}^{2} + \mathbf{y}^{\mathbf{d}(x_{2m+1})}Z\\
\label{ExchRel:X2m-1X2m+3Principal} 
X_{2m-1} X_{2m+3}&=&X_{2m+1}^{2}+ \mathbf{y}^{\mathbf{d}(x_{2m+2})}W
\end{eqnarray}
In particular, if $\PP=\text{Trop}(y_{1},y_{2},y_{3})$, these are precisely the exchange relations of the cluster algebra $\myAA_{\PP}=\myAA_{\bullet}(\Sigma)$ with principal coefficients at the seed $\Sigma=\Sigma_1$.
\end{corollary}

\section{Proof of theorem~\ref{Prop:BijectionDenominators}}\label{Sec:ProofPropBijection}

In this section we prove that the denominator vector map $\mathbf{d}: \myAA\rightarrow Q$ in the initial cluster $\{x_{1},x_{2},x_{3}\}$ restricts to a bijection between $\BB$ and the root lattice $Q$ of type $A_2^{(1)}$. 

Clearly the denominator vector of a cluster monomial $s_{1}^{a}s_{2}^{b}s_{3}^{c}$ is $a\mathbf{d}(s_{1})+b\mathbf{d}(s_{2} ) + c\mathbf{d}(s_{3})$. We also know that
$\mathbf{d}(u_{1})=\mathbf{d}(w) + \mathbf{d}(z)$ and $\mathbf{d}(u_{n})=n\mathbf{d}(u_{1})$. Hence the cone $\mathcal{C}_{Reg}=\ZZ_{\geq0} \mathbf{d}(w) + \ZZ_{\geq0} \mathbf{d}(z)$ generated by $\mathbf{d}(w)$ and $\mathbf{d}(z)$ is in bijection with the set $\{\mathbf{d}(u_{n}w^{k}), \mathbf{d}(u_{n}z^{k})|\, n, k \geq 0\}$.  To complete the proof of theorem~\ref{Prop:BijectionDenominators} it is enough to show the following:
\begin{eqnarray}\label{eq:Unimodular}&&
\text{For every  cluster $\{s_{1},s_{2},s_{3}\}$, the vectors $\mathbf{d}(s_1 )$, $\mathbf{d}(s_2)$ and }\\\nonumber&&
\text{$\mathbf{d}(s_3)$ form a $\ZZ$--basis of $Q$}.\\\label{eq:RealRoots}&&
\text{For every element $v$ of $Q$ which does not lie in the interior of $\myCC_{Reg}$}\\\nonumber&&
\text{there exists a unique cluster $\{s_{1},s_{2},s_{3}\}$ such that $v$ belongs to}\\\nonumber&&
\text{the cone $\myCC_{\{s_1 ,s_2 ,s_3 \}}:=\ZZ_{\geq0} \mathbf{d}(s_1 ) + \ZZ_{\geq0} \mathbf{d}(s_2)+\ZZ_{\geq0}\mathbf{d}(s_3)$. }
\end{eqnarray}
This is done in the subsequent sections.

\subsubsection*{Proof of \eqref{eq:Unimodular}}
By the explicit formulas of denominator vectors of cluster variables given in lemma~\ref{Lemma:DenominatorVectors}  one checks directly that the absolute value of the determinant $|\text{det}(\mathbf{d}(s_{1}),\mathbf{d}(s_{2}),\mathbf{d}(s_{3}))|=1$ for every cluster $\{s_{1},s_{2},s_{3}\}$.
 
\subsubsection*{Proof of \eqref{eq:RealRoots}}
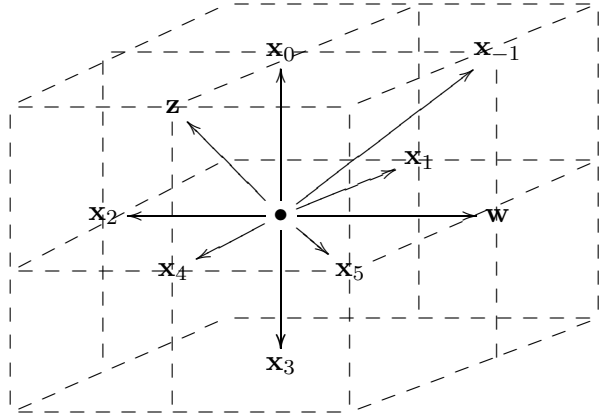
\begin{figure}[htbp]
\begin{center}
\xymatrix@C=8pt@R=8pt{
&&&& &&&&\ar@{--}[rrrrrr]&               &&\ar@{--}[dddddd]&&&\ar@{--}[dddddd]\\
&&&&
&&\ar@{--}[dddddd]\ar@{--}[rrrrrr]&&& \mathbf{x}_{0}              &&&\mathbf{x}_{-1}\ar@{--}[dddddd]&&\\
&&&&
\ar@{--}[rrrruu]\ar@{--}[dddddd]\ar@{--}[rrrrrr]&&&\mathbf{z}\ar@{--}[rrrruu]\ar@{--}[dddddd]&&               &\ar@{--}[rrrruu]\ar@{--}[dddddd]&&&&\\
&&&&
&&&&\ar@{--}[rrrrrr]&               &&\mathbf{x}_{1}&&&\\
&&&&
&&\mathbf{x}_{2}&&&\bullet\ar[rrr]\ar[lll]\ar[uuu]\ar[ddd]\ar[rrruuu]\ar[rru]\ar[dr]\ar[dll]\ar[ddd]\ar[lluu]&&&\mathbf{w}&&\\
&&&&
\ar@{--}[rrrruu]\ar@{--}[rrrrrr]&&&\mathbf{x}_{4}&&       &\mathbf{x}_{5}\ar@{--}[rrrruu]&&&&\\
&&&&
&&&&\ar@{--}[rrrrrr]&          &&&&&\\
&&&&
&&&&& \mathbf{x}_{3} &&&&&\\
&&&&
\ar@{--}[rrrrrr]\ar@{--}[rrrruu]&&&&&          &\ar@{--}[rrrruu]&&&&
}
\caption{Denominator vectors of cluster variables having at least
     one coordinate equal to zero. We wrote $x_{m}$ for $\mathbf{d}(x_{m} )$. The clusters involving here form a fan whose union is $Q \setminus Q_+$}
\label{Fig:InitialClusterVariables}
\end{center}
\end{figure}

We consider the basis of simple roots $\alpha_{1}$, $\alpha_{2}$, $\alpha_{3}$ of $Q$ and the corresponding coordinate system $(e_{1}, e_{2}, e_{3})$. In corollary~\ref{Cor:SignCoherence} it is shown that given a cluster $\myCC=\{s_1 , s_2 , s_3 \}$ the set of  corresponding denominator vectors $\{\mathbf{d}(s_1 ), \mathbf{d}(s_2 ), \mathbf{d}(s_3 
)\}$ are sign--coherent.  In particular if the initial cluster variable $x_{i}$ lies in the cluster $\myCC$, then the $i$--th coordinates of the other two elements of $\myCC$ are zero. By lemma~\ref{Lemma:DenominatorVectors} also the opposite holds: if the $i$--th coordinate of the denominator vector of  a cluster variable is zero, then it lies in the same cluster as $x_{i}$. There are hence precisely nine cluster variables with this property and the corresponding denominator vectors are shown in figure~\ref{Fig:InitialClusterVariables}. The corresponding cones $\myCC_{\{s_1 ,s_2 ,s_3\}}$ satisfy property \eqref{eq:RealRoots}, i.e. they do not overlap themselves. Moreover their union is the whole lattice except the interior of the positive octant
$Q_{+} = \ZZ_{\geq0} \mathbf{d}(x_{4} ) +\ZZ_{\geq 0} \mathbf{d}(w) +\ZZ_{\geq 0} \mathbf{d}(x_{0})$.

We hence consider the denominator vectors of cluster variables contained in $Q_+$. We suggest to  use figure~\ref{Fig:ClusterTriangulation} to visualize the situation.
By using lemma~\ref{Lemma:DenominatorVectors}, we notice that there are four affine lines in $Q_\mathbb{R}$ which contain the denominator vectors of all the cluster variables different from $x_2$. They contain respectively ``negative odd'', ``positive odd'', 
``negative even'' and ``positive even'' cluster variables: they are
$$
\begin{array}{cc}
\ell^{-}_{\text{odd}}:=\left\{\begin{array}{l}e_{2}=e_{3}\\ e_{1}=e_{2}-1\end{array}\right.;&
\ell^{+}_{\text{odd}}:=\left\{\begin{array}{l}e_{1}=e_{2}\\ e_{3}=e_{2}-1\end{array}\right.;\\
\ell^{-}_{\text{even}}:=\left\{\begin{array}{l}e_{1}=e_{2}\\ e_{3}=e_{2}+1\end{array}\right.;&
\ell^{+}_{\text{even}}:=\left\{\begin{array}{l}e_{2}=e_{3}\\ e_{1}=e_{2}+1\end{array}\right..
\end{array}
$$
We define the two-dimensional subspaces $P$ and $T$ of $Q_{\mathbb{R}}$ containing respectively both $\ell^{+}_{\text{odd}}$ and $\ell^{-}_{\text{even}}$ and both $\ell^{-}_{\text{odd}}$ and $\ell^{+}_{\text{even}}$ of equation: $P:= \{e_{1} = e_{2} \}$ and  $T:= \{e_{2} = e_{3} \}$.

 Let $\myCC_{P}$ be the (open) cone
inside $P \cap Q_+$ defined by $\myCC_{P}:=\{0 < e_3 < e_1\}\cup\{0\}$. By \eqref{eq:D(XmPositive)}, $\mathbf{d}(x_{2n+1}) \in\myCC_{P}$
for every $n\geq 2$. The vectors $v_{1}:=(1, 1, 0)^{t}=\mathbf{d}(x_{5})$ and $v_{2}= (0, 0, 1)^{t}= \mathbf{d}(x_{0})$ form a
$\ZZ$--basis of $P$ such that $\myCC_{P}$ is contained in $\ZZ_{\geq 0} v_{1} +\ZZ_{\geq 0} (v_1 + v_2 )$.  
In this basis $\mathbf{d}(x_{2n+1})=a_{n1} v_{1} + a_{n2} v_{2}$ where $a_{n1}=n-1$ and $a_{n2}=n-2$. The sequence $a_{n2}/a_{n1}$ is strictly increasing. It has limit
$\lim_{n\to\infty}\frac{a_{n2}}{a_{n1}}=1.$   We conclude that 
\begin{equation}\label{Eq:Cp=Union}
\myCC_{P}=\bigcup_{n\geq2}\myCC_{\{x_{2n+1},x_{2n+3}\}}
\end{equation}
(here and in the sequel we set $\myCC_{\{s_1,\cdots,s_k\}}:=\ZZ_{\geq0} \mathbf{d}(s_1) +\cdots+\ZZ_{\geq 0} \mathbf{d}(s_k)$). Moreover the interior of two different cones in the right hand side are disjoint. In particular we have
$$
\bigcup_{n\geq2}\myCC_{\{x_{2n+1},w,x_{2n+3}\}}=\ZZ_{\geq 0} \mathbf{d}(w) + \myCC_{P}
$$                               
and the cones in the left hand side have no common interior points. 

Similarly let $\myCC_{T}$ be the (open) cone inside $T \cap Q_+$ defined by $\myCC_{T}=\{0 < e_2 < e_1 \} \cup \{0\}$.  By \eqref{eq:D(XmPositive)}, $\mathbf{d}(x_{2n}) \in\myCC_{T}$ for every $n \geq 2$. The vectors $w_{1}=(1, 0, 0)^{t}=\mathbf{d}(x_{4})$ and $w_{2}=(0, 1, 1)^{t}=\mathbf{d}(x_{-1})$ form a $\ZZ$--basis of T such that $\myCC_{T}$ is contained in $\ZZ_{\geq 0} w_{1} +\ZZ_{\geq0}(w_{1}+w_{2})$. In this basis $\mathbf{d}(x_{2n})=b_{n1}w_{1}+b_{n2}w_{2}$ with $b_{n1}=n-1$ and
$b_{n2}=n-2$. 
The strictly increasing sequence $\{b_{n2} /b_{n1} \}$ has limit $1$ for $n\to\infty$. We conclude that 
\begin{equation}\label{Eq:Ct=Union}
\myCC_{T}=\bigcup_{n\geq2}\myCC_{\{x_{2n},x_{2n+2}\}}
\end{equation}
and the interiors of the cones in the right hand side are mutually disjoint. In particular we have
$$
\bigcup_{n\geq2}\myCC_{\{x_{2n} ,z,x_{2n+2}\}} =\ZZ_{\geq 0} \mathbf{d}(z)+\myCC_{T}.
$$ 
and the cones in the left hand side have no common interior points. 

We now prove that
$$
\myCC_{P}+\myCC_{T}=\bigcup_{m\geq4} \myCC_{\{x_{m} ,x_{m+1} ,x_{m+2}\}}
$$
and that the interiors of two different cones in the right hand side are disjoint. By definition $\myCC_{P}+\myCC_{T}$ is contained in $\ZZ\mathbf{d}(u_1)+\ZZ\mathbf{d}(x_4)+\ZZ\mathbf{d}(x_5)$ and consists of all integer vectors $v=(l,m,n)^t$ such that $l\geq m\geq n\geq0$ not all equals (see figure~\ref{Fig:ClusterTriangulation}). The expansion of such vector in the basis $\{\mathbf{d}(u_1),\mathbf{d}(x_4),\mathbf{d}(x_5)\}$ is:
\begin{equation}\label{Eq:VAlgorithm}
v=n\mathbf{d}(u_1)+(l-m)\mathbf{d}(x_4)+(m-n)\mathbf{d}(x_5).
\end{equation}
We are going to provide an algorithm which gives the precise expression of $v$ as an element of precisely one cone $\myCC_{\{x_m,x_{m+1},x_{m+2}\}}$ for some $m\geq4$. This is done in two steps. The first step provides the explicit  expression of an element of  $\myCC_{\{x_4,x_5,x_6,x_7\}}$ (this is the  quadrilateral at the bottom of figure~\ref{Fig:ClusterTriangulation}) in one of the two cones 
$\myCC_{\{x_4,x_5,x_6\}}$ or $\myCC_{\{x_5,x_6,x_7\}}$. The second step reduces to the case in which $v$ belongs to the cone  $\myCC_{\{x_4,x_5,x_6,x_7\}}$ .  

\textbf{First step}: the following lemma gives arithmetic conditions on $v=(l,m,n)^t$ to belong to $\myCC_{\{x_4,x_5,x_6,x_7\}}$.

\begin{lemma}\label{lemma:VInQuadrilateral}
An element $v=(l,m,n)^t\in\myCC_{P}+\myCC_{T}$ belongs to $\myCC_{\{x_4,x_5,x_6,x_7\}}$ if and only if $l\geq2n$. In this case  if $n\leq l-m$ then 
\begin{equation}\label{Eq:Alg1}
v=n\mathbf{d}(x_6)+(l-m-n)\mathbf{d}(x_4)+(m-n)\mathbf{d}(x_5)
\end{equation}
and $v\in\myCC_{\{x_4,x_5,x_6\}}$.
If $n\geq l-m$, then
\begin{equation}\label{Eq:Alg2}
v=(l-2n)\mathbf{d}(x_5)+(l-m)\mathbf{d}(x_6)+(n-l+m)\mathbf{d}(x_7)
\end{equation}
and $v\in\myCC_{\{x_5,x_6,x_7\}}$. 
\end{lemma}
\begin{proof}
If $l<2n$ then $l-n=(l-m)+(m-n)<n$ and there exists $k>0$ such that $n=(l-m)+(m-n)+k$. By \eqref{Eq:VAlgorithm}  we have 
$$
v=(k-1)\mathbf{d}(u_1)+(l-m)\mathbf{d}(x_8)+(m-n)\mathbf{d}(x_7).
$$
which, in view of corollary~\ref{Cor:D(x)+D(u)},  does not belong to $\myCC_{\{x_4,x_5,x_6,x_7\}}$. On the other hand let us assume $l\geq 2n$. Then $(l-m)+(m-n)>n$ and there are two cases: either $n\leq l-m$ or $n>l-m$. In the first case $v$ has the expansion \eqref{Eq:Alg1} and hence $v\in\myCC_{\{x_4,x_5,x_6\}}$; in the second case  $v$ has the expansion \eqref{Eq:Alg2} and hence $v\in\myCC_{\{x_5,x_6,x_7\}}$. 
\end{proof}
\begin{corollary}
The interiors of the two cones $\myCC_{\{x_4,x_5,x_6\}}$ and  $\myCC_{\{x_5,x_6,x_7\}}$ are disjoint 
\end{corollary}
\begin{proof}
Let $v=(l,m,n)^t$ be an element of both  $\myCC_{\{x_4,x_5,x_6\}}$ and  $\myCC_{\{x_5,x_6,x_7\}}$. Then, by lemma \eqref{lemma:VInQuadrilateral}, $l\geq 2n$, $n=l-m$ and $v=n\mathbf{d}(x_6)+(m-n)\mathbf{d}(x_5)$. It hence follows that $v$ does not belong to the interior of the two cones.  
\end{proof}
By lemma \eqref{lemma:VInQuadrilateral} we notice that an element $v$ of $\myCC_{\{x_4,x_5,x_6,x_7\}}$ has  the form:
\begin{equation}\label{Eq:VAlphaAlgorithm}
v=\alpha_1\mathbf{d}(x_4)+\alpha_2\mathbf{d}(x_5)+\alpha_3\mathbf{d}(x_6)+\alpha_4\mathbf{d}(x_7)
\end{equation}
for some non negative integers $\alpha_1$, $\alpha_2$, $\alpha_3$ and $\alpha_4$ such that $\alpha_1\alpha_4=0$ and $\alpha_1+\alpha_2+\alpha_3+\alpha_4=l-n$.

\textbf{Second step}: let $v=(l,m,n)^t$ be not an element of $\myCC_{\{x_4,x_5,x_6,x_7\}}$. By lemma~\ref{lemma:VInQuadrilateral} this implies that $l-n<n$. We divide $n$ by $l-n$ and we find $k$ and $r$ such that $0\leq r< (l-n)$ and $n=k(l-n)+r$.
We consider the vector
$$
v':=v-k(l-n)\mathbf{d}(u_1)=(l-k(l-n),m-k(l-n),n-k(l-n))^t.
$$
By using lemma~\ref{lemma:VInQuadrilateral}, one checks easily that $v'$ belongs to $\myCC_{\{x_4,x_5,x_6,x_7\}}$. Then by the first step $v'$ has the form \eqref{Eq:VAlphaAlgorithm}.  Then we get
\begin{eqnarray}\nonumber
v&=&v'+k(l-n)\mathbf{d}(u_1)\\\nonumber
  &=&\alpha_1(\mathbf{d}(x_4)+k\mathbf{d}(u_1))+\alpha_2(\mathbf{d}(x_5)+k\mathbf{d}(u_1))+\\\nonumber
&&+\alpha_3(\mathbf{d}(x_6)+k\mathbf{d}(u_1))+\alpha_4(\mathbf{d}(x_7)+k\mathbf{d}(u_1))\\\nonumber
  &=&\alpha_1\mathbf{d}(x_{4+2k})+\alpha_2\mathbf{d}(x_{5+2k})+\alpha_3\mathbf{d}(x_{6+2k})+\alpha_4\mathbf{d}(x_{7+2k})
\end{eqnarray}
which is the desired expansion. In particular we have that $v=(l,m,n)^t$ lies in $\myCC_{\{x_{4+k},x_{5+k},x_{6+k},x_{7+k}\}}$ if and only if $n=k(l-n)+r$
 for some $0\leq r<(l-n)$; in this case  $v$ belongs to exactly one cone, either $\myCC_{\{x_{4+k},x_{5+k},x_{6+k}\}}$ or $\myCC_{\{x_{5+k},x_{6+k},x_{7+k}\}}$.

By now in figure~\ref{Fig:ClusterTriangulation} we have obtained all the elements of the cone $\ZZ_{\geq0}\mathbf{d}(z)+\ZZ_{\geq0}\mathbf{d}(w)+\ZZ_{\geq0}\mathbf{d}(x_{4})=\myCC_{Reg}+\ZZ_{\geq0}\mathbf{d}(x_{4})$. We consider the orthogonal reflection $r_{Reg}$ with respect to the regular cone $\myCC_{Reg}$: this is the $\ZZ$--linear isomorphism of $Q$ which exchanges the first coordinate with the third one. In particular it fixes $\myCC_{Reg}$ pointwise. By remark~\ref{Rem:Symmetries}, $r_{Reg}$ sends $\mathbf{d}(x_{m})$ to $\mathbf{d}(x_{-m+4})$ for $m \geq 4$ and hence induces a bijection between $\myCC_{Reg}+\ZZ_{\geq 0} \mathbf{d}(x_{4})$ and $\myCC_{Reg}+\ZZ_{\geq 0} \mathbf{d}(x_{0})$. This concludes the proof of theorem~\ref{Prop:BijectionDenominators}.

\begin{remark}
The proof of theorem~\ref{Prop:BijectionDenominators} contains an algorithm to compute the ``virtual'' canonical decomposition of every element of $Q$. There are several more effective algorithms than this in much more generality  (see \cite{KacInfinite}, \cite{KacInfinite2}, \cite{SchofieldGeneric}, \cite{DWSubrepGen}, \cite{DW}). 
\end{remark}

\section{Proof of theorem~\ref{Thm:ExplicitExpressions}}\label{Sec:ExplicitExpressions}

The proof of theorem~\ref{Thm:ExplicitExpressions} is based on proposition~\ref{Prop:Homogeneity}: we find the explicit expression of the $F$--polynomial and of the $\mathbf{g}$--vector of every element of $\BB$ in every cluster of $\myAA_\PP$. Once again it is sufficient to consider only the two clusters $\{x_{1},x_{2},x_{3}\}$ and $\{x_{1},w,x_{3}\}$ (see remark~\ref{Rem:Symmetries}). The $F$--polynomials (resp. the $\mathbf{g}$--vectors) in the cluster $\{x_{1},x_{2},x_{3}\}$ and $\{x_{1},w,x_{3}\}$ are given respectively in proposition~\ref{Prop:FPolyInitial} (resp. proposition~\ref{Prop:GvectorsInitial}) and proposition~\ref{Prop:FpolyCyclic} (resp. proposition~\ref{prop:GVectorsCyclic}). 

\subsection{Proof of proposition~\ref{Prop:Homogeneity}}\label{subsec: Homogeneity} 

By \cite[corollary~6.3]{FZIV} the expansion of all cluster variables and hence of all cluster monomials in every cluster has the form \eqref{Eq:b=FbgbInEveryCluster}. In section ~\ref{Sec:FPolyGvector} it is shown that the $F$--polynomial of every cluster variable in every cluster has positive integer coefficients. It remains to deal with the $u_n$'s. We prove that for every $n\geq1$, $u_n$ has the form \eqref{Eq:b=Fbgb} in both the clusters $\{x_1,x_2,x_3\}$ and $\{x_1,w,x_3\}$; by the symmetry of the exchange relations this implies that $u_n$ has the form \eqref{Eq:b=Fbgb} in every cluster of $\myAA$.  

Let $\PP$ be a semifield. Let $\mathcal{M}$ be the set of all the elements $b$ of $\myFF_{\PP}$  that can be written in the form 
\begin{equation}\label{Eq:b=fyHat}
b=F_{b}(\hat{y}_{1},\hat{y}_{2},\hat{y}_{3})\mathbf{x}^{\mathbf{g}_{b}}
\end{equation} 
where $F_b$ is a polynomial with integer coefficients and
\begin{equation}\label{Eq:DefinitionYhatInitial}
\begin{array}{ccc}
\hat{y}_{1}:=\frac{y_{1}}{x_{2}x_{3}}=y_{1}\mathbf{x}^{\mathbf{h}_{1}}&\hat{y}_{2}:=\frac{y_{2}x_{1}}{x_{3}}=y_{2}\mathbf{x}^{\mathbf{h}_{2}}&\hat{y}_{3}:=y_{3}x_{1}x_{2}=y_{3}\mathbf{x}^{\mathbf{h}_{3}}
\end{array}
\end{equation}
Both the principal cluster variables $W$ and $Z$ belong to $\mathcal{M}$ and by  \eqref{Eq:DefinitionWZ} their $F$--polynomials and $\mathbf{g}$--vectors are respectively
\begin{equation}\label{eq:FwInitial}
 \begin{array}{cc} 
 F_{w}(\mathbf{y})=y_2 + 1,&\mathbf{g}_{w}=\left[\def\objectstyle{\scriptstyle}\def\labelstyle{\scriptstyle}\vcenter{\xymatrix@R=0pt@C=0pt{0\\-1\\1}}\right]\\
   F_{z} (\mathbf{y})= y_1 y_3 + y_1 + 1&\mathbf{g}_{z}=\left[\def\objectstyle{\scriptstyle}\def\labelstyle{\scriptstyle}\vcenter{\xymatrix@R=0pt@C=0pt{-1\\1\\0}}\right]
 \end{array}
\end{equation}
By definition~\ref{def:Un} $u_{1}=ZW-y_{1}y_{3}-y_{2}$. By inverting the equalities in \eqref{Eq:DefinitionYhatInitial} we get:
\begin{equation}\label{eq:U1Yhat}
u_{1}=[F_{z}(\hat{y}_{1},\hat{y}_{2},\hat{y}_{3})F_{w}(\hat{y}_{1},\hat{y}_{2},\hat{y}_{3})-\hat{y}_{1}\hat{y}_{3}-\hat{y}_{2}]\frac{x_{3}}{x_{1}}
\end{equation}
and hence $u_{1}$ belongs to $\mathcal{M}$ and its $F$--polynomial and its $\mathbf{g}$--vector are respectively
\begin{equation}\label{Eq:FGUn1}
\begin{array}{cc}
F_{u_{1}}(y_{1},y_{2},y_{3})=y_{1}y_{2}y_{3}+y_{1}y_{2}+y_{1}+1&
\mathbf{g}_{u_{1}}=\left[\def\objectstyle{\scriptstyle}\def\labelstyle{\scriptstyle}\vcenter{\xymatrix@R=0pt@C=0pt{-1\\0\\1}}\right].
\end{array}
\end{equation} 
Similarly, by definition~\ref{def:Un} and induction on $n$, we get:
\begin{eqnarray}\nonumber
u_{2}&=&[F_{u_{1}}(\hat{y}_{1},\hat{y}_{2},\hat{y}_{3})^{2}-2\hat{y}_{1}\hat{y}_{2}\hat{y}_{3}](\frac{x_{3}}{x_{1}})^{2}\\\label{Eq:Un=FGInitial}
u_{n+1}&=&[F_{u_{1}}F_{u_{n}}(\hat{y}_{1},\hat{y}_{2},\hat{y}_{3})-\hat{y}_{1}\hat{y}_{2}\hat{y}_{3}F_{u_{n-1}}(\hat{y}_{1},\hat{y}_{2},\hat{y}_{3})](\frac{x_{3}}{x_{1}})^{n+1}
\end{eqnarray}
We hence have that the $u_{n}$'s belong to $\mathcal{M}$ and their $F$--polynomial $F_{u_n}$ satisfies the initial condition \eqref{Eq:FGUn1} together with the recurrence relations for $n\geq2$:
\begin{eqnarray}\label{EqFU2Recursion}
F_{u_{2}}(y_1,y_2,y_3)&=&F_{u_{1}}^2(y_{1},y_{2},y_{3})-2y_{1}y_{2}y_{3}\\\label{EqFUNRecursion}
F_{u_{n+1}}(y_1,y_2,y_3)&=&F_{u_{1}}F_{u_{n}}(y_{1},y_{2},y_{3})-y_{1}y_{2}y_{3}F_{u_{n-1}}(y_{1},y_{2},y_{3})
\end{eqnarray}
It remains to prove $F_{u_n}$ is a polynomial with positive integer coefficients. This is done in proposition~\ref{Prop:FPolyInitial}.

We now prove that, for every $n\geq1$, $u_n$ has the form  \eqref{Eq:b=Fbgb}  in  the seed $\Sigma^{cyc}=\{H^{cyc},\{x_1,w,x_3\},\{\yy_1,\yy_2,\yy_3\}\}$ defined in \eqref{Def:CyclicSeed}. We introduce the elements $\hat{\yy}_{1}$, $\hat{\yy}_{2}$ and $\hat{\yy}_{3}$ of $\FF_{\PP}$ in analogy with \eqref{Eq:DefinitionYhatInitial} as follows:
\begin{equation}\label{Eq:DefinitionYhatCyclic}
\begin{array}{ccc}
\hat{\yy}_{1}:=\frac{\yy_{1}w}{x_{3}^{2}}&\hat{\yy}_{2}:=\frac{\yy_{2}x_{3}}{x_{1}}&\hat{\yy}_{3}:=\frac{\yy_{3}x_{1}^{2}}{w}
\end{array}
\end{equation}
We hence prove that for every $n\geq1$ there exists a polynomial $F^w_{u_n}(y_1,y_2,y_3)$ with positive integer coefficients and an integer vector $\mathbf{g}_{u_n}^w=(g_1,g_2,g_3)^t\in\ZZ^3$ such that the expansion of $u_n$ in $\Sigma^{cyc}$ is given by:
$$
u_n=F_{u_n}^w(\hat{\yy}_{1},\hat{\yy}_{2},\hat{\yy}_{3})x_1^{g_1}w^{g_2}x_3^{g_3}.
$$
By definition of the coefficient mutation \eqref{eq:CoefficientsMutation} in direction $2$, the coefficients $y_{1}$, $y_{2}$ and $y_{3}$ of the seed $\Sigma=\mu_2(\Sigma^{cyc})$ in the semifield $\PP$, are given by:
\begin{equation}
\begin{array}{ccc}
y_{1}=\frac{\yy_{1}\yy_{2}}{\yy_{2}\oplus1}&y_{2}=\frac{1}{\yy_{2}}&y_{3}=\yy_{3}(\yy_{2}\oplus1)
\end{array}
\end{equation}
 The following lemma shows that the elements $\{\hat{\yy}_{i}\}$ are  obtained from $\{\hat{y}_{i}\}$ by the mutation \eqref{eq:CoefficientsMutation} in direction $2$ (in the terminology of \cite{FZIV} this means that the families $\{\hat{y}_{i;\myCC}\}$ form a $Y$--pattern):
\begin{lemma}\label{lemma:YHatMut}
$$
\begin{array}{ccc}
\hat{y}_{1}=\frac{\hat{\yy}_{1}\hat{\yy}_{2}}{\hat{\yy}_{2}+1}&
\hat{y}_{2}=\frac{1}{\hat{\yy}_{2}}&
\hat{y}_{3}=\hat{\yy}_{3}(\hat{\yy}_{2}+1)
\end{array}
$$
\end{lemma}
\begin{proof}
By definition we have
$$
\begin{array}{cccc}
\hat{y}_{1}:=\frac{y_{1}}{x_{2}x_{3}}&\hat{y}_{2}=\frac{y_{2}x_{1}}{x_{3}}&\hat{y}_{3}:=y_{3}x_{1}x_{2}&x_{2}=\frac{x_{1}+\yy_{2}x_{3}}{(\yy_2\oplus1)w}
\end{array}
$$
and hence the proof follows by direct check.
\end{proof}
In view of lemma~\ref{lemma:YHatMut} the expansion of $u_n$ in $\Sigma^{cyc}$ is given by:
\begin{eqnarray}\nonumber
u_n&=&F_{u_n}(\hat{y}_1,\hat{y}_2,\hat{y}_3)(\frac{x_3}{x_1})^n\\\label{Eq:UnCyclic1}
&=&F_{u_n}(\frac{\hat{\yy}_{1}\hat{\yy}_{2}}{\hat{\yy}_{2}+1},\frac{1}{\hat{\yy}_{2}},\hat{\yy}_{3}(\hat{\yy}_{2}+1))(\frac{x_3}{x_1})^n
\end{eqnarray}
\begin{lemma}\label{lemma:FWPolynomial}
 For every $n\geq1$, $F_{u_n}(\frac{y_{1}y_{2}}{y_{2}+1},\frac{1}{y_{2}},y_{3}(y_{2}+1))\in\ZZ[y_1,y_2,y_3]$ 
\end{lemma}
\begin{proof}
By \eqref{Eq:FGUn1} the statement holds for $n=1$. By \eqref{EqFU2Recursion} and \eqref{EqFUNRecursion} an easy induction shows the  result. 
\end{proof}
\begin{definition}\label{def:Fw}
For every $n\geq1$ we define
\begin{equation}\label{Eq:FWUnFUn}
F_{u_n}^{w}(\yy_{1},\yy_{2},\yy_{3}):=F_{u_n}(\frac{\yy_{1}\yy_{2}}{\yy_{2}+1},\frac{1}{\yy_{2}},\yy_{3}(\yy_{2}+1)).
\end{equation}
We define the vector $\mathbf{g}_{u_n}^{w}:=(-n,0,n)^{t}\in\ZZ^{3}$
\end{definition}   
In view of lemma~\ref{lemma:FWPolynomial} and definition~\ref{def:Fw}, \eqref{Eq:UnCyclic1} is the desired expansion. Proposition~\ref{Prop:FpolyCyclic} provides the explicit formulas of $F_{u_n}^w$ which is hence a polynomial with positive integer coefficients. This concludes the proof of proposition~\ref{Prop:Homogeneity}.

\subsection{$F$--polynomials and $\mathbf{g}$--vectors of the elements of $\BB$}\label{Sec:FPolyGvector}
In this section we provides explicit formulas for the $F$--polynomials and the $\mathbf{g}$--vectors of every element of $\BB$ in every cluster of $\myAA$. 

\begin{proposition}\label{Prop:FPolyInitial}
 The $F$--polynomial $F_m$ of a cluster variable $x_m$ ($m\geq1$) in $\{x_1,x_2,x_3\}$ is: for every $m\geq 0$
  \begin{equation}\label{Eq:F2m+1}
   F_{2m+1}(\mathbf{y})=\sum_{\mathbf{e}}{e_{1}-1\choose e_{3}}{m-1-e_{2}\choose e_{1}-e_{2}}{m-1-e_{3}\choose e_2-e_{3}}\mathbf{y}^{\mathbf{e}}+1.
  \end{equation}
  \begin{equation}\label{Eq:F2m}
   F_{2m+2}(\mathbf{y})=\sum_{\mathbf{e}}{e_{1}-1\choose e_{3}}{m-e_2\choose e_1-e_2}{m-1-e_3\choose e_2- e_3}\mathbf{y}^{\mathbf{e}}+1.
  \end{equation}
\begin{equation}\label{Eq:F-(2m+1)}
   F_{-(2m+1)}(\mathbf{y})=\sum_{\mathbf{e}}{m-e_{3}\choose e_{1}-e_3}{e_{2}\choose e_{3}}{e_{1}+1\choose      e_{2}}\mathbf{y}^{\mathbf{e}}+y_1^my_2^{m+1}y_3^{m+1}.
  \end{equation}
  \begin{equation}\label{Eq:F-(2m)}
   F_{-2m}(\mathbf{y})=\sum_{\mathbf{e}}{m-e_3\choose e_1-e_{3}}{e_2+1\choose e_3}{e_1\choose e_2}\mathbf{y}^{\mathbf{e}}+y_1^my_2^my_3^{m+1}.
  \end{equation}
For every $n\geq 1$ the $F$--polynomial of $u_n$ is the following:
\begin{equation}\label{Eq:FUn}
F_{u_{n}}(\mathbf{y})=\mathbf{y}^{n\delta}+\sum_{\mathbf{e}}\textstyle{{e_{1}-e_{3}\choose e_{2}-e_{3}}}\left[\textstyle{{e_{1}-1\choose e_{3}}{n-e_{3}\choose n-e_{1}}+{e_{1}-1\choose e_{3}-1}{n-e_{3}-1\choose n-e_{1}}}\right]\mathbf{y}^{\mathbf{e}}+1.
\end{equation}
\end{proposition}

\begin{proof}
By \eqref{Eq:DefinitionWZ} the $F$--polynomial of $w$ and $z$ is respectively
$F_{w}(\mathbf{y})=y_2 + 1$, and $F_{z} (\mathbf{y})=y_1 y_3 + y_1 + 1$.
By \eqref{ExRelWX2mExplicit} and \eqref{ExRelZX2m+1Explicit}, the $F$--polynomials satisfy the following recurrence relations: for $m\geq 1$
\begin{eqnarray}\label{eq:FzF2m+1}
F_{2m+2}&=&F_{z}F_{2m+1}-y_{1}y_{3}F_{2m}\\
\label{eq:FwF2m}
F_{2m+1}&=&F_{w}F_{2m}-y_{2}F_{2m-1}
\end{eqnarray}
for which \eqref{Eq:F2m+1} and \eqref{Eq:F2m} follow by induction on $m\geq1$. By \eqref{ExRelWX2mExplicit} and \eqref{ExRelZX2m+1Explicit} we have that for every $m\geq1$
\begin{eqnarray}\label{eq:FzF-(2m+1)}
F_{-2m}&=&F_{z}F_{-(2m-1)}-y_{2}F_{-(2m-2)}\\
\label{eq:FwF-2m}
F_{-(2m+1)}&=&F_{w}F_{-2m}-y_1y_{3}F_{-(2m-1)}
\end{eqnarray}
from which \eqref{Eq:F-(2m+1)} and \eqref{Eq:F-(2m)} follow by induction on $m\geq1$.

In order to get \eqref{Eq:FUn} we proceed by induction on $n\geq1$. By direct check one verifies that the right--hand side of \eqref{Eq:FUn} satisfies the initial condition \eqref{Eq:FGUn1} together with the recurrence relations \eqref{EqFU2Recursion} and \eqref{EqFUNRecursion}.
\end{proof}

\begin{proposition}\label{Prop:GvectorsInitial}
The $\mathbf{g}$--vector $\mathbf{g}_m$ in $\{x_1,x_2,x_3\}$ of a cluster variable $x_m$ is given by:
for every $m\geq 0$:
 \begin{eqnarray}\label{eq:GvectorInitialPositive}
 \mathbf{g}_{2m+1}=\left[\!\!\def\objectstyle{\scriptscriptstyle}\def\labelstyle{\scriptscriptstyle}\vcenter{\xymatrix@R=0pt@C=0pt{1-m\\0\\m}}\! \!\right]&&
 \mathbf{g}_{2m+2}=\left[\!\!\def\objectstyle{\scriptscriptstyle}\def\labelstyle{\scriptscriptstyle}\vcenter{\xymatrix@R=0pt@C=0pt{-m\\1\\m}}\!\! \right]\\\label{eq:GvectorInitialNegative}
 \mathbf{g}_{-(2m+1)}=\left[\!\!\def\objectstyle{\scriptscriptstyle}\def\labelstyle{\scriptscriptstyle}\vcenter{\xymatrix@R=0pt@C=0pt{-m\\-1\\m}}\!\!\right]&&
 \mathbf{g}_{-2m}=\left[\!\!\def\objectstyle{\scriptscriptstyle}\def\labelstyle{\scriptscriptstyle}\vcenter{\xymatrix@R=0pt@C=0pt{-m\\0\\m-1}}\!\!\right]
 \end{eqnarray}
For every $n\geq1$ the $\mathbf{g}$--vector of $u_n$ is the following
\begin{equation}\label{Eq:Gun}
\mathbf{g}_{u_{n}}=\left[\def\objectstyle{\scriptstyle}\def\labelstyle{\scriptstyle}\vcenter{\xymatrix@R=0pt@C=0pt{-n\\0\\n}}\right]
\end{equation}
\end{proposition}

\begin{proof}
By the exchange relation \eqref{ExchRelPrincipal} the family $\{\mathbf{g}_{m}:\,m\geq1\}\subset\ZZ^{3}$ satisfies the initial conditions $\mathbf{g}_{i}=\mathbf{e}_{i}$, for $i=1,2,3$, together with the recurrence relations: for $m\geq1$
\begin{equation}\label{eq:RecurrenceGm}
\mathbf{g}_{m+3}+\mathbf{g}_{m}=\mathbf{g}_{m+1}+\mathbf{g}_{m+2}
\end{equation}
The proof is hence by induction on $m\geq1$ and $m\leq1$. The equality \eqref{Eq:Gun} follows from \eqref{Eq:Un=FGInitial}.
\end{proof}

\begin{proposition}\label{Prop:FpolyCyclic}
For every $m\in\ZZ$ the $F$--polynomial $F_m^w$ in $\{x_1,w,x_3\}$  of a cluster variable $x_m$ is given by: for $m\geq0$
\begin{eqnarray}\label{eq:F2m+1Cyclic}
F_{2m+1}^{w}(\mathbf{y})&=&\sum_{e_{1},e_{3}}\scriptstyle{{e_{1}-1\choose e_{3}}{m-1-e_{3}\choose e_{1}-e_3}}y_{1}^{e_{1}}y_{3}^{e_{3}}+1;\\\label{eq:F2m+2Cyclic}
F_{2m+2}^{w}(\mathbf{y})&=&\sum_{\mathbf{e}}\scriptstyle{{e_{1}-1\choose e_{3}}{m-1-e_{3}+e_{2}\choose e_{1}-e_{3}}{1\choose e_{2}}}\mathbf{y}^{\mathbf{e}}+y_{2}+1.
\\\label{eq:F-2m+1Cyclic}
F_{-(2m+1)}^{w}(\mathbf{y})&=&\sum_{e_{1},e_{3}}\scriptstyle{{m-e_{3}\choose e_{1}-e_3}{e_{1}+1\choose e_{3}}}y_{1}^{e_{1}}y_{3}^{e_{3}}+y_{1}^{m}y_{3}^{m+1};\\
\label{eq:F-2mCyclic}
F_{-2m}^{w}(\mathbf{y})&=&\sum_{\mathbf{e}}\scriptstyle{{m-e_{3}\choose e_{1}-e_3}{e_{1}+1-e_{2}\choose e_{3}-e_{2}}{1\choose e_{2}}}\mathbf{y}^{\mathbf{e}}+y_{1}^{m}y_{3}^{m+1}(y_{2}+1).
\end{eqnarray}
\begin{equation}\label{eq:FzCyclic}
F_{z}^{w}(\mathbf{y})=y_{1}y_{2}^{2}y_{3}+y_{1}y_{2}y_{3}+y_{1}y_{2}+y_{2}+1.
\end{equation}
For every $n\geq1$:
\begin{equation}\label{eq:FunCyclic}
F_{u_{n}}^{w}(\mathbf{y})=y_{1}^{n}y_{3}^{n}+\sum_{e_{1},e_{3}}\left[\scriptstyle{{n-e_{3}\choose n-e_{1}}{e_{1}-1\choose e_{3}}+{n-e_{3}-1\choose n-e_{1}}{e_{1}-1\choose e_{3}-1}}\right]y_{1}^{e_{1}}y_{3}^{e_{3}}+1.
\end{equation}
\end{proposition}

\begin{proof}
Let $\PP=\textrm{Trop}(\yy_1,\yy_2,\yy_3)$ and let $\myAA_\bullet(\Sigma^{cyc})$ be the cluster algebra with principal coefficients at the seed $\Sigma^{cyc}=\{H^{cyc},\{x_1,w,x_3\},\{\yy_1,\yy_2,\yy_3\}\}$ defined in \eqref{Def:CyclicSeed}. In view of proposition~\ref{Prop:Homogeneity} and lemma~\ref{lemma:YHatMut} the expansion of a cluster variable $x$ in this seed is given by:
\begin{eqnarray}\nonumber
x&=&\frac{F_{x}(\hat{y}_1,\hat{y}_2,\hat{y}_3)}{F_{x}|_\PP(y_1,y_2,y_3)}x_1^{g_1}x_2^{g_2}x_3^{g_3}\\\label{Eq:XW}
&=&\frac{F_{x}(\frac{\hat{\yy}_{1}\hat{\yy}_{2}}{\hat{\yy}_{2}+1},\frac{1}{\hat{\yy}_{2}},\hat{\yy}_{3}(\hat{\yy}_{2}+1))}{F_{x}|_\PP(\yy_{1}\yy_{2},\frac{1}{\yy_{2}},\yy_{3})}x_1^{g_1}(\frac{x_{1}+\yy_{2}x_{3}}{w})^{g_2}x_3^{g_3}
\end{eqnarray}
where $\mathbf{g}_x:=(g_1,g_2,g_3)^t$ is the $\mathbf{g}$--vector of $x$ in the cluster $\{x_1,x_2,x_3\}$.
In this expression we replace $x_1 = w = x_3 = 1$ and we get the relation:
\begin{equation}\label{eq:Fw=F}
F_{x}^{w}(\mathbf{\yy})=\frac{F_{x}(\frac{\yy_{1}\yy_{2}}{1+\yy_{2}},\frac{1}{\yy_{2}},\yy_{3}(1+\yy_{2}))}{F_{x}|_{\PP}(\yy_{1}\yy_{2},\frac{1}{\yy_{2}},\yy_{3})}\cdot(1+\yy_{2})^{g_{2}}
\end{equation}
By direct check, using proposition~\ref{Prop:FPolyInitial}, we get that 
\begin{equation}\label{eq:Fb|P}
F_{x}|_{\PP}(\yy_{1}\yy_{2},\frac{1}{\yy_{2}},\yy_{3})=\left\{
\begin{array}{cc}
\frac{1}{\yy_{2}}&\text{ if }x=x_{-(2m+1)}\;m\geq0,\textrm{ or }x=w\\
1&\text{ otherwise.}
\end{array}
\right.
\end{equation}
By the explicit formulas for the $\mathbf{g}$--vectors given in proposition~\ref{Prop:GvectorsInitial} and in \eqref{eq:FwInitial} we hence have that for every cluster variable $x$ the $F$--polynomials $F_x$ and $F_x^w$ are related by the following formula: 
\begin{equation}\label{eq:Fw=FII}
F_{x}^{w}(\mathbf{\yy})=\left\{\begin{array}{cc}
F_{x}(\frac{\yy_{1}\yy_{2}}{1+\yy_{2}},\frac{1}{\yy_{2}},\yy_{3}(1+\yy_{2}))\cdot\frac{\yy_{2}}{1+\yy_{2}}&\text{if }x=x_{-(2m+1)} \\&\textrm{or }x=w,\\
F_{x}(\frac{\yy_{1}\yy_{2}}{1+\yy_{2}},\frac{1}{\yy_{2}},\yy_{3}(1+\yy_{2}))\cdot(1+\yy_{2})^{\mathbf{g}_2}&\text{ otherwise.}
\end{array}\right.
\end{equation}
The proof of \eqref{eq:F2m+1Cyclic}--\eqref{eq:FzCyclic} now follows from proposition~\ref{Prop:FPolyInitial} by direct check.

Equation \eqref{eq:FunCyclic} follows from \eqref{Eq:FUn} by using \eqref{Eq:FWUnFUn}.
\end{proof}

\begin{corollary}\label{Cor:FbW(0)=1}
For every element $b$ of $\BB$, $F_{b}^{w}$ has constant term $1$.
\end{corollary}
\begin{remark}\label{Rem:Counterexample} 
We notice that
$$
F_{z}^{w}|_{\text{Trop}(\yy_{1},\yy_{2},\yy_{3})}(\frac{1}{\yy_{1}}.\frac{1}{\yy_{2}},\frac{1}{\yy_{3}})=\frac{1}{\yy_{1}\yy_{2}^{2}\yy_{3}}.
$$
In \cite[conjecture~7.17]{FZIV} it was expected the right-hand side to be $\mathbf{\yy}^{-\mathbf{d}^{w}(z)}=\frac{1}{\yy_{1}\yy_{2}\yy_{3}}$.
This counterexample appears also in \cite{BMR} and in \cite{FuKeller}.
\end{remark}

\begin{proposition}\label{prop:GVectorsCyclic}
For every $m\in\ZZ$ the $\mathbf{g}$-vector $\mathbf{g}_{m}^w$ of a cluster variable $x_m$ in the cluster $\{x_1,w,x_3\}$ is the following: for every $m\geq0$
\begin{eqnarray}\label{Eq:GPositiveW}
\mathbf{g}_{2m+1}^{w}=\left[\def\objectstyle{\scriptstyle}\def\labelstyle{\scriptstyle}\vcenter{\xymatrix@R=0pt@C=0pt{1-m\\0\\m}}\right]&&
\mathbf{g}_{2m+2}^{w}=\left[\def\objectstyle{\scriptstyle}\def\labelstyle{\scriptstyle}\vcenter{\xymatrix@R=0pt@C=0pt{1-m\\-1\\m}}\right]\\\label{Eq:GNegativeW}
\mathbf{g}_{-(2m+1)}^{w}=\left[\def\objectstyle{\scriptstyle}\def\labelstyle{\scriptstyle}\vcenter{\xymatrix@R=0pt@C=0pt{-m\\1\\m-1}}\right]&&
\mathbf{g}_{-2m}^{w}=\left[\def\objectstyle{\scriptstyle}\def\labelstyle{\scriptstyle}\vcenter{\xymatrix@R=0pt@C=0pt{-m\\0\\m-1}}\right]
\end{eqnarray}
For every $n\geq1$
\begin{equation}\label{eq:GvectorUnW}
\begin{array}{ccc}
\mathbf{g}_{w}^{w}=\left[\def\objectstyle{\scriptstyle}\def\labelstyle{\scriptstyle}\vcenter{\xymatrix@R=0pt@C=0pt{0\\1\\0}}\right]&
\mathbf{g}_{z}^{w}=\left[\def\objectstyle{\scriptstyle}\def\labelstyle{\scriptstyle}\vcenter{\xymatrix@R=0pt@C=0pt{0\\-1\\0}}\right]&
\mathbf{g}_{u_{n}}^w=\left[\def\objectstyle{\scriptstyle}\def\labelstyle{\scriptstyle}\vcenter{\xymatrix@R=0pt@C=0pt{-n\\ 0\\ n}}\right]
\end{array}
\end{equation}
\end{proposition}
\begin{proof}
Let $x$ be a cluster variable. The following lemma shows a formula which relates the $\mathbf{g}$--vector of $x$ in $\{x_1,x_2,x_3\}$ and in $\{x_1,w,x_3\}$. 
\begin{lemma}\label{lemma:Gw=G}
For every $b\in\BB$ the $\mathbf{g}$--vector $\mathbf{g}_b=(g_1,g_2,g_3)^t$ of $b$ in $\{x_1,x_2,x_3\}$ and the $\mathbf{g}$--vector
$\mathbf{g}_b^w=(g_1^w,g_2^w,g_3^w)^t$ of $b$ in $\{x_1,w,x_3\}$ are related by the following formula
\begin{equation}\label{eq:RecursionGCyclic}
\begin{array}{ccc}
g_1^w=g_{1}+g_{2}-\textrm{min}(g_2,0),&g_2^w=-g_{2},&g_3^w=g_{3}+\textrm{min}(g_2,0)
\end{array}
\end{equation}
\end{lemma}
\begin{proof}[Proof of lemma~\ref{lemma:Gw=G}]
Let $\PP=\textrm{Trop}(\yy_1,\yy_2,\yy_3)$ and let $\myAA_\bullet(\Sigma^{cyc})$ be the cluster algebra with principal coefficients at the seed $\Sigma^{cyc}$ defined in \eqref{Def:CyclicSeed}. Let $x$ be a cluster variable of $\myAA_\bullet(\Sigma^{cyc})$. We expand $x$ in $\{x_1,w,x_3\}$  as in both \eqref{Eq:XW} and \eqref{Eq:b=Fbgb} and we find the equality
$$
\frac{F_{x}(\frac{\hat{\yy}_{1}\hat{\yy}_{2}}{\hat{\yy}_{2}+1},\frac{1}{\hat{\yy}_{2}},\hat{\yy}_{3}(\hat{\yy}_{2}+1))}{F_{x}|_\PP(\yy_{1}\yy_{2},\frac{1}{\yy_{2}},\yy_{3})}x_1^{g_1}(\frac{x_{1}+\yy_{2}x_{3}}{w})^{g_2}x_3^{g_3}=F_x^w(\hat{\yy}_{1}, \hat{\yy}_{2}, \hat{\yy}_{2})x_1^{g_1^w}w^{g_2^w}x_3^{g_3^w}.
$$
By \eqref{eq:Fb|P} and \eqref{eq:Fw=FII} we get:
$$
\left\{
\begin{array}{cc}
 \yy_2x_1^{g_1}(\frac{x_1+\yy_2x_3}{w})^{g_2}x_3^{g_3}=\frac{\hat{\yy}_2}{1+\hat{\yy}_2}x_1^{g_1^w}w^{g_2^w}x_3^{g_3^w}&\textrm{ if }x=x_{-(2m+1)}\textrm{ or }w\\
x_1^{g_1}(\frac{x_1+\yy_2x_3}{w})^{g_2}x_3^{g_3}=(1+\hat{\yy}_2)^{g_2}x_1^{g_1^w}w^{g_2^w}x_3^{g_3^w}&\textrm{ otherwise}
\end{array}
\right.
$$
from which \eqref{eq:RecursionGCyclic} follows by using proposition~\ref{prop:GVectorsCyclic}. Let now $b=s_1^as_2^bs_3^c$ be a  cluster monomial. In proposition~\ref{prop:GVectorsCyclic} we notice that the second entry of the $\mathbf{g}$--vector of both $s_1$, $s_2$ and $s_3$ have the same sign and hence the transformation \eqref{eq:RecursionGCyclic} is linear. Then $\mathbf{g}_b^w$ is given by \eqref{eq:RecursionGCyclic}. The same argument works if $b=u_nw^k$ or $b=u_nz^k$. 
\end{proof}
The proof of proposition~\ref{prop:GVectorsCyclic} follows from proposition~\ref{Prop:GvectorsInitial} by  lemma~\ref{lemma:Gw=G}.
\end{proof}
Formula \eqref{eq:RecursionGCyclic} between $\mathbf{g}$--vectors of cluster monomials in  two adjacent clusters  was conjectured in \cite{FZIV} and proved in \cite{FuKeller} in much more generality.

%

\subsection{Proof of proposition~\ref{Prop:gbDb}}
We denote by  
$$
f:=\left(\def\objectstyle{\scriptscriptstyle}\def\labelstyle{\scriptscriptstyle}\vcenter{\xymatrix@R=0pt@C=0pt{-1 & 0 & 0\\[?] & -1 & 0\\[?] & [?] & -1}}\right)
$$
the map $f:Q\rightarrow Q$ which acts on $Q$ as follows
$$
f\cdot(a,b,c)^t:=(-a, -b+[a]_+,-c+[a]_++[b]_+)
$$
where $[b]_+:=\textrm{max}(b,0)$.  Proposition~\ref{Prop:gbDb} says that for every element $b$ of $\BB$ the corresponding $\mathbf{g}$--vector $\mathbf{g}_b$ and denominator vector $\mathbf{d}(b)$ are related by 
\begin{equation}\label{Eq:Gb=FD}
\mathbf{g}_b=f\cdot \mathbf{d}(b)
\end{equation}
We hence prove \eqref{Eq:Gb=FD}.  By the explicit formulas for the g--vectors given in proposition~\ref{Prop:GvectorsInitial} and from the explicit formulas for the denominator vectors given in \eqref{eq:D(XmPositive)},  formula \eqref{Eq:Gb=FD} holds for cluster variables and for the $u_{n}$'s. 

By corollary~\ref{Cor:SignCoherence} denominator vectors of cluster variables belonging
to the same cluster are sign--coherent. It is clear that if $v_1$ and $v_2$ are sign--coherent then $f\cdot(v_1+v_2)=f\cdot v_1+f\cdot v_2$. Moreover $f$ is injective and hence  \eqref{Eq:Gb=FD} holds for cluster monomials. By lemma~\ref{lemma:DenominatorsCyclic} denominator
vectors of the $u_{n}$ 's, $w$ and $z$ lie in the positive octant $Q_+$ in which $f$ is linear. The
claim is hence true for $u_{n}w^{k}$ and $u_{n}z^{k}$, $n,k \geq 1$.

\subsection{Proof of proposition~\ref{Prop:GVectorPar}}
By proposition~\ref{Prop:gbDb} the map $\mathbf{d}(b)\mapsto\mathbf{g}_b$ is bijective. By theorem~\ref{Prop:BijectionDenominators} the map $b\mapsto\mathbf{d}(b)$ is a bijection between $\BB$ and $Q$. Then the composition $b\mapsto\mathbf{d}(b)\mapsto\mathbf{g}_b$ is a bijection between $\BB$ and $Q$. 
 
The map $\mathbf{g}_b\mapsto \mathbf{g}_b^w$ given by \eqref{eq:RecursionGCyclic} is bijective and hence the map $b\mapsto \mathbf{g}_b^w$ is bijective.

\section{Proof of theorem~\ref{Thm:CanonicalBasis}}\label{Sec:CanonicalBasis}

Let $\PP$ be a tropical semifield. In this section we prove that the set $\BB$ of cluster monomials and of the elements $\{u_{n}w^{k},u_{n} z^{k}|n \geq 1, k \geq 0\}$ of the cluster algebra $\myAA_{\PP}$ has the following properties:
\begin{itemize}
      \item $\BB$ is a linearly independent set over $\ZZ\PP$ (section~\ref{sec:LinearIndipendence});
      \item the elements of $\BB$ are positive (section~\ref{subsec:Positivity});
      \item $\BB$ spans $\myAA_{\PP}$ over $\ZZ\PP$ (section~\ref{subsec:StraighteningRelations});
      \item the elements of $\BB$ are positive indecomposable (section~\ref{sec:PositiveIndecomposability}).
\end{itemize}
and hence $\BB$ is an atomic basis   of $\myAA_{\PP}$.

\subsection{Linear independence of $\BB$}\label{sec:LinearIndipendence}
Let $\PP$ be an arbitrary semifield and let $\myAA_{\PP}$ be the cluster algebra with initial seed $\Sigma$ given by \eqref{Eq:InitialSeed}. In view of  proposition~\ref{Prop:Homogeneity} the expansion of an element $b$ of $\BB$ in the seed $\Sigma$ has the form:
\begin{equation}\label{Eq:B=FGProof}
b=\frac{F_{b}(y_{1}\mathbf{x}^{\mathbf{h}_{1}},y_{2}\mathbf{x}^{\mathbf{h}_{2}},y_{3}\mathbf{x}^{\mathbf{h}_{3}})\mathbf{x}^{\mathbf{g}_{b}}}{F_{b}|_{\PP}(y_{1},y_{2},y_{3})}
\end{equation}
where $F_{b}$ and $\mathbf{g}_{b}$ are respectively the polynomial and the vector given in section~\ref{subsec: Homogeneity}, $\mathbf{h}_{i}$ is the $i$--th column vector of the exchange matrix $H$ of the seed $\Sigma$. Moreover the polynomial $F_{b}$ has the form: $F_{b}(y_{1},y_{2},y_{3})=1+\sum_{\mathbf{e}}\chi_{\mathbf{e}}(b)\mathbf{y}^{\mathbf{e}}$ where the sum is over a finite set of non--negative integer vectors $E(b)=\{\mathbf{e}\in\ZZ_{\geq0}^{3}\setminus\{0\}|\chi_{\mathbf{e}}(b)\neq0\}$ and the coefficients $\{\chi_{\mathbf{e}}(b):\mathbf{e}\in E(b)\}$ are positive integer numbers. The expansion of $b$ has hence the form:
$$
b=\frac{\mathbf{x}^{\mathbf{g}_{b}}+\sum_{\mathbf{e}=(e_{1},e_{2},e_{3})}\chi_{\mathbf{e}}(b)\mathbf{y}^{\mathbf{e}}\mathbf{x}^{\mathbf{g}_{b}+\sum_{i=1}^{3}e_{i}\mathbf{h}_{i}}}{1\oplus\bigoplus_{\mathbf{e}}\chi_{\mathbf{e}}(b)\mathbf{y}^{\mathbf{e}}}
$$
We say that a monomial $\mathbf{x}^{\mathbf{c}}$ is a \emph{summand} of an element $b$ of $\BB$ in the cluster $\{x_1,x_2,x_3\}$ if it appears with non zero coefficients in the expansion \eqref{Eq:B=FGProof} of $b$ in the cluster $\{x_1,x_2,x_3\}$.

We introduce in $\ZZ^{3}$ the following binary relation: given $\mathbf{s}, \mathbf{t}\in\ZZ^{3}$ we say that $\mathbf{s}\leq_{H}\mathbf{t}$ if and only if there exist non--negative integers $e_{1},e_{2},e_{3}$ such that $\mathbf{t}=\mathbf{s}+\sum_{i=1}^{3}e_{i}\mathbf{h}_{i}$. The vectors $\mathbf{h}_{1}$, $\mathbf{h}_{2}$ and $\mathbf{h}_{3}$ form a pointed cone, i.e. a non--negative linear combination of them is zero if and only if all the coefficients are zero. The relation $\leq_{H}$ is hence a partial order on $\ZZ^{3}$. The map $b\mapsto\mathbf{g}_{b}$ between $\BB$ and $\ZZ^{3}$ is injective (actually bijective) by proposition~\ref{Prop:GVectorPar} and hence the partial order $\leq_{H}$ induces a partial order on $\BB$ given by: 
\begin{equation}\label{def:OrderingOnB}
b \leq b' \Longleftrightarrow\mathbf{g}_b \leq_{H} \mathbf{g}_{b'}.
\end{equation}
In particular every finite subset $\BB'$ of $\BB$ has a minimal object $b_{0}$. Then the monomial $\mathbf{x}^{\mathbf{g}_{b_{0}}}$ is not a summand of any other element of $\BB'$. We conclude that $\BB$ is a linearly independent set over $\ZZ\PP$.
\begin{remark}
Linear independence of cluster monomials for cluster algebras with an ``acyclic'' seed was proved in \cite{CK2} and \cite{GLS08}.  
\end{remark}

\subsection{Positivity of the elements of $\BB$}\label{subsec:Positivity}  
In this section we show that the elements of the set $\BB$ defined in theorem~\ref{Thm:CanonicalBasis} are positive, i.e. their Laurent expansion in every cluster of $\myAA_{\PP}$ has coefficients in $\ZZ_{\geq0}\PP$. In view of remark~\ref{Rem:Symmetries} it is sufficient to show that they have such property only in the two clusters $\{x_{1},x_{2},x_{3}\}$ and $\{x_{1},w,x_{3}\}$. By theorem~\ref{Thm:ExplicitExpressions} the elements $w$, $z$ , $u_{n}$ and $x_{m}$ with $m,n\geq1$ have such property. Since a product of positive elements is positive, we conclude that  the cluster monomials and  the elements $\{u_nw^k,u_nz^k:\,n,k\geq1\}$ are positive.

\subsection{The set $\BB$ spans $\myAA_{\PP}$ over $\ZZ\PP$}\label{subsec:StraighteningRelations}
In this section we show the set $\BB$ defined in theorem~\ref{Thm:CanonicalBasis} spans $\myAA_{\PP}$ over $\ZZ\PP$. The strategy of the proof is the following: since $\BB$ contains cluster variables, the monomials in its elements span $\myAA_{\PP}$ over $\ZZ\PP$. It is then sufficient to express every such monomial as a $\ZZ\PP$--linear combination of elements of $\BB$. We introduce the following multi--degree on monomials in the elements of $\BB$: the generic monomial $M$ has the form $M=u_{n_{1}}^{a_{1}}\cdots u_{n_{s}}^{a_{s}}x_{m_{1}}^{b_{1}}\cdots x_{m_{t}}^{b_{t}}w^{c}z^{d}$ where $0 < n_{1}< \cdots < n_{s}$, $m_{1}< \cdots < m_{t}$ and the exponents are positive integers. We define the multi-degree
$\mu(M ) = (\mu_{1}(M), \mu_{2} (M), \mu_{3}(M)) \in\ZZ_{\geq0}^3$ by setting
\begin{equation}\label{eq:Multidegree}
\left\{
\begin{array}{l}
\mu_{1}(M):=\sum_{i=1}^{s}a_{i}+\sum_{j=1}^{t}b_{j}+c+d\\
\mu_{2}(M):=m_{t}-m_{1}\\
\mu_{3}(M):=b_{1}+b_{t}
\end{array}
\right.
\end{equation}
The lexicographic order of $\ZZ^3$ makes it into a well ordered set (i.e., every non-empty
subset of $\ZZ_{\geq0}^3$ has the smallest element). In section~\ref{SubSec:ProofSpanProperty} we  show that every such monomial can be expressed as a linear combination of monomials of  (lexicographically) smaller multi-degree. In section~\ref{SubSubSec:Straightening} we find the minimal monomials which do not belong to $\BB$ and express them as linear combinations of elements of $\BB$.  We refer to such expressions as straightening relations.

\subsubsection{Straightening relations}\label{SubSubSec:Straightening}
In this section we express the  monomials in the elements of $\BB$ which do not belong to $\BB$ and are minimal with respect to the multi--degree \eqref{eq:Multidegree} as $\ZZ\PP$--linear combinations of elements of $\BB$. We notice that  the multi--degree \eqref{eq:Multidegree} does not depend on coefficients, and hence a cluster variable $s$ and the relative principal cluster variable $S:=F_s|_\PP(y_1,y_2,y_3)s$ in the initial seed $\Sigma$  have the same multi--degree.  It is convenient to consider principal cluster variables in the initial seed $\Sigma$ since the exchange relations are simpler and are given by \eqref{ExchRelPrincipalX0}, \eqref{ExchRelPrincipal}, \eqref{ExRelWX2mExplicit}, \eqref{ExRelZX2m+1Explicit}, \eqref{ExchRel:X2m-2X2m+2Principal} and \eqref{ExchRel:X2m-1X2m+3Principal}. 

Such minimal monomials are the following:
\begin{equation}\label{eq:MinimalMonomials}
\begin{array}{cccccc}
u_{n}u_{p};&u_{n}X_{m};&X_{m}X_{m+2+n};&X_{2m}W;&X_{2m+1}Z;&ZW
\end{array}
\end{equation}
for every $n, p \geq 1$ and $m \in\ZZ$. Indeed every monomial $M$ in \eqref{eq:MinimalMonomials} satisfies $\mu_1(M)=2$  and hence they are minimal (it follows from the definition that $\mu_{1}(M)=1$ if and only if $M$ is either a cluster variable or one of the $u_{n}$'s). Moreover they are the only monomials not belonging to $\BB$ with this property.

The straightening relations for the monomials $X_{2m}W$, $X_{2m+1}Z$, $ZW$ are given respectively by \eqref{ExRelWX2mExplicit}, \eqref{ExRelZX2m+1Explicit} and \eqref{Eq:DefinitionU0u_1U2}.
Propositions~\ref{Prop:UnUp} and \ref{Prop:Straightening} give the remaining ones.

\begin{proposition}\label{Prop:UnUp}
For every $n,p \geq 1$:
\begin{equation}\label{eq:UnUp}
u_{n}u_{p}=\left\{
\begin{array}{cc} u_{n+p}+\mathbf{y}^{p\mathbf{\delta}}u_{n-p}&\text{ if }n>p\\
u_{2n}+2\mathbf{y}^{n\mathbf{\delta}}&\text{ if }n=p
\end{array}
\right.
\end{equation}
where $\delta:=(1, 1, 1)^{t}$.
\end{proposition}
\begin{proof}
We use the definition of the $u_{n}$'s given in \eqref{Eq:DefinitionUn+1}. For simplicity we assume now
that $u_{0}:=2$, so that the relation $u_{1}u_{n}=u_{n+1}+\mathbf{y}^{\delta}u_{n-1}$ holds for every $n \geq 1$ (instead
of holding only for $n \geq 2$). Moreover, with this convention, we have to prove that for every $p:1\leq p \leq n$ we have:
\begin{equation}\label{eq:UnUpReformulation}
u_{n}u_{p}=u_{n+p}+\mathbf{y}^{p\mathbf{\delta}}u_{n-p}
\end{equation}
If $n = p = 1$ then \eqref{eq:UnUpReformulation} is the definition \eqref{Eq:DefinitionU0u_1U2} of $u_{2}$; we assume $n \geq 2$ and we proceed by induction on $p \geq 1$: if $p = 1$, then \eqref{eq:UnUpReformulation} is just the definition \eqref{Eq:DefinitionUn+1} of $u_{n+1}$. We then assume $2\leq p+1\leq n$ and we get:
\begin{eqnarray}\nonumber
u_{n}u_{p+1}=u_{n}[u_{1}u_{p}-\mathbf{y}^{\delta}u_{p-1}]=&&\\\nonumber
=u_{1}[u_{n+p}+\mathbf{y}^{p\delta}u_{n-p}]-\mathbf{y}^{\delta}[u_{n+p-1}+\mathbf{y}^{(p-1)\delta}u_{n-p+1}]=&&\\\nonumber
=u_{n+1+p}+\mathbf{y}^{\delta}u_{n+p-1}+\mathbf{y}^{p\delta}[u_{n+1-p}+\mathbf{y}^{\delta}u_{n-p-1}]\\\nonumber-\mathbf{y}^{\delta}[u_{n+p-1}+\mathbf{y}^{(p-1)\delta}u_{n-p+1}]=&&\\\nonumber
=u_{n+p+1}+\mathbf{y}^{(p+1)\delta}u_{n-(p+1)}&&
\end{eqnarray}
\end{proof} 

In order to give the remaining straightening relations we need to introduce the following coefficients. We use the notation:
$$
\mathbf{y}^{\mathbf{e}}\stackrel{\text{min}}{\oplus}\mathbf{y}^{\mathbf{d}}:=y_1^{\text{min}(e_1,d_1)}y_2^{\text{min}(e_2,d_2)}y_3^{\text{min}(e_3,d_3)}
$$  
\begin{definition}\label{Def:CoefficientsStraightening}
For every $m\in\ZZ$ we define
\begin{equation}\label{def:Csi}
 \xi_{m}:=\left\{\begin{array}{cc}
   \mathbf{y}^{\mathbf{d}(x_{m+3})} = y_{1;m}&\text{ if } m \geq1\\ 
   \bfy^{\mathbf{d}(x_{m})} = y_{1;m-3}&\text{ if } m \leq 0                                     
\end{array}\right.
\end{equation}
and also
\begin{equation}\label{def:Zeta+-}
\begin{array}{cc}
\zeta_{n}^{-}(m)=\left\{\begin{array}{cc}                         
\xi_{m} \!\!\!\stackrel{\text{min}}{\oplus}\!\!\!\bfy^{n\delta}\!\!\!&\!\!\!\text{ if }m \geq 1\\
       1\!\!\!&\!\!\!\text{ if } m \leq 0
\end{array}\right.\!\!\!;&
\zeta_{n}^{+}(m)=\left\{\begin{array}{cc}                         
       1\!\!\!&\!\!\!\text{ if }m\geq1\\
\xi_{m}\!\!\!\stackrel{\text{min}}{\oplus}\!\!\! \bfy^{n\delta}\!\!\!&\!\!\!\text{ if } m \leq 0.\\
\end{array}\right.
\end{array}
\end{equation}
For every integer $k \geq 0$ we define
\begin{equation}\label{def:gamma123}
  \begin{array}{ccc}
\gamma_1(k)=y_1^{\lceil\frac{k}{2}\rceil}y_2^{\lfloor\frac{k}{2}\rfloor}y_3^{\lceil\frac{k}{2}\rceil};&\!\!\!\! \gamma_2(k)= y_1^{\lfloor\frac{k}{2}\rfloor}y_2^{\lceil\frac{k}{2}\rceil}y_3^{\lfloor\frac{k}{2}\rfloor};&\!\!\!\!  \gamma_3 (k)=\!\!\left\{\begin{array}{cc}\mathbf{y}^{\frac{k}{2}\delta}&\!\!\text{if }k\text{ is even},\\
    0&\!\!\text{if }k\text{ is odd}.\end{array}\right.
 \end{array}
\end{equation}
and we define for $i = 1, 2, 3$ the corresponding elements of $\myAA_{\PP}$ :
$$
\Gamma_i (n)=\sum_{k\geq0}(\lfloor \frac{k}{2} \rfloor + 1)\cdot\gamma_{i}(k)\cdot u_{n-k}.
$$
We also define for every $m \in \ZZ$ and $m_{1}\geq 0$:
$$
\eta_{m;m_{1}}^{-}:=\left\{
\begin{array}{cc}
   \xi_{m} \stackrel{\text{min}}{\oplus} \xi_{m+m_{1}}&\text{ if } m \leq 0 < m + m_1\\
   1&\text{ otherwise}
\end{array}\right.
$$
and
$$
\eta_{m;m_{1}}^{+}:=\left\{
\begin{array}{cc}
   1&\text{ if } m \leq 0 < m + m_1\\
   \xi_{m}  \stackrel{\text{min}}{\oplus} \xi_{m+m_{1}}&\text{ otherwise}
\end{array}\right.
$$
\end{definition}

\begin{proposition}\label{Prop:Straightening}
With the notations of definition~\ref{Def:CoefficientsStraightening} the following equalities hold.
\begin{itemize}
      \item[(i)]: For every $m \in \ZZ$ and $n \geq 1$:
\begin{equation}\label{Eq:UnXm}
 u_{n}X_{m}= \zeta_{n}^{-}(m) X_{m-2n} + \zeta_{n}^{+}(m)X_{m+2n}
\end{equation}
      \item[(ii)]: For every $m \in \ZZ$ even and $n \geq 0$:
\begin{equation}\label{eq:XmXm+2n+3Even}
X_{m} X_{m+2n+3}=\eta_{m;2n+3}^{-} X_{m+n+1}X_{m+n+2} + \eta_{m;2n+3}^{+} \Gamma_1 (n)
\end{equation}
      \item[(iii)]: For every $m \in \ZZ$ odd and $n \geq 0$:
\begin{equation}\label{eq:XmXm+2n+3Odd}
X_{m} X_{m+2n+3} = \eta_{m;2n+3}^{-} X_{m+n+1} X_{m+n+2} + \eta_{m;2n+3}^{+} \Gamma_2(n)
\end{equation}
     \item[(iv)]: For every $m \in \ZZ$ even and $n \geq 2$:
\begin{equation}\label{eq:XmXm+2nEven}
X_{m}X_{m+2n}=\eta_{m;2n}^{-}X_{m+2\lfloor\frac{n}{2}\rfloor}X_{m+2\lceil\frac{n}{2}\rceil}+\eta_{m;2n}^{+}\Gamma_3(n-2)Z
\end{equation}
     \item[(v)]: For every $m \in \ZZ$ odd and $n \geq 2$:
\begin{equation}\label{eq:XmXm+2nOdd}
X_{m}X_{m+2n}=\eta_{m;2n}^{-}X_{m+2\lfloor\frac{n}{2}\rfloor}X_{m+2\lceil\frac{n}{2}\rceil}+\eta_{m;2n}^{+}\Gamma_3(n-2)W
\end{equation}
\end{itemize}
\end{proposition}
\begin{proof}
We prove part $(i)$ by induction on $n\geq 1$.  By using relations \eqref{ExRelWX2mExplicit} and  \eqref{ExRelZX2m+1Explicit} it is easy to verify that for every $m \in \ZZ$ we have:
\begin{equation*}
u_1X_m=\left\{\begin{array}{cc}
\mathbf{y}^\delta X_{m-2}+X_{m+2}&\text{ if }m\geq2\\
y_1X_{-1}+X_3&\text{ if }m=1\\
X_{-2}+y_3X_2&\text{ if }m=0\\
X_{-3}+y_2y_3X_1&\text{ if }m=-1\\
X_{m-2}+\mathbf{y}^\delta X_{m+2}&\text{ if }m\leq-2
\end{array}\right.
\end{equation*}
and \eqref{Eq:UnXm} holds for $n=1$.
We now proceed by induction on $n \geq 1$. We use the convention that $u_{0}=2$ so that the relation $u_{n+1}=u_1u_n-\bfy^{\delta}u_{n-1}$ (given in definition~\ref{def:Un}) holds for every $n \geq 1$.
Moreover, with this convention, since $\zeta_{0}^{\pm}(m)= 1$, \eqref{Eq:UnXm} still holds for n = 0. We
have
\begin{eqnarray}\nonumber
u_{n+1}X_{m} = u_1u_n X_{m} - \bfy^{\delta}u_{n-1} X_{m} =&&\\\nonumber
u_n[\zeta_1^{-}(m)X_{m-2}+\zeta_1^{+}(m)X_{m+2} ]+&&\\\nonumber
-\bfy^{\delta} [\zeta_{n-1}^{-}(m)X_{m-2n+2}+\zeta_{n-1}^{+}(m)X_{m+2n-2}]=&&\\\nonumber
&&\\\nonumber
\zeta_1^{-}(m)[\zeta_{n}^{-}(m - 2)X_{m-2-2n}+\zeta_{n}^{+}(m-2)X_{m-2+2n}]+&&\\\nonumber
+\zeta_1^{+}(m)[\zeta_{n}^{-}(m+2)X_{m+2-2n}+\zeta_{n}^{+}(m+2)X_{m+2+2n}]+&&\\\nonumber
-\bfy^{\delta}\zeta_{n-1}^{-}(m)X_{m+2-2n}-\bfy^{\delta}\zeta_{n-1}^{+}(m)X_{m-2+2n}=&&\\\nonumber
&&\\\nonumber
X_{m\!-\!2\!-\!2n} [\zeta_1^{-}(m)\zeta_{n}^{-}(m-2)]\!+\!X_{m-2+2n}[\zeta_1^{-}(m)\zeta_{n}^{+}(m-2)\!-\!\bfy^{\delta} \zeta_{n-1}^{+}(m)]+&&\\\nonumber
+X_{m\!+\!2\!-\!2n}[\zeta_1^{+}(m)\zeta_{n}^{-}(m+2)-\bfy^{\delta}\zeta_{n-1}^{-}(m)]\!+\!X_{m\!+\!2\!+\!2n}[\zeta_1^{+} (m)\zeta_{n}^{+}(m+2)]&&
\end{eqnarray}
The claim follows by lemma~\ref{lemma2Straightening} below.
\begin{lemma}\label{lemma2Straightening}
For every $m\in\ZZ$ and $n\geq 1$ the following equalities hold:
\begin{enumerate}
    \item $\zeta_1^{-}(m)\zeta_{n}^{-}(m - 2) = \zeta_{n+1}^{-}(m)$;
    \item $\zeta_1^{-}(m)\zeta_{n}^{+}(m - 2)-\bfy^{\delta}\zeta_{n-1}^{+}(m) = 0$;
    \item $\zeta_1^{+}(m)\zeta_{n}^{-}(m + 2)-\bfy^{\delta}\zeta_{n-1}^{-}(m) = 0$;
    \item $\zeta_1^{+}(m)\zeta_{n}^{+}(m + 2)=\zeta_{n+1}^{+}(m)$.
\end{enumerate}    
\end{lemma}
The proof of lemma~\ref{lemma2Straightening} is by direct check.

We prove part $(ii)$ and $(iii)$ together. It is convenient to prove that the 
following relation holds for every $m\in\ZZ$, $n\geq 0$ and $i=i_{m}=1$ if $m$ is even and $2$ if m is
odd:
\begin{equation}\label{Eq:StraightReformulationii}
\eta_{m;2n+3}^{+}\Gamma_i (n) = X_{m} X_{m+2n+3} - \eta_{m;2n+3}^{-} X_{m+n+1}X_{m+n+2}
\end{equation}
We proceed by induction on $n \geq 0$. We first prove \eqref{Eq:StraightReformulationii} for $n = 0$. In this case $\Gamma_1 (0) = \Gamma_2 (0) = 1$. By the exchange relations \eqref{ExchRelPrincipalX0} and \eqref{ExchRelPrincipal} we know that for every $m \in\ZZ$ the following relation holds:
\begin{equation}\label{eq:StraighteningExchangeRelation}
X_{m} X_{m+3} =\left\{\begin{array}{cc}
y_{m+3} X_{m+1} X_{m+2}+1&\text{ if } m = 0,-1, -2;\\
X_{m+1} X_{m+2} + y_{1;m}&\text{ otherwise.}\end{array}\right.
\end{equation}

We need the following lemma.

\begin{lemma}\label{lemma:PropertiesCoeffStraightening}
With notations of definition~\ref{Def:CoefficientsStraightening} we have the following:
\begin{enumerate}
\item For every $m \in \ZZ$ and $k \geq 0$:
\begin{eqnarray}\label{eq:XimXim+k}
\xi_{m}\stackrel{\text{min}}{\oplus}\xi_{m+k}&=&\left\{
\begin{array}{cc}
\xi_{m}&\text{ if } m \geq 1,\\
\xi_{m+k}&\text{ if } m+k \leq 0.                                
\end{array}\right.\\\label{eq:XiYDelta}
\xi_{m} \stackrel{\text{min}}{\oplus}\bfy^{\delta}&=&\left\{
\begin{array}{cr}                                                  
\xi_{m}&\text{ if }  -1\leq m \leq 2,\\
\bfy^{\delta}&\text{ otherwise. }
\end{array}
\right.
\end{eqnarray}
If $m \geq 0$ and $n\geq 1$ we get
\begin{equation}\label{eq:Xi-mXin}
\xi_{-m}\stackrel{\text{min}}{\oplus}\xi_{n}=\left\{
\begin{array}{cc}
\xi_{-m}&\text{ if }m<n-1,\\
\bfy^{k\delta}&\text{ if }m=n-1=2k,\\
y_{2}\bfy^{k\delta}&\text{ if }m=n-1=2k+1,\\
\xi_{n}&\text{ if }m>n-1.
\end{array}
\right.
\end{equation}
\item For every $m \in \ZZ$ and $n \geq 1$ the following relation holds:
\begin{equation}\label{eq:Zeta+Zeta-}       
\begin{array}{cc}                               
 \zeta_{n}^{+}(m) = (13)\zeta_{n}^{-}(1-m);& \zeta_{n}^{-}(m)= (13)\zeta_{n}^{+}(1-m)
\end{array}
\end{equation}
where $(13)$ is the automorphism of $\PP$ that exchanges $y_1$ with $y_3$.
\item For every $n \geq 1$ and $i =1, 2, 3$ we have:
\begin{equation}\label{eq:u_1Gammai}
u_{1}\Gamma_i (n) = \Gamma_i (n + 1) + \bfy^{\delta}\Gamma_i (n-1) -\gamma_{i} (n+1)
\end{equation}                      
\end{enumerate}
\end{lemma}
\begin{proof}[Proof of lemma~\ref{lemma:PropertiesCoeffStraightening}] \eqref{eq:XimXim+k} and \eqref{eq:XiYDelta} follow directly by the definition of $\xi_{m}$ and
$\xi_{m+k}$ by using lemma~\ref{Lemma:DenominatorVectors}: indeed one can see that $\mathbf{d}(x_{m}) \leq \mathbf{d}(x_{m+k} )$ (resp. $\mathbf{d}(x_{m} ) \geq \mathbf{d}(x_{m+k}))$ if $m\geq 1$ (resp. $m+k \leq 0$). (Here $\leq$ is understood term by term). We now  compute $\xi_{-m} \stackrel{\text{min}}{\oplus}\xi_{n}:=\bfy^{\mathbf{d}(x_{-m})}\stackrel{\text{min}}{\oplus}\bfy^{\mathbf{d}(x_{n+3})}$. By remark~\ref{Rem:Symmetries},
$\mathbf{d}(x_{-m})=(13)\mathbf{d}(x_{m+4})$, where $(13)$ is the linear operator on $\ZZ^3$ that exchanges
the first entry with the third one. We now consider all the possible cases:
\begin{itemize}
      \item[If] $m + 4 < n + 3$ then $\mathbf{d}(x_{m+4})\leq \mathbf{d}(x_{n+3})$; since $m+4$ and $n+3$ are positive
integers, $\mathbf{d}(x_{m+4})$ and $\mathbf{d}(x_{-m})$ have respectively the form $(d_{3}+1,d_{2},d_{3})$ and $(d'_{3}+1,d'_{2},d'_{3})$ for some $d_{2}$, $d_{3}$ , $d'_{2}$ , $d'_{3}$ $\geq 0$; in particular $\mathbf{d}(x_{-m})=(d_3 , d_2 , d_3 +1)$. Since by hypothesis $d_3 < d'_3$ and $d_2 < d'_2$, we conclude $\mathbf{d}(x_{-m}) \leq \mathbf{d}(x_{n+3})$
so that $\bfy^{\mathbf{d}(x_{-m}})\stackrel{\text{min}}{\oplus} \bfy^{\mathbf{d}(x_{n+3})} = \bfy^{\mathbf{d}(x_{-m})}$.
        \item[If] $m+4=n+3 = 2k + 4$ for some $k \geq 0$, then by \eqref{eq:D(XmPositive)}, $\mathbf{d}(x_{m+4})=(k+1,k,k)^{t}$ so that $(13)\mathbf{d}(x_{m+4})=(k, k, k+1)^{t}$; then $\bfy^{\mathbf{d}(x_{-m})}\stackrel{\text{min}}{\oplus}\bfy^{\mathbf{d}(x_{n+3})}=\bfy^{k\delta}$.
      \item[If] $m + 4 = n + 3 = 2k + 5$ for some $k \geq 0$, then by \eqref{eq:D(XmPositive)}, $\mathbf{d}(x_{m+4} ) =(k + 1, k + 1, k)^{t}$ so that $(13)\mathbf{d}(x_{m+4} ) = (k, k + 1, k + 1)^{t}$ ; then $\bfy^{\mathbf{d}(x_{-m} )}\stackrel{\text{min}}{\oplus}\bfy^{\mathbf{d}(x_{n+3})}=y_2 \bfy^{k\delta}$.
      \item[If] $m + 4 > n + 3$ then $\mathbf{d}(x_{m+4} ) \geq \mathbf{d}(x_{n+3})$, then also $(13)\mathbf{d}(x_{m+4} ) \geq \mathbf{d}(x_{n+3})$.
      \end{itemize}
and \eqref{eq:Xi-mXin} is proved. Formula\eqref{eq:Zeta+Zeta-} follows from the definition and from remark~\ref{Rem:Symmetries}.   Formula \eqref{eq:u_1Gammai} follows by \eqref{eq:UnUp}.
\end{proof}
By part~1 of lemma~\ref{lemma:PropertiesCoeffStraightening}, it is immediate to verify that
$$
\begin{array}{cc}
\eta_{m;3}^{-} =\left\{\begin{array}{cc}y_{m+3}&\text{ if } m = 0, -1, -2\\ 1&\text{ otherwise} \end{array}\right.;&
\eta_{m;3}^{+} =\left\{\begin{array}{cc}1&\text{ if } m = 0, -1, -2\\ y_{1;m}&\text{ otherwise} \end{array}\right.
\end{array}
$$
so that \eqref{Eq:StraightReformulationii} specializes to \eqref{eq:StraighteningExchangeRelation} when $n = 0$, i.e. for every $m \in \ZZ$ the following relation holds
\begin{equation}\label{eq:ExchRelStraightening}            
 X_{m} X_{m+3} = \eta_{m;3}^{-} X_{m+1} X_{m+1} + \eta_{m;3}^{+}.
\end{equation}

We now assume $n \geq 1$. In this case, by the inductive hypothesis we have:
\begin{eqnarray}\nonumber
&\Gamma_i(n+1)=u_1\Gamma_i(n)-\bfy^{\delta}\Gamma_i (n -1)+\gamma_{i}(n+1)=&\\\nonumber
&&\\\nonumber
&\frac{u_1}{\eta_{m;2n+3}^{+}}\cdot[X_{m} X_{m+2n+3}-\eta_{m;2n+3}^{-}X_{m+n+1}X_{m+n+2}]+&\\\nonumber
&-\frac{\bfy^{\delta}}{\eta_{m;2n+1}^{+}}\cdot[X_{m}X_{m+2n+1}-\eta_{m;2n+1}^-X_{m+n} X_{m+n+1}]+ \gamma_{i}(n+1)=&\\\nonumber
&&\\\nonumber
&\frac{X_{m}}{\eta_{m;2n+3}^{+}}\cdot [\zeta_1^{-}(m+2n+3)X_{m+2n+1}+\zeta_1^{+}(m+2n+3)X_{m+2n+5}] +&\\\nonumber
&-\frac{\eta_{m;2n\!+\!3}^{-}}{\eta_{m;2n\!+\!3}^{+}} \cdot [\zeta_1^{-}(m\!+\!n\!+\!2)X_{m\!+\!n} X_{m\!+\!n\!+\!1}+\zeta_1^{+} (m\!+\!n\!+\!2)X_{m\!+\!n\!+\!1}X_{m\!+\!n\!+\!4}]+&\\\nonumber
&-\frac{\bfy^{\delta}}{\eta_{m;2n+1}^{+}} \cdot [X_{m} X_{m+2n+1} - \eta_{m;2n+1}^{-} X_{m+n}X_{m+n+1}]+\gamma_{i}(n+1)=&\\\nonumber
&X_{m} X_{m+2n+1} [\frac{\zeta_1^{-}(m+2n+3)\eta_{m;2n+1}^{+} -\bfy^{\delta}\eta_{m;2n+3}^{+}}{\eta_{m;2n+1}^+ \eta_{m;2n+3}^+}]\!+\!X_{m} X_{m+2n+5}[\frac{\zeta_{1}^{+}(m+2n+3)}{\eta_{m;2n+3}^{+}}]+&\\\nonumber
&+X_{m+n} X_{m+n+1} [\frac{\bfy^{\delta} \eta_{m;2n+1}^{-}\eta_{m;2n+3}^{+}-\eta_{m;2n+3}^{-}\zeta_1^{-}(m+n+2)\eta_{m;2n+1}^{+}}{\eta_{m;2n+1}^{+}\eta_{m;2n+3}^{+}}+&\\\nonumber                    
&-\frac{\zeta_1^{+}(m+n+2)\eta_{m;2n+3}^{-}}{\eta_{m;2n+3}^{+}}X_{m+n+1} X_{m+n+4}+ \gamma_{i} (n+1)=&
\end{eqnarray}
\begin{eqnarray}\nonumber
& X_{m} X_{m+2n+1} [\frac{\zeta_1^{-}(m+2n+3)\eta_{m;2n+1}^{+} -\bfy^{\delta}\eta_{m;2n+3}^{+}}{\eta_{m;2n+1}^{+} \eta_{m;2n+3}^{+}}+X_{m} X_{m+2n+5}[\frac{\zeta_1^{+}(m+2n+3)}{\eta_{m;2n+3}^{+}}]+&\\\nonumber
&+X_{m+n} X_{m+n+1} [\frac{\bfy^{\delta} \eta_{m;2n+1}^{-}\eta_{m;2n+3}^{+}-\eta_{m;2n+3}^{-}\zeta_1^{-}(m+n+2)\eta_{m;2n+1}^{+}}{\eta_{m;2n+1}^{+} \eta_{m;2n+3}^{+}}]+&\\\nonumber                                      
&\frac{\zeta_1^{+}(m+n+2)\eta_{m;2n+3}^{-}}{\eta_{m;2n+3}^{+}}[\eta_{m+n+1;3}^{+}X_{m+n+2} X_{m+n+3}+\eta_{m+n+1;3}^{+}]+\gamma_{i}(n+1)&
\end{eqnarray}
Lemma~\ref{lemma:Straightening3} below shows that this polynomial is equal to
$$
\frac{1}{\eta_{m;2n+5}^{+}}[X_{m} X_{m+2n+5}-\eta_{m;2n+5}^{-}X_{m+n+2} X_{m+n+3}]
$$
and we are done.

We prove $(iii)$ and $(iv)$ together. In order to do that we introduce the variable
$c=c(m)$ depending on $m\in\ZZ$ in the following way: $c$ is $w$ if $m$ is odd and $c$ is $z$ if $m$ is even. With this convention, both \eqref{eq:XmXm+2nEven} and \eqref{eq:XmXm+2nOdd} are equivalent to the following:
\begin{equation}\label{eq:StraighteningC}
c(m)\Gamma_3(n-2)=\frac{1}{\eta_{m;2n}^{+}}[X_{m} X_{m+2n} - \eta_{m;2n}^{-} X_{m+2\lfloor\frac{n}{2}\rfloor}X_{m+2\lceil\frac{n}{2}\rceil}].
\end{equation}
In order to prove \eqref{eq:StraighteningC} we proceed by induction on $n\geq 2$. We verify directly the formula for $n=2$ and $n=3$. We then assume $n\geq 4$. By using \eqref{eq:u_1Gammai} and the inductive hypothesis we get the following equality:
\begin{eqnarray}\nonumber
c(m)\Gamma_3 (n - 2) =X_{m}X_{m+2n+4}[\frac{\zeta_1^{-}(m+2n-2)\eta_{m;2n-4}^{+}-\bfy^{\delta}\eta_{m;2n-2}^{+}}{\eta_{m;2n-2}^{+}\eta_{m;2n-4}^{+}}]+&&\\\nonumber
X_{m} X_{m+2n} [\frac{\zeta_1^{+}(m+2n-2)}{\eta_{m;2n-2}^{+}}]+&&\\\nonumber
X_{m+\!2\lfloor\frac{n-2}{2}\rfloor}X_{m+\!2\lceil\frac{n-2}{2}\rceil}[\frac{\bfy^{\delta}\eta_{m;2n\!-\!4}^{-} \eta_{m;2n\!-\!2}^{+}\!-\!\eta_{m;2n\!-\!2}^{-}\eta_{m;2n\!-\!4}^{+}\zeta_1^{-}(m\!+\!2\lceil\frac{n-1}{2}\rceil)}{\eta_{m;2n\!-\!2}^{+} \eta_{m;2n\!-\!4}^{+}}]&&\\\nonumber
-\frac{\eta_{m;2n-2}^{-}\zeta_1^{+}(m+2\lceil\frac{n-1}{2}\rceil)}{\eta_{m;2n-2}^{+}}X_{m+2\lfloor\frac{n-1}{2}\rfloor}X_{m+2\lceil\frac{n-1}{2}\rceil}+ c(m)\gamma_3 (n-2)&&
\end{eqnarray}
Lemma~\ref{lemma:Straightening3} concludes the proof.
\end{proof}
\begin{lemma}\label{lemma:Straightening3}
For every $n\geq 1$, $m_1 \geq 3$ and $m \in \ZZ$ the following equalities hold in $\ZZ\PP$:
\begin{enumerate}    
    \item $\zeta_1^{-}(m+m_1+2)\eta_{m;m_1}^{+}-\bfy^{\delta}\eta_{m;m_1+2}=0$;
    \item $\zeta_1^{+}(m+m_1)=\eta_{m;m_1}^{+} /\eta_{m;m_1+2}^{+}$;
    \item $\bfy^{\delta}\eta_{m;m_1}^{-}\eta_{m;m_1 +2}^{+}-\eta_{m;m_1+2}^{-}\eta_{m;m_1}^{+}\zeta_1^{-}(m +\lceil\frac{m_1+2}{2}\rceil)=0$;
    \item $\zeta_1^{+}(m+n+2)\eta_{m;2n+3}^{-} \eta_{m+n+1;3}^{-}=\eta_{m;2n+3}^{+}\eta_{m;2n+5}^{-} /\eta_{m;2n+5}^{+}$;
    \item For $i=i(m):=1$ if $m$ is even and $i=i(m):=2$ if $m$ is odd we have for every $n\geq 1$:
$$\gamma_{i}(n+1)\eta_{m;2n+3}^{+}-\zeta_1^{+}(m+n+2)\eta_{m;2n+3}^{-}\eta_{m+n+1;3}^{+}=0;$$
    \item $\frac{\eta_{m;2n-2}^{-}\zeta_1^{+}(m+2\lceil\frac{n-1}{2}\rceil)}{\eta_{m;2n-2}^{+}}X_{m+2\lfloor\frac{n-1}{2}\rfloor}X_{m+2\lceil\frac{n-1}{2}\rceil}- c\gamma_3 (n - 2) =\\ 
    \frac{\eta_{m;2n}^{-}}{\eta_{m;2n}^{+}}X_{m+2\lfloor\frac{n}{2}\rfloor}X_{m+2\lceil\frac{n}{2}\rceil}$                                                     
\end{enumerate}
\end{lemma}                                                                                      
The proof of lemma~\ref{lemma:Straightening3} follows by direct check.

\subsubsection{Span property}\label{SubSec:ProofSpanProperty}

In this section we prove that a monomial of the form $M = u^{a_{1}}_{n_{1}}\cdots u^{a_{s}}_{n_{s}}x^{b_{1}}_{m_{1}}\cdots x^{b_{t}}_{m_{t}}w^{c}z^{d}$ in the elements of $\BB$ is a $\ZZ\PP$--linear combination of elements of $\BB$. 
We proceed by induction on the multi--degree $\mu(M)$ (defined in \eqref{eq:Multidegree}). If $\mu_1(M)=1$ then $M$ is either a cluster variable or one of the $u_n$'s.   
If $\sum_{i=1}^{s}a_{i} \geq 2$ (resp. $1$) then one can apply \eqref{eq:UnUp} (resp. \eqref{Eq:UnXm}), expressing $M$ as a linear combination of monomials with smaller value of $\mu_{1}$. So we can assume that $M=x^{b_{1}}_{m_{1}}\cdots x^{b_{t}}_{m_{t}}w^{c}z^{d}$. If both $c$ and $d$ are positive, by using the fact that $ZW=u_1+y_{1}y_{3}+y_{2}$, one obtains again a sum of two monomials with smaller value of $\mu_{1}$. So we assume that $d = 0$ (resp. $c = 0$) and that we can apply the  relation \eqref{ExRelWX2mExplicit} (resp. \eqref{ExRelZX2m+1Explicit}), i.e. some
$m_{i}$ is odd (resp. even). We again obtain a sum of two monomials having smaller
value of $\mu_{1}$ than the initial one. So we can assume that $M$ has one of the following
forms: $M_1:= (\prod_{m_{i}\text{ even }}x^{b_{i}}_{m_{i}})w^{c}$ or $M_{2}:= (\prod_{m_{i}\text{ odd }} x^{b_{i}}_{m_{i}})z^{d}$ or $M_{3}:=x_{m_{1}}^{b_{1}}\cdots x_{m_{t}}^{b_{t}}$ with
$m_{t}-m_1 \geq 3$. We apply either   \eqref{eq:XmXm+2nEven} or \eqref{eq:XmXm+2nOdd}  or \eqref{eq:ExchRelStraightening} to the product $x_{m_{1}}x_{m_{t}}$. By
inspection, in the resulting expression for both $M_1$ and $M_{2}$, all the monomials except at most one that has smaller value of $\mu_{1}$ have the same value of $\mu_{1}$. By further
inspection, for every such monomial $M'$, if $\text{min}(b_{1},b_{t})=1$ (resp. $\text{min}(b_{1},b_{t}) \geq 2$) then $\mu_{2} (M') < \mu_{2} (M)$ (resp. $\mu_{2}(M')=\mu_{2}(M)$ and $\mu_{3}(M')= \mu_{3} (M)-2$). Analogously in the resulting expression for $M_{3}$ , there is precisely one monomial $M$ with $\mu_{1} (M')=\mu_{1}(M)$, while the rest of the terms have smaller value of $\mu_{1}$. Moreover if $\text{min}(b_{1},b_{t})=1$ (resp. $\text{min}(b_{1},b_{t})\geq 2$) then $\mu_{2}(M') < \mu_{2}(M)$ (resp. $\mu_{2} (M')=\mu_{2}(M)$ and $\mu_{3} (M')=\mu_{3}(M)-2$).

\subsection{The elements of $\BB$ are positive indecomposable}\label{sec:PositiveIndecomposability}
In the previous sections we proved that the set $\BB=\{\textrm{cluster monomials}\}\cup\{u_nw^k,\,u_nz^k:\,n\geq1,\,k\geq0\}$ is a $\ZZ\PP$--basis of $\myAA_\PP$ (for every semifield $\PP$) and its elements are positive, i.e. their Laurent expansion in every cluster has coefficients in $\ZZ_{\geq0}\PP$. In this section we prove that given a positive element of $\myAA_\PP$ its expansion in $\BB$ has coefficients in $\ZZ_{\geq0}\PP$. Let hence $p$ be a positive element of $\myAA_\PP$. We express $p=\sum_{b\in \BB'} a_bb$ as a $\ZZ\PP$--linear combination of elements of a (finite) subset $\BB'\subset\BB$. 
\begin{definition}
Let $\myCC=\{s_1,s_2,s_3\}$ be a cluster of $\myAA_\PP$. A Laurent monomial $s_1^as_2^bs_3^c$ is called \emph{proper} if either $a<0$ or $b<0$ or $c<0$.
\end{definition}
\begin{lemma}\label{Lemma:IndecClMon}
Let $\myCC$ be a cluster of $\myAA_\PP$ and let $b\in\BB$ be not a cluster monomial in $\myCC$. The Laurent expansion of $b$ in $\myCC$ is a $\ZZ\PP$--linear combination of proper Laurent monomials.  
\end{lemma}
The proof of lemma~\ref{Lemma:IndecClMon} will be given in section~\ref{SubSec:ProofLemmaIndClMon}. Now suppose that a cluster monomial $b$ in some cluster $\myCC$ appears in the expansion of $p$ in $\BB$ with coefficient $a_b$. We expand $p$ in the cluster $\myCC$ and in view of lemma~\ref{Lemma:IndecClMon} the monomial $b$ appears with coefficient $a_b$ in this expansion. Since $p$ is positive we conclude that $a_b\in\ZZ_{\geq0}\PP$. It remains to deal with elements $u_n$, $u_nw^k$ and $u_nz^k$. Without lost of generality we can assume that the cluster monomials in $\BB'$ are of the form $x_m^as^bx_{m+2}^c$ with $m\geq1$ and $s=x_{m+1}$ or $w$ or $z$. We have the following lemma.
\begin{lemma}\label{Lemma:PosIndUn}
Let $b$ be an element of $\BB$ which is not divisible by cluster variables $x_m$ with $m\leq0$. 
\begin{enumerate}
\item The (proper) Laurent monomial $\frac{x_1^n}{x_3^n}$ is a summand  of $u_n$ in the cluster $\{x_1,x_2,x_3\}$ but it is not a summand of $b$ in $\{x_1,x_2,x_3\}$. Moreover the coefficient of $\frac{x_1^n}{x_3^n}$ in this expansion is an element $\mathbf{y}$ of $\PP$.
\item The (proper) Laurent monomial $\frac{x_1^nw^k}{x_3^n}$ is a summand  of $u_nw^k$ in the cluster $\{x_1,w,x_3\}$ but it is not a summand of $b$ in $\{x_1,w,x_3\}$. Moreover the coefficient of $\frac{x_1^nw^k}{x_3^n}$ in this expansion is an element $\mathbf{y}$ of $\PP$.
\item The (proper) Laurent monomial $\frac{x_0^nz^k}{x_2^n}$ is a summand of $u_nz^k$ in the cluster $\{x_0,z,x_2\}$ but it is not a summand of $b$ in $\{x_0,z,x_2\}$. Moreover the coefficient of $\frac{x_0^nz^k}{x_2^n}$ in this expansion is an element $\mathbf{y}$ of $\PP$.
\end{enumerate}  
\end{lemma}
The proof of lemma~\ref{Lemma:PosIndUn} will be given in section~\ref{SubSubSec:ProofLemmaInUn}. Now assume that $u_n$ (resp. $u_nw^k$, $u_nz^k$) appears with coefficient $a$ in the expansion of $p$ in $\BB$. We expand $p$ in the cluster $\{x_1,x_2,x_3\}$ (resp. $\{x_1,w,x_3\}$, $\{x_0,z,x_2\}$) and in view of lemma~\ref{Lemma:PosIndUn} (1) (resp. (2), (3)) we find that the Laurent monomial $\frac{x_1^n}{x_3^n}$ (resp. $\frac{x_1^nw^k}{x_3^n}$, $\frac{x_0^nz^k}{x_2^n}$) has coefficient  $a\mathbf{y}$ in this expansion. Since $p$ is positive we conclude that $a\mathbf{y}\in\ZZ_{\geq0}\PP$. Since $\mathbf{y}\in\PP$ we conclude that $a\in\ZZ_{\geq0}\PP$. 

In order to prove lemmas~\ref{Lemma:IndecClMon} and \ref{Lemma:PosIndUn} we use ``Newton polytopes'' of the elements of $\BB$ in every cluster of $\myAA_\PP$. This is the subject of the next section.  

\subsubsection{Newton polytopes of the elements of $\BB$}\label{subsec:Newton}

The Newton polytope of a Laurent polynomial $x\in\ZZ[s_{1}^{\pm1}, s_{2}^{\pm1},s_{3}^{\pm1}]$ with respect to the ordered set $\myCC=\{s_1 , s_2 , s_3\}$ is the convex hull in $Q_{\RR} = \RR\mathbf{e}_{1} \oplus\RR\mathbf{e}_{2}\oplus\RR\mathbf{e}_{3}$ of all lattice points $\mathbf{g}=(g_1 , g_2 , g_3 )^t$ such that the monomial $\mathbf{s}^{\mathbf{g}}:= s_1^{g_1}s_2^{g_2}s_3^{g_3}$ appears with a non-zero coefficient in $x$. We denote it by $\text{Newt}_{\myCC}(x)$.  In this section we find the Newton polytopes of the elements of $\BB$ in every cluster of $\myAA_\PP$. By the symmetry of the exchange relations it is sufficient to consider only the two clusters $\{x_1,x_2,x_3\}$ and $\{x_1,w,x_3\}$. Moreover in such clusters the Newton polytope of the cluster
variable $x_{-m}$ is obtained from the Newton polytope of the cluster variable $x_{m+4}$ by the automorphism $(13)$ of $Q_\RR$ that exchanges the first coordinate with the third one (see remark~\ref{Rem:Symmetries}). It is hence sufficient to consider only cluster variables $x_m$ with $m\geq 2$. Before doing that we notice the following interesting fact.
\begin{lemma}\label{rem:Gradation}
The algebra $\myAA_{\PP}$ is $\ZZ$--graded by the following grading:  $\mathrm{deg}(w)=2$, $\mathrm{deg}(x_{2m+1})=1$, $\mathrm{deg}(u_n) =0$, $\mathrm{deg}(x_{2m})=-1$, $\mathrm{deg}(z)=-2$ for $m\in\ZZ$,  $\mathrm{deg}(y)=0$ for every $y\in\PP$
\end{lemma}
\begin{proof}
The exchange relations \eqref{ExchRel:XmXm+3}--\eqref{ExchRel:X2m-1X2m+3} are homogeneous with respect to such grading. The fact that $\mathrm{deg}(u_n)=0$ follows from definition~\ref{def:Un}.
\end{proof}
The elements of $\BB$ are homogeneous with respect to the grading given in lemma~\ref{rem:Gradation}. This implies that the Newton polytopes of the elements of $\BB$ are actually polygons. Indeed let $\myCC=\{s_1,s_2,s_3\}$ be a cluster of $\myAA_\PP$, $b\in\BB$ and let $P_{b}^\myCC:=\{(e_{1},e_{2},e_{3})\in Q_{\RR}| \mathrm{deg}(s_1)e_{1}+\mathrm{deg}(s_2)e_{2}+\mathrm{deg}(s_3)e_{3}=\mathrm{deg}(b)\}$. Then $\mathrm{Newt}_{\myCC} (b) \subset P_{b}^{\myCC}$.

The following proposition gives the Newton polygons of the elements of $\BB$ in the cluster $\{x_{1},x_{2},x_{3}\}$ .
\begin{proposition}\label{Prop:Newton}
For every $m \geq 2$ and $n \geq 1$ we have:
$$
\text{Newt}_{\{x_{1},x_{2},x_{3}\}}(x_{2m+1})=\mathrm{Conv}\{{\scriptscriptstyle\left[\!\!\begin{array}{c}1-m \\0 \\m\end{array}\!\!\right],\left[\!\!\begin{array}{c}1-m \\1-m \\1\end{array}\!\!\right],\left[\!\!\begin{array}{c}0 \\1-m \\2-m\end{array}\!\!\right],\left[\!\!\begin{array}{c}m-2 \\-1\\2-m\end{array}\!\!\right]}\}
$$
\begin{equation}\label{eq:NewtX2m}
\text{Newt}_{\{x_{1},x_{2},x_{3}\}}(x_{2m})=\mathrm{Conv}\{{\scriptscriptstyle\left[\!\!\begin{array}{c}1-m \\1 \\m-1\end{array}\right],\left[\begin{array}{c}1-m \\2-m \\0\end{array}\right],\left[\begin{array}{c}-1 \\2-m \\2-m\end{array}\right],\left[\begin{array}{c}m-3 \\0\\2-m\end{array}\right]}\}
\end{equation}
\begin{equation}\label{eq:NewtUn}
\text{Newt}_{\{x_{1},x_{2},x_{3}\}}(u_{n})=\mathrm{Conv}\{{\scriptscriptstyle\left[\begin{array}{c}-n\\0 \\n\end{array}\right],\left[\begin{array}{c}-n \\-n \\0\end{array}\right],\left[\begin{array}{c}0 \\-n \\-n\end{array}\right],\left[\begin{array}{c}n \\0\\-n\end{array}\right]}\}
\end{equation}
\begin{equation}\label{eq:NewtW}
\text{Newt}_{\{x_{1},x_{2},x_{3}\}}(w)=\mathrm{Conv}\{{\scriptscriptstyle\left[\begin{array}{c}1\\-1 \\0\end{array}\right], \left[\begin{array}{c}0\\-1 \\1\end{array}\right]}\}
\end{equation}
\begin{equation}\label{eq:NewtZ}
\text{Newt}_{\{x_{1},x_{2},x_{3}\}}(z)=\mathrm{Conv}\{{\scriptscriptstyle\left[\begin{array}{c}0\\1\\-1\end{array}\right], \left[\begin{array}{c}-1\\0 \\-1\end{array}\right], \left[\begin{array}{c}-1\\1\\0\end{array}\right]}\}
\end{equation}
where $\mathrm{Conv}$ means convex hull in $Q_{\RR}$.
\end{proposition}
\begin{proof} Formulas \eqref{eq:NewtW} and \eqref{eq:NewtZ} follows from \eqref{Eq:DefinitionWZ}. By proposition~\ref{Prop:Homogeneity}, up to a factor in $\PP$ which does not modify the Newton polygon, every element $b$ of $\BB$ has the form
$$
b=\sum_{\mathbf{e}\in E(b)}\chi_{\mathbf{e}}(b)\mathbf{y}^{\mathbf{e}}\mathbf{x}^{\mathbf{g}_{b}+H\mathbf{e}}
$$
where $\chi_{\mathbf{e}}(b)$ is the coefficient of $\mathbf{y}^{\mathbf{e}}$ in the $F$--polynomial of $b$ given in proposition~\ref{Prop:FPolyInitial}, $E(b):=Conv\{e \in\ZZ^3 | \chi_{\mathbf{e}}(b)\neq 0\}$ is the convex hull in $Q_{\RR}$ of the support of $\mathbf{e}\mapsto\chi_{\mathbf{e}}(b)$, $\mathbf{g}_{b}$ is the $\mathbf{g}$--vector of $b$ given in proposition~\ref{Prop:GvectorsInitial} and $H$ is the exchange matrix given in \eqref{Eq:InitialSeed}. The affine map $N_{b} : \mathbf{e}\mapsto\mathbf{g}_b+H\mathbf{e}$
sends convex sets to convex sets. In particular if $E(b)=\mathrm{Conv}\{\mathbf{e}_{1},\cdots,\mathbf{e}_{n}\}$,
then $\text{Newt}(b)=\mathrm{Conv}\{N_{b}(\mathbf{e}_{1}), \cdots,N_{b}(\mathbf{e}_{n})\}$. The proof is based on the following lemma.
\begin{lemma}\label{Lemma:E}
For every $m\geq2$ and $n\geq1$ we have:
$$
E(x_{2m+1})=\mathrm{Conv}\{{\scriptscriptstyle
\left[\!\!\!\begin{array}{c}0 \\0 \\0\end{array}\!\!\!\right],
\left[\!\!\!\begin{array}{c}m-1 \\0 \\0\end{array}\!\!\!\right],
\left[\!\!\!\begin{array}{c}m-1 \\m-1 \\0\end{array}\!\!\!\right],
\left[\!\!\!\begin{array}{c}m-1 \\m-1 \\m-2\end{array}\!\!\!\right]}\};
$$
$$
E(x_{2m})=
\mathrm{Conv}\{{\scriptscriptstyle
\left[\!\!\!\begin{array}{c}0 \\0 \\0\end{array}\!\!\!\right],
\left[\!\!\!\begin{array}{c}m-1 \\0 \\0\end{array}\!\!\!\right],
\left[\!\!\!\begin{array}{c}m-2 \\m-2 \\0\end{array}\!\!\!\right],
\left[\!\!\!\begin{array}{c}m-1 \\m-2 \\0\end{array}\!\!\!\right],
\left[\!\!\!\begin{array}{c}m-1 \\m-2 \\m-2\end{array}\!\!\!\right],
\left[\!\!\!\begin{array}{c}m-2 \\m-2 \\m-3\end{array}\!\!\!\right]}\};
$$
$$
E(u_{n})=\mathrm{Conv}\{{\scriptscriptstyle
\left[\!\!\!\begin{array}{c}0 \\0 \\0\end{array}\!\!\!\right],
\left[\!\!\!\begin{array}{c}n \\0 \\0\end{array}\!\!\!\right],
\left[\!\!\!\begin{array}{c}n \\ n \\0\end{array}\!\!\!\right],
\left[\!\!\!\begin{array}{c}n \\ n \\ n\end{array}\!\!\!\right]}\}.
$$
\end{lemma}
\begin{proof} The proof follows from proposition~\ref{Prop:FPolyInitial} by direct check.
\end{proof}
In order to finish the proof of proposition~\ref{Prop:Newton} we apply the affine transformation $N_{b}$ to every generator of $E(b)$ given in lemma~\ref{Lemma:E}. If $b=x_{2m+1}$ or $b=u_n$ we find the desired  expression. If $b =x_{2m}$ we apply $N_{x_{2m}}$ to $E(x_{2m})$ and we get
$$
\text{Newt}(x_{2m})\!\!=\!\!\mathrm{Conv}\!\{\!{\scriptscriptstyle
\left[\!\!\!\begin{array}{c}1-m \\ 1 \\ m-1\end{array}\!\!\!\right],
\left[\!\!\!\begin{array}{c}1-m \\2-m \\0\end{array}\!\!\!\right],
\left[\!\!\!\begin{array}{c} -1 \\ 3-m \\ 3-m\end{array}\!\!\right],
\left[\!\!\!\begin{array}{c} -1 \\2-m \\ 2-m\end{array}\!\!\!\right],
\left[\!\!\!\begin{array}{c}m-3 \\0 \\ 2-m\end{array}\!\!\!\right],
\left[\!\!\!\begin{array}{c}m-4 \\0 \\3-m\end{array}\!\!\!\right]}\!\}.
$$
To conclude \eqref{eq:NewtX2m} we show that the third and the last generators are convex combinations of the others: indeed let $v_{1}, \cdots,v_{6}$ be the generators enumerated from left
to right. Then $v_{3}=\frac{1}{2m-3}v_{1}+\frac{2m-5}{2m-3}v_{4}+\frac{1}{2m-3}v_{5}$ and 
$v_{6}=\frac{1}{(m-1)(2m-4)}[(m-2)v_{1}+v_{2}+(m-1)(2m-5)v_{5}]$.
\end{proof}  
The following proposition gives the Newton polygons of the elements of $\BB$ in the cluster $\{x_1,w,x_3\}$.                                                     
\begin{proposition}\label{Prop:NewtonCyclic}
For every $m \geq 2$ and $n \geq 1$ we have:
$$
\mathrm{Newt}_{\{x_1,w,x_3\}}(x_{2m+1})=
\mathrm{Conv}\{{\scriptscriptstyle
\left[\!\!\!\begin{array}{c}1-m \\0 \\ m\end{array}\!\!\!\right],
\left[\!\!\!\begin{array}{c}m-3 \\ 1 \\ 2-m\end{array}\!\!\!\right],
\left[\!\!\!\begin{array}{c}1-m \\m-1 \\ 2-m\end{array}\!\!\!\right]}\};
$$
\begin{equation}\label{Eq:NewtWX2m+2}                                                          
\mathrm{Newt}_{\{x_1,w,x_3\}}(x_{2m+2})=
\mathrm{Conv}\{{\scriptscriptstyle
\left[\!\!\!\begin{array}{c}1-m \\-1 \\ m\end{array}\!\!\!\right],
\left[\!\!\!\begin{array}{c}m-2 \\ 0 \\ 1-m\end{array}\!\!\!\right],
\left[\!\!\!\begin{array}{c}-m \\-1 \\ m+1\end{array}\!\!\!\right],
\left[\!\!\!\begin{array}{c}-m \\ m-1 \\ 1-m\end{array}\!\!\!\right]}\};
\end{equation}
\begin{equation}\label{Eq:NewtWUCyclic}
\mathrm{Newt}_{\{x_1,w,x_3\}}(u_{n})=
\mathrm{Conv}\{{\scriptscriptstyle
\left[\!\!\!\begin{array}{c} n \\0 \\ -n\end{array}\!\!\!\right],
\left[\!\!\!\begin{array}{c}-n \\ 0 \\ n\end{array}\!\!\!\right],
\left[\!\!\!\begin{array}{c}-n \\ n \\ -n\end{array}\!\!\!\right]}\};
\end{equation}
\begin{equation}\label{Eq:NewtWZCyclic}
\mathrm{Newt}_{\{x_1,w,x_3\}}(z)=
\mathrm{Conv}\{{\scriptscriptstyle
\left[\!\!\!\begin{array}{c} 1 \\-1 \\ -1\end{array}\!\!\!\right],
\left[\!\!\!\begin{array}{c}-1\\ 0 \\ -1\end{array}\!\!\!\right],
\left[\!\!\!\begin{array}{c}-1 \\ -1 \\ 1\end{array}\!\!\!\right]}\};
\end{equation}
\begin{equation}\label{Eq:NewtWX2}
\mathrm{Newt}_{\{x_1,w,x_3\}}(x_2)=
\mathrm{Conv}\{{\scriptscriptstyle
\left[\!\!\!\begin{array}{c} 1 \\-1 \\0\end{array}\!\!\!\right],
\left[\!\!\!\begin{array}{c}0\\-1 \\1\end{array}\!\!\!\right]}\};
\end{equation}
where $\mathrm{Conv}$ means convex hull in $Q_{\RR}$.
\end{proposition}
\begin{proof} Formula \eqref{Eq:NewtWZCyclic} follows from \eqref{Eq:ZCyclic}. Formula \eqref{Eq:NewtWX2} follows from \eqref{ExchRel:WX2m}. By proposition~\ref{Prop:Homogeneity}, up to a factor in $\PP$, every element $b$ of $\BB$ has the form
$$
b=\sum_{\mathbf{e}\in E^w(b)}\chi_{\mathbf{e}}^{w}(b)\mathbf{\yy}^{\mathbf{e}}\mathbf{x}^{\mathbf{g}_{b}^{w}+H^{Cyc}\mathbf{e}}
$$
where $\chi_{\mathbf{e}}^{w}(b)$ is the coefficient of $\mathbf{\yy}^{\mathbf{e}}$ in the $F$--polynomial  $F_{b}^{w}$ of $b$ given in proposition~\ref{Prop:FpolyCyclic}, $E^{w}(b)=\mathrm{Conv}\{\mathbf{e} \in\ZZ^3 | \chi_{\mathbf{e}}^{w}(b)\neq 0\}$ is the convex hull in $Q_{\RR}$ of the support of $\mathbf{e}\mapsto\chi_{\mathbf{e}}^{w}(b)$, $\mathbf{g}_{b}^{w}$ is the $\mathbf{g}$--vector in the cluster $\{x_{1},w,x_{3}\}$ of $b$ given in proposition~\ref{prop:GVectorsCyclic} and $H^{Cyc}$ is the  exchange matrix given in \eqref{Def:CyclicSeed}. 

The affine map $N_{b}^{w}:e\mapsto\mathbf{g}_b^{w}+H^{Cyc}\mathbf{e}$ 
sends convex sets to convex sets. In particular if $E^{w}(b)=\mathrm{Conv}\{\mathbf{e}_{1},\cdots,\mathbf{e}_{n}\}$,
then $\text{Newt}_{\{x_{1},w,x_{3}\}}(b)=\mathrm{Conv}\{N_{b}^{w}(\mathbf{e}_{1}), \cdots,N_{b}^{w}(\mathbf{e}_{n})\}$. The proof is hence based on the following lemma.

\begin{lemma}\label{Lemma:EW}
For every $m\geq2$ and $n\geq1$ we have:
$$
E^w(x_{2m+1})=\mathrm{Conv}\{{\scriptscriptstyle
\left[\!\!\!\begin{array}{c}0 \\0 \\0\end{array}\!\!\!\right],
\left[\!\!\!\begin{array}{c}m-1 \\0 \\ m-2\end{array}\!\!\!\right],
\left[\!\!\!\begin{array}{c}m-1 \\ 0 \\0\end{array}\!\!\!\right]}\};
$$
$$
E^w(x_{2m+2})=
\mathrm{Conv}\{{\scriptscriptstyle
\left[\!\!\!\begin{array}{c}0 \\0 \\0\end{array}\!\!\!\right],
\left[\!\!\!\begin{array}{c}m-1 \\0 \\ m-2\end{array}\!\!\!\right],
\left[\!\!\!\begin{array}{c}m-1 \\0 \\0\end{array}\!\!\!\right],
\left[\!\!\!\begin{array}{c}0 \\1 \\0\end{array}\!\!\!\right],
\left[\!\!\!\begin{array}{c}m \\1 \\0\end{array}\!\!\!\right],
\left[\!\!\!\begin{array}{c}m\\1 \\m-1\end{array}\!\!\!\right]}\};
$$
$$
E^w(u_{n})=\mathrm{Conv}\{{\scriptscriptstyle
\left[\!\!\!\begin{array}{c}0 \\0 \\0\end{array}\!\!\!\right],
\left[\!\!\!\begin{array}{c}n \\0 \\0\end{array}\!\!\!\right],
\left[\!\!\!\begin{array}{c}n \\ 0 \\n\end{array}\!\!\!\right]}\}.
$$
\end{lemma}
\begin{proof} It follows from  proposition~\ref{Prop:FpolyCyclic} by a case by case inspection.
\end{proof}
In order to finish the proof of proposition~\ref{Prop:NewtonCyclic} we apply the affine transformation $N_{b}^w$ to every generator of $E^w(b)$ given in lemma~\ref{Lemma:EW}. If $b=x_{2m+1}$ or $b=u_n$ we find the desired  expression. If $b =x_{2m}$ we apply $N^w_{x_{2m}}$ to $E(x_{2m})$ and we get
$$
\mathrm{Newt}(x_{2m})\!\!=\!\!\mathrm{Conv}\!\{\!{\scriptscriptstyle
\left[\!\!\!\begin{array}{c}1-m \\ -1 \\ m\end{array}\!\!\!\right],
\left[\!\!\!\begin{array}{c}m-3 \\0 \\2-m\end{array}\!\!\!\right],
\left[\!\!\!\begin{array}{c} 1-m \\ m-2 \\ 2-m\end{array}\!\!\right],
\left[\!\!\!\begin{array}{c} -m \\-1 \\ m+1\end{array}\!\!\!\right],
\left[\!\!\!\begin{array}{c}-m \\ m-1 \\ 1-m\end{array}\!\!\!\right],
\left[\!\!\!\begin{array}{c}m-2 \\0 \\ 1-m\end{array}\!\!\!\right]}\!\}.
$$
To conclude \eqref{Eq:NewtWX2m+2} we show that the second and the third generators are convex combinations of the others: indeed let $v_{1}, \cdots,v_{6}$ be the generators enumerated from left
to right. Then $v_{2}=\frac{1}{2m}v_{4}+\frac{1}{2m(m-1)}v_{5}+\frac{2m-3}{2m-2}v_{6}$ and 
$v_{3}=\frac{1}{2m-1}v_{1}+\frac{2m-3}{2m-1}v_{5}+\frac{1}{2m-1}v_{6}$.
\end{proof}  

\subsubsection{Proof of lemma~\ref{Lemma:IndecClMon}}\label{SubSec:ProofLemmaIndClMon}
By the symmetry of the exchange relations we prove the lemma only for the two clusters $\{x_1,x_2,x_3\}$ and $\{x_1,w,x_3\}$. Let hence $\myCC$ be either the cluster $\{x_1,x_2,x_3\}$ or $\{x_1,w,x_3\}$. Let $b\in\BB$ be not a cluster monomial in $\myCC$. We prove that the Newton polygon of $b$ in $\myCC$ does not intersect the positive octant $Q_+$. Since $\mathrm{Newt}(s_{1}^{p}s_{2}^{q}s_{3}^{r})= p\mathrm{Newt}(s_1)+q\mathrm{Newt}(s_2)+r\mathrm{Newt}(s_3)$, it is sufficient to find a non--zero linear form $\varphi_{b}^\myCC:Q\rightarrow\mathbf{R}$, $\varphi^\myCC_b(e_1,e_2,e_3)=\alpha e_1+\beta e_2+\gamma c_3$, such that $\alpha,\beta, \gamma\geq0$ and that takes negative values on the vertexes of both $\mathrm{Newt}(s_1)$, $\mathrm{Newt}(s_2)$
and $\mathrm{Newt}(s_3)$. If $b$ is a  cluster monomial not divisible by cluster variables $x_m$ with $m\leq0$, tables~\ref{Tab:LinearFormInitial} and \ref{Tab:LinearFormCyclic}  show a linear form with the desired property (this can be checked directly by using lemmas~\ref{Prop:Newton} and  \ref{Prop:NewtonCyclic}). For $m\geq0$ $Newt_\myCC(x_{-m})$ is obtained from $Newt_\myCC(x_{m+4})$ by exchanging the first coordinate with the third one. In particular if $Newt_\myCC(x_{m+4})$ does not intersect the positive octant the same holds for $Newt_\myCC(x_{-m})$.

Let $\myCC=\{x_1,x_2,x_3\}$. If $b=u_nw^k$ or $b=u_nz^k$, $n,k>0$, table~\ref{Tab:LinearFormInitial} shows a linear form $\varphi_b$ which takes negative values on the vertexes of $\mathrm{Newt}_\myCC(b)$. We notice that both such forms satisfy $\varphi_b(\gamma)\leq0$ for every vertex $\gamma$ of $\mathrm{Newt}_{\myCC}(u_n)$. It hence remains to check that $1$ is not a summand of $u_n$ in $\{x_1,x_2,x_3\}$. This is done by direct check using theorem~\ref{Thm:ExplicitExpressions}. 

Let $\myCC=\{x_1,w,x_3\}$. We notice that the grading of a monomial in the elements of $\myCC$ is non--negative (see lemma~\ref{rem:Gradation}). Now $\mathrm{deg}(u_nz^k)=-2k$  and hence, for $k>0$, b is a sum of proper Laurent monomials in $\mathcal{C}$. If $b=u_n$ then $deg(u_n)=0$ and hence the only monomial in $\myCC$ that could appear in the expansion of $u_n$ in $\myCC$ is $1$. By using theorem~\ref{Thm:ExplicitExpressions} we check that this is not the case. Let hence  $b=u_nw^k$ with $k>0$. In view of \eqref{Eq:NewtWUCyclic} the linear form $\varphi_{b}^w(e_1,e_2,e_3):=e_1+e_3$ satisfies $\varphi_{b}^w(e_1,e_2,e_3)(\gamma)\leq0$ for every vertex $\gamma$ of $\mathrm{Newt}_{\{x_1,w,x_3\}}(u_n)$ and $\varphi_{b}^w(0,1,0)=0$. It follows that the only possible  monomial in $\myCC$ in the expansion of $b$ in $\myCC$  is $1$. But $1$ is a summand of $u_nw^k$ in $\myCC$ if and only if $(0,-k,0)^t$ is an element of $\mathrm{Newt}_{\{x_1,w,x_3\}}(u_n)$ which is not the case in view of \eqref{Eq:NewtWUCyclic}.  

\begin{table}[htdp]
\begin{center}
\begin{tabular}{|c|c|c|}
\hline
$b$&$\varphi_{b}$&\\
\hline\hline
$x_{2m+1}^{p}x_{2m+2}^{q}x_{2m+3}^{r}$&$[m(m-1),m(m-1),(m^{2}-2m+\frac{1}{2})]$&$m\geq2$\\\hline
$x_{3}^{p}x_{4}^{q}x_{5}^{r}$&$[1,1,0]$&\\\hline
$x_{2m}^{p}x_{2m+1}^{q}x_{2m+2}^{r}$&$[m(m-1),(m-\frac{1}{2})(m-\frac{3}{2}),(m^{2}-2m+\frac{1}{2})]$&$m\geq4$\\\hline
$x_{2}^{p}x_{3}^{q}x_{4}^{r}$&$[1,0,0]$&$r>0$\\\hline
$x_{4}^{p}x_{5}^{q}x_{6}^{r}$&$[6,3,2]$&\\\hline
$x_{6}^{p}x_{7}^{q}x_{8}^{r}$&$[9,5,5]$&\\\hline
$x_{2m+1}^{p}w^{q}x_{2m+3}^{r}$&$[m,2m,
(m-2)]$&$m\geq2$\\\hline
$x_{1}^{p}w^{q}x_{3}^{r}$&$[0,1,0]$&$q>0$\\\hline
$x_{3}^{p}w^{q}x_{5}^{r}$&$[0,1,0]$&$q>0$\\\hline
$x_{2m}^{p}z^{q}x_{2m+2}^{r}$&$[m(m-1),\frac{m}{4}(m-1),m(m-\frac{3}{2})]$&$m\geq2$\\\hline
$x_{2}^{p}z^{q}x_{4}^{r}$&$[1,0,0]$&\\\hline
$u_{n}w^{k}$&$[1,2,1]$&\!\!$n,k>0$\!\!\\\hline
$u_{n}z^{k}$&$[2,1,2]$&\!\!$n,k>0$\!\!\\\hline
\end{tabular}
\end{center}
\caption{Every summand $s=x_{1}^{e_{1}}x_{2}^{e_{2}}x_{3}^{e_{3}}$ of $b$ in $\{x_1,x_2,x_3\}$ satisfies $\varphi_{b}(e_{1},e_{2},e_{3})<0$. In the second column it is written the (row) vector that defines $\varphi_{b}$. In the first column $b$ is assumed to be not divisible by $x_{m}$, $m\leq0$; it hence satisfies conditions given in the third column on the right. We abbreviate $\varphi_b:=\varphi_b^{\{x_1,x_2,x_3\}}$.}
\label{Tab:LinearFormInitial}
\end{table}

\begin{table}[htdp]
\begin{center}
\begin{tabular}{|c|c|c|}
\hline
$b$&$\varphi^w_{b}$&\\
\hline\hline
$x_{2m+1}^{p}x_{2m+2}^{q}x_{2m+3}^{r}$&$[m(m-1),\frac{1}{4}(m-1),(m^{2}-2m+\frac{1}{2})]$&$m\geq2$\\\hline
$x_{1}^{p}x_{2}^{q}x_{3}^{r}$&$[0,1,0]$&$q>0$\\\hline
$x_{3}^{p}x_{4}^{q}$&$[1,1,0]$&$q>0~$\\\hline
$x_{3}^{p}x_{4}^{q}x_{5}^{r}$&$[1,0,0]$&$r>0$\\\hline
$x_{2m}^{p}x_{2m+1}^{q}x_{2m+2}^{r}$&$[m(m-1),(m-1),(m^{2}-2m+\frac{1}{2})]$&$m\geq2$\\\hline
$x_{2}^{p}x_{3}^{q}x_{4}^{r}$&$[0,1,0]$&\\\hline
$x_{2m+1}^{p}w^{q}x_{2m+3}^{r}$&$[m(m-1),0,
(m^{2}-2m+\frac{1}{2})]$&$m\geq2$\\\hline
$x_{3}^{p}w^{q}x_{5}^{r}$&$[1,0,0]$&\\\hline
$x_{2m}^{p}z^{q}x_{2m+2}^{r}$&$[1,2,1]$&$m\geq1$\\\hline
\end{tabular}
\end{center}
\caption{Every summand $s=x_{1}^{e_{1}}w^{e_{2}}x_{3}^{e_{3}}$ of  $b$ in $\{x_1,w,x_3\}$ satisfies $\varphi_{b}(e_{1},e_{2},e_{3})<0$. In the second column it is written the (row) vector that defines $\varphi_{b}$. In the first column $b$ is assumed to be not divisible by $x_m$, $m\leq0$; it hence satisfies conditions given in the third column on the right. We abbreviate $\varphi_b^w:=\varphi_b^{\{x_1,w,x_3\}}$.}
\label{Tab:LinearFormCyclic}
\end{table}

\subsubsection{Proof of lemma~\ref{Lemma:PosIndUn}}\label{SubSubSec:ProofLemmaInUn}
Part $(3)$ follows by part $(2)$ by the symmetry of the exchange relations. By the explicit formula of $u_n$ in both the clusters $\{x_1,x_2,x_3\}$ and  $\{x_1,w,x_3\}$ given in theorem~\ref{Thm:ExplicitExpressions}  we find that the Laurent monomial $x_1^n/x_3^n$  has coefficient in $\PP$ in these expansions and it is not a summand of $u_p$ for $p\neq n$. In particular the monomial $x_1^nw^k/x_3^n$ appears in the expansion of $u_nw^k$ in $\{x_1,w,x_3\}$ with coefficient in $\PP$ and it is not a summand of $u_pw^k$ for $p\neq n$.  The  monomial $x_1^n/x_3^n$  is not a summand  of $u_{p}w^q$ and $u_{p}z^q$ if $q> 0$, because $\mathrm{deg}(x_1^n/x_3^n)=0$  whereas
$\mathrm{deg}(u_pw^q)=2q$ and $\mathrm{deg}(u_pz^q)=-2q$ (see lemma~\ref{rem:Gradation}). Similarly the monomial  $x_1^nw^k/x_3^n$ is not a summand  of $u_p$, $u_{p}w^r$ and $u_{p}z^q$ if $p,r,q> 0$ and $r\neq k$, because $\mathrm{deg}(x_1^nw^k/x_3^n)=2k$ whereas
$\mathrm{deg}(u_p)=0$, $\mathrm{deg}(u_pw^r)=2r$ and $\mathrm{deg}(u_pz^q)=-2q$. Finally the monomial $x_1^n/x_3^n$ (resp. $x_1^nw^k/x_3^n$) is not a summand of a cluster monomial $b$ not divisible by $x_{m}$, $m\leq0$, in the cluster $\{x_1,x_2,x_3\}$ (resp. $\{x_1,w,x_3\}$) because after a glance at table~\ref{Tab:LinearFormInitial} (resp. table~\ref{Tab:LinearFormCyclic})  $\varphi_b(n,0,-n)\geq0$ (resp. $\varphi_b^w(n,k,-n)>0$). This concludes the proof.  

\section*{Acknowledgements}
This paper is part of my phd thesis \cite{Thesis} developed at the ``Universit\`a degli studi di Padova'' under the supervision of Professor Andrei Zelevinsky. I thank the director of the doctoral school, Professor Bruno Chiarellotto, for his support. I thank Professor Andrei Zelevinsky for the patience he had in introducing me to this subject and for many advices. I thank Alberto Tonolo and Silvana Bazzoni, for many conversations about this topic. I thank the Department of mathematics of the Northeastern University of Boston for its kind hospitality during the three semesters I spent there; in particular Sachin, Shih--Wei and Daniel for several conversations about cluster algebras.  I thank the referee of the previous version of this paper for his deep suggestions and comments. 

\bibliographystyle{plain}
\bibliography{Bibliografia10}
\end{document}